\documentclass[12pt]{article}
\usepackage{amssymb,amsmath,amsthm,epsfig}
\usepackage{latexsym, enumerate}
\usepackage{eepic}
\usepackage{epic}
\usepackage{graphicx}
\usepackage{subfigure}
\usepackage{color}
\usepackage{ifpdf}
\usepackage{dsfont}
\usepackage{multirow}
\usepackage{array}
\usepackage{indentfirst}
\usepackage{booktabs}
\usepackage{bm}
\usepackage{booktabs}
\usepackage{threeparttable}
\usepackage{algorithm}

\usepackage{tikz}
\usetikzlibrary{shapes.geometric, arrows}

\theoremstyle{plain}

\theoremstyle{definition}

\newtheorem{rem}{Remark}[section]

\tikzstyle{startstop} = [rectangle, rounded corners, minimum width=4cm, minimum height=1cm,text centered, draw=black, fill=red!30]
\tikzstyle{model} = [rectangle, rounded corners, minimum width=3cm, minimum height=1cm,text centered, draw=black, fill=purple!30]
\tikzstyle{re} = [rectangle, rounded corners, minimum width=2cm, minimum height=2cm,text centered, draw=black, fill=green!30]
\tikzstyle{io} = [trapezium, trapezium left angle=70, trapezium right angle=110, minimum width=3cm, minimum height=1cm, text centered, draw=black, fill=blue!30]
\tikzstyle{process} = [rectangle, minimum width=3cm, minimum height=1cm, text centered, draw=black, fill=orange!30]
\tikzstyle{decision} = [rectangle, minimum width=2cm, minimum height=1cm, text centered, draw=black, fill=green!30]
\tikzstyle{arrow} = [thick,->,>=stealth]
\tikzstyle{arrow1} = [thick,<-,>=stealth]

\newcommand{\bff}{\mathbf}
\newcommand{\bb}{\mathbb}
\newcommand{\call}{\mathcal}
\newcommand{\bt}{\bm{\theta}}

\begin{document}
\title{ \large\bf  A two-stage ensemble Kalman filter based on multiscale model reduction  for inverse problems in time fractional diffusion-wave equations}

\author{
Yuming Ba\thanks{College of Mathematics and Econometrics, Hunan University, Changsha 410082, China. Email:yumingb@hnu.edu.cn.}
\and
Lijian Jiang\thanks{Institute  of Mathematics, Hunan University, Changsha 410082, China. Email: ljjiang@hnu.edu.cn. Corresponding author}
\and
Na Ou\thanks{College of Mathematics and Econometrics, Hunan University, Changsha 410082, China. Email: oyoungla@hnu.edu.cn.}
}

\date{}
\maketitle

\begin{center}{\bf ABSTRACT}
\end{center}\smallskip

 Ensemble Kalman filter (EnKF) has been widely used in state estimation  and parameter estimation for the dynamic system where observational data is obtained sequentially in time.
 Very burdened  simulations for the forward problem are  needed to update a large number of  EnKF ensemble samples.  This will slow down the  analysis efficiency  of the EnKF for large-scale and
 high dimensional models.  To reduce uncertainty and accelerate posterior inference, a two-stage ensemble Kalman filter is presented to  improve the sequential  analysis of EnKF.
  It is known that the final posterior  ensemble may be concentrated in a small portion of the entire support of the initial prior ensemble.
   It will be much more efficient if we first build  a new prior by some  partial  observations, and construct  a surrogate only over the significant region of the new prior.
  To this end,  we construct a very  coarse model using generalized multiscale finite element method (GMsFEM) and generate a new prior ensemble in the first stage.
   GMsFEM provides a set of hierarchical multiscale basis functions supported in coarse blocks. This gives flexibility and adaptivity  to choosing degree of freedoms
   to construct a reduce model. In  the second stage,  we build an initial surrogate model based on the new prior by using GMsFEM and sparse generalized polynomial chaos (gPC)-based stochastic collocation methods.
  To improve the initial surrogate model,  we dynamically update the surrogate model, which is adapted to the sequential availability of data and the updated  analysis.
  The two-stage EnKF can achieve a better  estimation than  standard EnKF,  and   significantly improve the efficiency to update the ensemble analysis (posterior exploration).
  To enhance the applicability and flexibility in Bayesian inverse problems,  we  extend the two-stage EnKF to non-Gaussian models and hierarchical models.
  In the paper, we focus on the time fractional diffusion-wave models in porous media and investigate their  Bayesian inverse problems using the proposed two-stage EnKF.
  A few numerical examples  are carried out to demonstrate the performance of the two-stage EnKF method by taking account of parameter and structure inversion in permeability fields and source functions.

\smallskip
{\bf keywords:} GMsFEM, time fractional diffusion-wave equation, two-stage EnKF, inverse problems

\section{Introduction}

The  model inputs (parameters, source, domain geometry and system structure, et. al.) in many practical systems are often unknown.  We need to identify or estimate these inputs by partial and noisy  observations to
construct predictive models and calibrate the models.  This results in inverse problems.  In the paper, we consider the inverse problems in anomalous diffusion models.
The anomalous diffusion can be roughly classified into two categories:  subdiffusion $(0<\gamma<1)$ and superdiffusion $(1<\gamma<2)$. Here $\gamma$ denotes the fraction derivative with respect to time.
The anomalous diffusion equations are also called fractional diffusion-wave equations.
 Theoretical results such as  existence, uniqueness, stability and numerical error estimates are presented in \cite{lz} for some type of anomalous diffusion equations.
  The relationship between  anomalous diffusion equations and regular diffusion equations is discussed in \cite{OK}.
  The fractional order diffusion-wave equation as a typical fractional partial differential equation \cite{WW}, is a generalization of the classical diffusion and wave equation and can be used to better characterize anomalous diffusion phenomena in various  fields. The fractional diffusion-wave equations can model porous media applications,   viscoelastic mechanics, power-law phenomenon in fluid and complex network, allometric scaling laws in biology and ecology,  quantum evolution of complex systems and fractional kinetics \cite{MKS}.

In practice, the inputs and  parameters in the anomalous diffusion models are often unknown and need to be identified based on some observation data and prior information.
The problem of identifying unknown inputs in mathematical models has been intensively studied in the framework of inverse problems and various numerical methods have been developed \cite{anb, GLi, LMi}. The mathematical model of inverse problem  is featured with quantities which renders  useful simulation prediction obtained by imperfect model equations and measurements.  The inverse problem  is usually ill-posed.
Many methods \cite{kaipio2006statistical,  Han} such as regularization or penalty can be used to overcome the ill-posedness.  The unknown inputs (e.g., permeability field) of the anomalous diffusion models in porous media may have multiscale structure, complex geometry  patterns and uncertainty. This significantly increases the challenge of the inverse problems for these models.

Practical models  usually involve  uncertainty. Moreover,  the prior information for unknown parameters and observations are often characterized by random variables. Thus, it is desirable to treat the computational model and its
inverse problem in  statistical perspective.   Once of statistical approach for inverse problems is  Bayesian inference. The Bayesian approach \cite{kaipio2006statistical, stu10}  incorporates uncertainties in noisy observations and prior information,  and can derive  the posterior probability density of the parameters, which enables us to quantify the uncertainty in the parameters.
 The popular sampling methods in Bayesian inversion are
 Markov chain Monte Carlo (MCMC) \cite{Gil} method and its variants \cite{SD,JG},  which require  costly  computation to achieve convergence and explore the whole state space in  high dimension sample spaces.
 MCMC simulation has to run a long enough chain to give an accurate estimate, and  entails repeated solutions of the forward model.  This leads to great  challenge for solving the Bayesian inverse problem.

The ensemble Kalman filter (EnKF) is another Bayesian method. It can be seen as a reduced-order Kalman filter (KF) or a Monte Carlo implementation of KF \cite{Ge, XI}.
Since its introduction by Evenson in \cite{Ge},
 EnKF has been applied in many fields such as  oceanography, numerical weather prediction, hydrology and petroleum reservoir history matching.
   EnKF can be used for both data assimilation, where the goal is to estimate model states by incorporating  dynamical observation data,
 and inverse problem, where the objective is estimate unknown parameters appearing in models.
 Ensemble members of EnKF are forecasted in time by solving forward models  and updated by an approximate Kalman filter scheme.
 EnKF has the significant advantage that its inherent recursive process is adapted to the sequential availability of observation data in dynamical systems.
 Unlike MCMC, the ensemble samples are updated independently each other in EnKF and it is not necessary to propose samplers in a very tricky manner. Thus,  EnKF and other filter methods  have attracted much more attention in community of Bayesian inversion \cite{ihs13, OGE}.
For nonlinear models,
 the linearization of model can be avoided in EnKF using  ensemble covariance as an approximation of the posterior error covariance.
  The EnKF methods provide the first and second moments of random parameter, which are approximated by ensemble mean and ensemble covariance, respectively. Thus, the EnKF algorithm makes Gaussian approximation in a sequential manner. These approximation is accurate  for Gaussian prior models. Some insight analysis of EnKF for inverse problems has been made in recent years \cite{ogh15, ss17}.
  For the non-Gaussian models, a normal-score ensemble Kalman filter is proposed in \cite{HZ}, where the normal-score transformation is applied to transform unknown non-Gaussian  parameters to Gaussian and make the  parameters follow marginal Gaussian distributions.

   As a Bayesian sampling method,  EnKF   needs to compute  the forward problem repeatedly.
When the forward model is computationally intensive, such as multiscale models, a direct
application of EnKF forecast  with full order model  would be computationally prohibitive. In order to significantly improve the simulation efficiency, seeking more efficient sampling from the previous posterior and  building  surrogates of the forward models \cite{efendiev2005efficient, martin2012stochastic} are necessary to accelerate the EnKF analysis (posterior exploration).
Multiscale models can be solved efficiently and accurately by the numerical multiscale methods in a coarse grid instead of resolving all scales in very fine grid.
 As a numerical multiscale method,   Generalized Multiscale Finite Element Method (GMsFEM) \cite{YE, hou1997multiscale} can provide a reduce model with a good trade-off between accuracy and
 computation efficiency.  The main idea of GMsFEM is
to use a variational formulation to capture the impact of small-scale features on the coarse-scale by using multiscale basis functions.  The small-scale information is integrated into multiscale basis functions, which can be used repeatedly for  different source terms and boundary conditions of the model \cite{jl17}.
GMsFEM has been developed to solve multiscale models with complex multiscale structures and  its convergence  is independent of the contrastness of the multiscales \cite{YE}.

In the framework of EnKF, the output of model depends on  random parameters.
 We  use generalized polynomial chaos (gPC)-based stochastic collocation methods to propagate prior uncertainty through the forward model in a sequential manner.
 The gPC stochastic collocation methods require only a few number of uncoupled deterministic simulations with no reformulation of the governing equations
of the forward model. We assume that the model's output is a stochastic field and admits a gPC expansion. Then we select  a set of
collocation nodes and use least-squares methods to determine the coefficients of the gPC basis functions.
 To taking account of the potential sparsity of the gPC expansion, we  use $l_1$ regularization  to the least-squares problem. This allows
 using much fewer samples to construct a gPC surrogate model.  The idea has been employed in Bayesian inverse problems \cite{jo16, jo18, lm14, yan2015stochastic}.
 To accelerate computation for the the $l_1$ regularized  least-squares problem, the lagged diffusivity fixed point method is used.
The gPC surrogate is usually constructed based on a prior density \cite{jo16, lm14, yan2015stochastic} . However, the posterior is concentrated in a small portion of the entire prior support in many inference problems \cite{jo18}.
In EnKF, the prior is sequentially updated by incorporating new data information.
 Thus, it may be much more efficient to build surrogates only over the important region of the updated  prior than the initial  prior support.

In this paper, we propose a two-stage EnKF to take care of  the challenges and concerns mentioned above.  In the first stage,  we construct a coarse GMsFEM model with very few multiscale
basis functions, and build a new prior using standard EnKF based on the first partial  measurement data information in time.  The initial ensemble samples are drawn from the new prior for the second stage of EnKF.
By integrating GMsFEM and sparse gPC stochastic collocation method based on the new prior,  we build an initial surrogate model for the second stage.
 Because the ensemble samples are  updated by the new analysis of  EnKF,  this also sequentially updates the prior based on the new ensemble samples.
To improve the initial surrogate model,  we dynamically update the  surrogate model based on the updated prior in each assimilation step.
We note that the surrogates are constructed efficiently in EnKF procedure by using GMsFEM and sparse gPC stochastic collocation method.
By virtue of building  new priors,    we exclude the unimportant region of the posterior.
  We may use some other methods such as  ensemble smoother (ES) \cite{Em} to build  the new prior. In general, ES is used to estimate the parameters and states  when simulation models are typically stable functions.
   To extend the two-stage EnKF to non-Gaussian models, we integrate the proposed EnKF with normal-score transformation to broaden the applicability.
   The two-stage EnKF is also explored in hierarchical Bayesian inverse problems.  This increases  flexibility  in prior modeling for the Bayesian inference.

The structure of the paper is as follows. We begin with the general framework of EnKF for inverse problems.
In Section $\ref{GMs}$, we focus on the time fractional diffusion-wave models and  the surrogate model construction using GMsFEM and sparse gPC.
 Section $\ref{Sur}$ is devoted to presenting  the two-stage EnKF based on the surrogate model.
 In Section \ref{exam}, we present a few numerical examples to illustrate the performance of proposed EnKF  with applications of inverse problems for  time fractional diffusion-wave equations.
 Some conclusions and comments are made finally.


\section{Ensemble Kalman filter  for inverse problems}\label{EnKF}

  Let $\cal{U}$ be  a Hilbert space  and   $\mathcal{N}$  a generic  forward operator on  $\cal{U}$ for some physical system.
  We assume that the forward operator  describes the relation of parameter $\bt$, state $\bf{u}$ and source term $f$, i.e.,
\begin{equation}\label{ll}
  \mathcal{N}({\bf u};\bt)=f,
\end{equation}
where $\bf{u}\in\mathcal{U}$, $f\in\mathcal{U}^*$, the dual space  of $\mathcal{N}$. We assume that the solution of the problem has the form
\[
  \mathbf{u}=\mathcal{X}(\bt;f).
\]
 Let  $\bb{H}$ be the observation  operator mapping  the model state $\bf{u}\in \cal{U}$ to the observation space $\cal{Y}$
\begin{equation*}
  y=\bb{H}\big(\bf{u}(\bt)\big):=\bb{H} (\bt) \in\cal{Y}.
\end{equation*}
Let $\varepsilon$ be  additive Gaussian noise for observation.  Then  the observation data can be expressed by
\begin{equation}\label{me}
  {\bf d}={\bb H}(\bt)+\varepsilon.
\end{equation}
In the paper, we assume that  $\varepsilon$ is  independent of $\bt$. In practical setting,  the observation data is in a finite dimensional space and can be expressed by
 \[
 {\bf d}={\bb H}(\bt)+\varepsilon\in {\mathbb R}^{n_d},
 \]
 where $n_d$ is the dimension of observations.

\subsection{Bayesian inference using EnKF}

  The EnKF was introduced by Evensen \cite{Ge} as  a powerful  data assimilation method.
  Kalman filter is used for sequential update for states in  linear dynamical systems and Gaussian distribution. It provides  the mean and covariance information of the posterior distribution. When the prior is Gaussian, the filter gets the posterior Gaussian distribution from the joint Gaussian observation and the parameter. But for nonlinear dynamical system, EnKF has been widely used for data assimilation. In the paper,
  we use EnKF for Bayesian  inverse problems. This is a particular application of EnKF in recent years \cite{OGE, AS}.

 Given some observation data,  we want to estimate the parameter $\bt$.  In  Bayesian context, both $\bt$ and $\bf d$ are random variables. Thus Bayes rule gives  the posterior probability density for $\bt$ by
\begin{equation*}
p({\bt}|{\bf d}) \propto p({\bf d}|{\bt})p({\bt}),
\end{equation*}
where $p(\bt)$ is the prior distribution  before the data is observed. The data enter the Bayesian inference through the likelihood function $p(\bf d|\bt)$.

 If the information of observations is incomplete,  the covariance of the noise observation may be unknown.  For this situation, we  need to estimate  the covariance of observation noise.
 Let $\varepsilon$  be independent and identically distributed (i.i.d.) Gaussian random vector with zero mean  and variance  $\sigma^2$, i.e.,
\[
\varepsilon\sim N(0,\sigma^2\bm{I}),
\]
where ${\bm I}$ is the $n_d\times n_d$ identity matrix. Thus the likelihood function $p(\bf d|\bt)$ obeys the Gaussian distribution.
Let $\|\cdot\|$ be the Euclidean norm and $\|\cdot\|_{\bm P}=\|{\bm P}^{-\frac{1}{2}}\cdot\|$ the weighted norm, where ${\bm P}$ is the prior's covariance matrix.
 If the prior $p(\bt)$ is also Gaussian distribution, then
\begin{equation}\label{post}
p({\bt}|{\bf d}) \propto \exp\bigg(-\frac{\|{\bf d}-{\bb H}({\bt})\|^2}{2\sigma^2}-\frac{\|{\bt}-{\bt}^b\|_{\bm P}^2}{2}\bigg),
\end{equation}
where ${\bt}^b$ is the mean of  prior (background information).   When $\sigma^2$ is unknown, $\sigma^2$ is a hyperparameter in the hierarchical Bayesian model. Then the
 corresponding posterior
\begin{equation*}\label{hypost}
p({\bt}, \sigma^2|{\bf d}) \propto p({\bf d}|{\bt}, \sigma^2)p({\bt})p(\sigma^2).
\end{equation*}
The marginal posterior of $\sigma^2$ is
\[
p(\sigma^2|{\bt}, {\bf d}) \propto p({\bf d}, {\bt}|\sigma^2)p(\sigma^2).
\]
Because the likelihood
\begin{equation}\label{lik}
p({\bf d}, {\bt}|\sigma^2)=\frac{1}{(2\pi\sigma^2)^{\frac{n_d}{2}}}\exp\bigg(-\frac{\|{\bf d}-\bb{H}(\bt)\|^2}{2\sigma^2}\bigg)
\end{equation}
belongs to the inverse-gamma family, the conjugate prior $p(\sigma^2)$ can be the inverse-gamma distribution
\begin{equation}\label{hyper}
p(\sigma^2) \propto (\sigma^2)^{-(\alpha+1)}e^{\beta/{\sigma^2}}.
\end{equation}
From $(\ref{lik})$ and $(\ref{hyper})$, we get
\begin{equation}\label{albe}
 \sigma^2|{\bt}, {\bf d}\sim \text{Inv-gamma}\bigg(\alpha+\frac{n_d}{2}, \beta+\frac{\|{\bf d}-\bb{H}(\bt)\|^2}{2}\bigg).
\end{equation}
As in \cite{AG}, we choose two numbers   $n_s$  ($n_s$ often small and between $0.01$ and $1$ ) and
\[
\sigma_s^2=\frac{\|{\bf d}-\bb{H}(\bt)\|^2}{n_d-n_p},
\]
where $\bt\in {\bb R}^{n_p}$, such that $\alpha=n_s/2$ and $\beta=\sigma_s^2\alpha$.
  Once the posterior distribution of $\bt$ is inferred, we can extract the posterior mean or the maximum a posteriori (MAP) of the unknown parameter $\bt$.
  We note that the MAP estimate is equivalent to the solution of a regularized  optimization  problem.  In fact,  maximizing  the right hand side of $(\ref{post})$  is equivalent to the minimization problem
\begin{equation}\label{cos}
\min_{\bt\in{\bb R}^{n_p}}\bigg(\frac{1}{2}\big({\bf d}-{\bb H}(\bt)\big)^T{\bm R}^{-1}\big({\bf d}-{\bb H}(\bt)\big)+\frac{1}{2}({\bt}-{\bt}^b)^T{\bm P}^{-1}({\bt}-{\bt}^b)\bigg),
\end{equation}
where ${\bm R}=\sigma^2{\bm I}$. The hyperparameter $\sigma^2$ can be drawn from Inv-gamma distribution of (\ref{albe}).
 When the observation operator $\bb H$ is linear,  Kalman filter (KF)  method can be derived from ($\ref{cos}$) by completing the squares on the variable $\bt$ and gives the following analysis
\begin{equation}\label{KF}
\left\{
\begin{aligned}
\bt^a&=\bt^b+{\bm K}({\bf d}-{\bb H}(\bt^b)),\\
{\bm P}^a&=({\bm I}-{\bm K\bm H}){\bm P},
\end{aligned}
\right.
\end{equation}
where $\bm{H}=\bb{H}$ and  the Kalman gain  $\bm K$ is  given by
\[
{\bm K}={\bm P}{\bm H}^T({\bm H}{\bm P}{\bm H}^T+{\bm R})^{-1}.
\]

Let $k$ be an artificial time for the data assimilation  in a dynamic system.  We denote unknown parameter, hyperparameter and observation  data as $\bt_k$ $\in \mathbb{R}^{n_p}$, $\sigma_k^2\in \mathbb{R}$ and ${\bf d}_k\in\mathbb{R}^{n_d}$ at artificial time step $k$, respectively.  Then we define the artificial discrete  dynamic system
\[
 \left\{
 \begin{aligned}
  \bt_k &=\bt_{k-1}, \\
  {\bm d}_k &={\bb H}(\bt_k)+\varepsilon_k,\\
  \sigma_k^2 &\sim \text{Inv-gamma}(\alpha+\frac{n_d}{2}, \beta+\frac{\|{\bf d}_k-{\bb H}(\bt_{k})\|^2}{2}).
  \end{aligned}
\right.
\]
In the framework of EnKF, an estimate for $\bt$ is updated in each data assimilation step. The sequential update needs  forecast steps and analysis steps, which transport information of the current time to the next observation time in the forecast step. At the time step $k$,   we denote the forecast by $\bt^{f}_{k}$, the analysis by $\bt^{a}_{k}$, forecast error covariance matrix by $\bm P_k^f$  and analysis error covariance matrix by $\bm P_k^a$. Then we have
 \begin{equation*}
 \left\{
 \begin{aligned}
  \bt^f_k &=\bt^a_{k-1}, \\
  \sigma_k^2 &\sim \text{Inv-gamma}(\alpha+\frac{n_d}{2}, \beta+\frac{\|{\bf d}_k-{\bb H}(\bt_k^f)\|^2}{2}),\\
  {\bm P}^f_k &={\bm P}^a_{k-1}.
\end{aligned}
\right.
\end{equation*}
As in (\ref{KF}),  the posterior is the weighted sum of observations and forecast in the analysis step $k$, i.e.,
\begin{equation*}
 \left\{
 \begin{aligned}
  \bt^a_k &=\bt^f_{k}+{\bm K}_k({\bm d}_k-\bb{H}(\bt^f_k)), \\
  {\bm P}^a_k &=({\bm I}-{\bm K}_k{\bm H}){\bm P}^f_k,
\end{aligned}
\right.
\end{equation*}
where  Kalman gain
\[
\bm{K}_k=\bm{P}^f_k\bm{H}^{\rm T}(\bm{HP}^f_k\bm{H}^{\rm T}+\bm{R}_k)^{-1},\quad  \bm{R}_k=\sigma_k^2{\bm I}.
\]
 Here  $\bm{H}$ is the Jacobian matrix of $\bb{H}$ in Extended Kalman Filter (EKF) when $\bb H$ is nonlinear.

 It is well known that the KF and EKF are numerically scarcely affordable and the storage of a few state vectors is impossible in the high-dimensional systems. To overcome the difficulty,
    EnKF is desirable for  nonlinear data assimilation problems in high-dimensional space. The advantage of EnKF is that we apply a useful approximation to the Kalman filter to avoid propagating the first and second order statistical moments. To this end,  Monte Carlo method is  used to propagate an ensemble of realizations from the prior distribution. In  EnKF, we just update the propagating  ensemble and the Kalman gain matrix is approximated by
\[
\bm{K}_k=\text{Cov}({\Theta}^f_k,{Z}_k)\text{Cov}({Z}_k,{Z}_k)^{-1},
\]
where ${\Theta}^f_k$ is the forecast ensemble and ${Z}_k$ is the ensemble of simulated observations. Thus, the forecast error covariance matrix and analysis error covariance matrix are not necessary to compute. The true mean and covariance are approximated by ensemble mean and ensemble covariance, respectively. In the paper, we make use of the stochastic analysis ensemble generation method, where the simulated observations are perturbed by simulated observation error $\varepsilon^f$. The  $\varepsilon^f$ is independent of  $\varepsilon$.

 Let $\{\bm{d}_1, \bm{d}_2, \cdots, \bm{d}_k, \cdots, \bm{d}_I\}$ be a time series of observations  and ${\bf d}_k\in\mathbb{R}^{n_d}$ ($k=1,\cdots I$). We assume that the prior distribution  of $\bt$ is $\bm{\mu}_{0}$ (Gaussian), and ${\bm\varepsilon}^f_k\sim{\emph N}(0,\sigma_k^2{\bm I})$, where $\sigma_k^2$ is unknown.  We initially pick  $M$  ensemble members for EnKF,  through which we obtain the  analysis ensemble ${\Theta}^{a}_I$. Furthermore, the mean and covariance of ${\Theta}^{a}_I$ can be used to estimate the unknown parameter.
 The pseudo-code of EnKF algorithm is presented in Algorithm \ref{Seq}.
\begin{algorithm}
\caption{Sequential  EnKF algorithm with unknown $\sigma^2$}
  $\textbf{Input}:$ number of ensemble members $M$, initial ensemble members $\{\bt_{1,0},\cdots, \bt_{M,0}\}$ drawn from prior $\bm{\mu}_{0}$, the number of data assimilation steps $I$, observations $\{\bm{d}_1, \bm{d}_2, \cdots, \bm{d}_k, \cdots, \bm{d}_I\}$.\\
  $\textbf{Output}:$ ${\Theta}^{a}_{I}$.\\
  $\bf1$. ${\Theta}^a_0=(\bt_{1,0},\cdots,\bt_{M,0})$\\
  $\bf2$. \textbf{for} $k=1:I$\\
  $~~~~~$ $(1)$. $\emph{Forecast/predictor:}$ Generate ensemble of $\bm{z}_k$ by\\
  $~~~~~~~~~~~$ ${\Theta}^{f}_{k}={\Theta}^{a}_{k-1}$\\
  $~~~~~~~~~~~$ \textbf{for} $j=1:M$\\
    $~~~~~$ $~~~~~~~~$ $\bt^f_{j, k}={\Theta}^{f}_{k}(:, j)$, \quad ${Z}(:, j)=\bb{H}(\bt^f_{j, k})$, \rule{0pt}{0.8cm}
    $\sigma_s^2=\frac{\|{\bf d}-\bb{H}(\bt^f_{j, k})\|^2}{n_d-n_p}$,\quad $\beta=\sigma_s^2n_s$,\\
    $~~~~~$ $~~~~~~~~$ $S(j)\sim \text{Inv-gamma}\bigg(\alpha+\frac{n_d}{2}, \beta+\frac{\|{\bf d}_k-{\bb H}(\bt^f_{j, k})\|^2}{2}\bigg)$.\\
  $~~~~~~~~~~$ \textbf{end for}  \\
  $~~~~~~~~~~$ $\sigma_k^2=\frac{\sum_{j=1}^M S(j)}{M}$, \quad $E(:, j)={\bm\varepsilon}^f_k\sim{\emph N}(0,\sigma_k^2{\bm I}) \quad(j=1,\cdots,M)$, \quad $Z_k=Z+E$.\\
  $~~~~~$ $(2)$. $\emph{Analysis/corrector:}$ Update the previous ensemble ${\Theta}^a_{k-1}=(\bt^a_{1,k-1},\cdots,\bt^a_{M,k-1})$ by\\
    $~~~~~$ $~~~~$ ${\Theta}^{a}_{k}={\Theta}^{f}_{k}+\bm{K}_k({D}_k-{Z}_k)$, \rule{0pt}{0.8cm}\\
    $~~~~~~~~~~$ where ${D}_k=[\bm{d}_k, \cdots, \bm{d}_k]$ $\in \mathbb{R}^{n_d \times M}$ and $\bm{K}_k=\text{Cov}({\Theta}^f_k,{Z}_k)\text{Cov}({Z}_k,{Z}_k)^{-1}$. \rule{0pt}{0.8cm}\\
  \textbf{end for}\\
  \label{Seq}
 \end{algorithm}
\begin{rem}\label{EES}
 We  denote $(\bm{d}_1^T, \bm{d}_2^T, \cdots, \bm{d}_k^T, \cdots, \bm{d}_I^T) \in \mathbb{R}^{(I\cdot n_d)\times1}$  by $\bm d^{1:I}$.  We replace ${\bb H}$ by ${\bb H}^{1:I}$  in ($\ref{me}$).
  We can use ensemble smoother (ES) to do  a single global update. Then the analysis in ES is
 \begin{eqnarray*}
         {\Theta}^{a}={\Theta}^{f}+\bm{K}({D}-{Z}),
\end{eqnarray*}
where $D=[{\bm d}^{1:I}, \cdots, {\bm d}^{1:I}]$ $\in \mathbb{R}^{(I\cdot n_d) \times M}$.
\end{rem}

\subsection{EnKF for non-Gaussian model using  normal-score transformation}

In general,  EnKF is a Gaussian approximation for the estimated parameter  because it  reproduces the mean and covariance. If the target distribution is Gaussian  and unimodal, EnKF inherently gives an accurate  estimation. However, if the the target distribution is non-Gaussian  or  multimodal, the approximation may  not capture the properties of target distribution. In this situation, we can invoke the normal-score  transformation, which maps non-Gaussian into Gaussian and is invertible \cite{Pres}. We perform the normal-score  transformation  after each forecast step.  Let  ${\bb F}_k$ be  the normal-score operator at the assimilation  step $k$ and satisfy
\[
{\bm q}_{k}={\bb F}_k(\bt_{k})\sim N(0, \bm I).
\]
Because  the support of cumulative distribution function (CDF) is  $[0, 1]$, the transformation can be fulfilled by   CDF.
The normal-score transformation  renders the Gaussian random variables one by one, and the  multivariate properties of  parameter vector are also changed but  not necessary to be multi-Gaussian \cite{HZ}.

We want to incorporate the  normal-score  transformation  into EnKF for non-Gaussian cases.   Let ${\bm q}^f_{k}$ and ${\bm q}^a_{k}$ be the forecast and analysis after the  normal-score transformation. Then the forecast step of EnKF implies that
\[
 \left\{
 \begin{aligned}
  \bt^f_k &=\bt^a_{k-1}, \\
  {\bm q}^a_{k-1} &={\bb F}_{k-1}(\bt^a_{k-1}), \\
  {\bm q}^f_{k} &={\bm q}^a_{k-1}. \\
\end{aligned}
\right.
\]
The analysis step is followed by
\[
 \left\{
 \begin{aligned}
  {\bm q}^a_k &={\bm q}^f_{k-1}+{\bm K}_k({\bm d}_k-\bb{H}(\bt^f_k)), \\
  \bt^a_k &={\bb F}^{-1}_k({\bm q}^a_k),
\end{aligned}
\right.
\]
where the Kalman gain matrix $\bm{K}_k$ is approximated by
\[
\bm{K}_k=\text{Cov}({\Xi}^f_k,{Z}_k)\text{Cov}({Z}_k,{Z}_k)^{-1}.
\]
 Here ${\Xi}^f_k$ is the forecast ensemble after the transformation and ${Z}_k$ is the ensemble of simulated observations.
 Let $\bm{\mu}_{0}$ be a non-Gaussian distribution and ${\bm\varepsilon}^f_k\sim{\emph N}(0,\sigma^2{\bm I})$. We describe the normal-score EnKF (NS-EnKF) in Algorithm \ref{non}.
\begin{algorithm}
\caption{NS-EnKF  algorithm}
  \textbf{Input}: number of ensemble members $M$, initial ensemble members $\{\bt_{1,0},\cdots, \bt_{M,0}\}$ drawn from prior $\bm{\mu}_{0}$, the number of  data assimilation steps $I$, observations $\{\bm{d}_1, \bm{d}_2, \cdots, \bm{d}_k, \cdots, \bm{d}_I\}$.\\
  \textbf{Output}: ${\Theta}^{a}_{I}$.\\
 $\bf1$. ${\Theta}^a_0=(\bt_{1,0},\cdots,\bt_{M,0})$,   $\Xi^a_{0}={\bb F}_0({\Theta}^{a}_{0})$\\
  $\bf2$.\textbf{for} $k=1:I$\\
  $~~~~~$ $(1)$. $\emph{Forecast/predictor:}$ Generate ensemble of $\bm{z}_k$ by\\
  $~~~~~~~~~~~$ ${\Theta}^{f}_{k}={\Theta}^{a}_{k-1}$,\quad  $\Xi^a_{k-1}={\bb F}_{k-1}({\Theta}^{a}_{k-1})$, \quad $\Xi^f_{k}={\Xi}^{a}_{k-1}$.\\
  $~~~~~~~~~~~$ \textbf{for} $j=1:M$\\
    $~~~~~$ $~~~~~$ $~~~$ $\bt^f_{j, k}={\Theta}^{f}_{k}(:, j)$, \quad ${\bm q}^f_k=\Xi^f_k(:,j)$,\quad ${Z}(:, j)=\bb{H}(\bt^f_{j, k})$. \rule{0pt}{0.8cm}\\
  $~~~~~~~~~~~$ \textbf{end for}  \\
  $~~~~~~~~~~$ $E(:, j)={\bm\varepsilon}^f_k\sim{\emph N}(0,\sigma^2{\bm I}) \quad(j=1,\cdots,M)$, \quad $Z_k=Z+E$.\\
  $~~~~~$ $(2)$. $\emph{Analysis/corrector:}$ Update the previous ensemble ${\Theta}^a_{k-1}$ and  $\Xi^a_{k-1}$ by\\
    $~~~~~$ $~~~~~$ ${\Xi}^{a}_{k}={\Xi}^{f}_{k}+\bm{K}_k({D}_k-{Z}_k)$, \quad ${\Theta}^a_k={\bb F}^{-1}_k(\Xi^a_k)$,\rule{0pt}{0.8cm}\\
    $~~~~~~~~~$ where ${D}_k=[\bm{d}_k, \cdots, \bm{d}_k]$ $\in \mathbb{R}^{n_d \times M}$ and $\bm{K}_k=\text{Cov}({\Xi}^f_k,{Z}_k)\text{Cov}({Z}_k,{Z}_k)^{-1}$. \rule{0pt}{0.8cm}\\
  \textbf{end for}\\
  \label{non}
\end{algorithm}

\section{Surrogate model construction using GMsFEM and sparse gPC}\label{GMs}

 For the EnKF methods presented in Algorithm \ref{Seq} and  Algorithm \ref{non},  we need to repeatedly compute the forward model for all ensemble members.  This computation is very expensive when the  forward model is a complex PDE model and the number of ensemble members is large.  In order to significantly  accelerate the forward model computation, we construct a surrogate model for the forward model using model reduction methods.

 The goal is to approximate a large-scale problem in  a low dimensional space. To this end, the key idea is to choose a set of appropriate basis functions,
 which can span a good approximation space for the solution. If equation   ($\ref{ll}$) is linear with respect to ${\bf u}$,
  we can derive an algebraic system for ($\ref{ll}$) as follows by applying  suitable discretization method
\begin{equation}\label{origin}
{\bf K}({\bt}){\bf u}={\bf f},
\end{equation}
where ${\bf u}\in \bb{R}^{N_h}$ is the numerical solution vector and $\bff{f}\in \bb{R}^{N_h}$ the source vector. The $N_h$ is the number of spatial degree of freedoms and is usually very large if we
straightforwardly solve the equation in fine grid.  We can use a model reduction method and reduce the number of basis functions to improve the efficiency.
Then we can get a reduced  algebraic system for ($\ref{ll}$),
\[
{\bf K}_r({\bt}){\bf u}_r={\bf f}_r.
\]
Let  $R\in{\bb R}^{N_h\times M_v}$ ($M_v\ll N_h$) be the matrix comprised of the $M_v$ reduced basis functions. Then a projection reduce method implies
\[
{\bf K}_r({\bt})=R^T{{\bf K}}(\bt)R,\quad \quad {\bf f}_r=R^T{\bf f}.
\]

In order to accelerate evaluations of the posterior density for each updated parameter ensemble, we use stochastic response surface methods to construct surrogate. The solution ${\bf u_r(\bt)}$ of the reduced model can be expressed by stochastic basis functions such as polynomial chaos \cite{xk02}, radial basis functions \cite{ASV}, and wavelet basis functions \cite{SMA}.
The surrogate model is constructed through the stochastic collocation method by solving  a  $l_1$ penalized least-square problem.
We use the lagged diffusivity fixed point method for the $l_1$ optimization problem and get a sparse representation for ${\bf u_r(\bt)}$ using fewer samples.

\subsection{GMsFEM}

In the paper, we consider the following time fractional PDE model
\begin{equation}\label{examplelc}
  ^{c}D_t^{\gamma} {u}-\nabla(k(x)\nabla u)=f(x,t),\quad \quad x\in\Omega, \quad t\in (0, T]
\end{equation}
 where $\gamma\in(0, 1) \cup(1, 2)$  and is the fractional order of the derivative with respect to  time. Here we consider  the Caputo fractional derivative defined by
\begin{equation}\label{caputo}
  ^{c}D_t^{\gamma} {u}=\frac{1}{\Gamma(m-\gamma)}\int_0^{t}(t-\tau)^{m-\gamma-1}\frac{\partial^{m} u(x, \tau)}{\partial \tau^{m}}d\tau,\qquad m-1<\gamma<m ,
\end{equation}
where $\Gamma(\cdot)$ is the Gamma function and $m$ is a positive integer. The equation (\ref{examplelc}) has a close form subject to a suitable boundary condition and initial condition.

\begin{figure}[htbp]
  \centering
  \includegraphics[width=3.3in, height=2in]{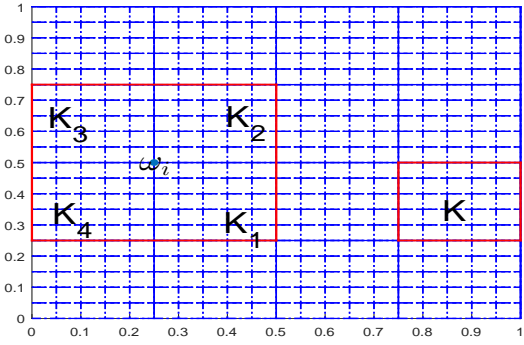}
  \caption{Illustration of a coarse neighborhood and a coarse block}\label{coarse-cell}
 \end{figure}

In the model equation,  $k(x)$ usually refers to a permeability field in porous media applications.  The permeability field has heterogeneous and multiscale structure inherently and
results in a multiscale model.  We will use general multiscale finite element method (GMsFEM) to reduce the model and get a coarse GMsFEM model. This
can achieve a good trade-off between  efficiency and  accuracy for  simulating the forward model.
 We will apply   GMsFEM presented in  \cite{YE}  to the time-fractional diffusion-wave equation (\ref{examplelc}).
For GMsFEM, we need to pre-compute a set of multiscale basis functions. To this end, the first step is to construct a snapshot space $V^{\omega_i}_{\text{snap}}$ for multiscale basis by solving  local eigenvalue problem on each coarse block $\omega_i$,
\begin{eqnarray}
\label{Mbasis}
\begin{cases}
& -\text{div}(k(x, \bt_j)\nabla \varphi_{l,j}) =\lambda_{l,j} k(x, \bt_j) \varphi_{l,j}\ \ \text{in}\ \omega_i,\\
& k(x, \bt_j)\nabla \varphi_{l,j}\cdot \vec{n} = 0\ \ \text{on}\ \partial \omega_i,
\end{cases}
\end{eqnarray}
where the samplers $\{\bt_j\}^{N_\theta}_j$ are drawn from the prior distribution of $\bt$.  By a finite element method discretization on underlying fine grid, the local eigenvalue problem can be formulated as an algebraic system,
\[
A(\bt_j)\varphi_{l,j}=\lambda_{l,j} S(\bt_j)\varphi_{l,j},
\]
where
\[
[A(\bt_j)]_{mn}=\int_{\omega_i} k(x, \bt_j)\nabla v_n\nabla v_m, \quad \quad
[S(\bt_j)]_{mn}=\int_{\omega_i} k(x, \bt_j) v_nv_m,
\]
and $v_n$ are the basis functions in fine grid.  We take the first $M^i_{\text{snap}}$ eigenfunctions corresponding to the dominant eigenvalues for each coarse neighborhood $\omega_i$ (see Figure \ref{coarse-cell}),
 $i=1,2,\cdots,N_H$,  where $N_H$ is the number of coarse nodes. Hence we  construct the space of snapshots by
\[
V^{\omega_i}_{\text{snap}}=\text{span}\{\varphi_{l,j}, 1 \leq j \leq N_\theta, 1 \leq l \leq M^i_{\text{snap}}\}.
\]
The snapshot functions can be stacked into a matrix as
\[
R_\text{snap}=[\varphi_1,\cdots,\varphi_{M_{\text{snap}}}],
\]
where $M_{\text{snap}}=N_{\theta}\times M^i_{\text{snap}}$ denotes the total number of snapshots used in the construction.
The second step  is to solve the following  local problems in the snapshot space
  \begin{eqnarray}\label{MS-basis}
   \begin{cases}
   & -\text{div}(k(x, \bar{\bt})\nabla \psi^i_k) =\lambda_k k(x, \bar{\bt}) \psi^i_k\ \ \text{in}\ \omega_i,\\
   & k(x, \bar{\bt})\nabla \psi^i_k\cdot \vec{n} = 0\ \ \text{on}\ \partial \omega_i,
   \end{cases}
  \end{eqnarray}
where $\bar{\bt}=\frac{1}{N_{\theta}}\sum_{j=1}^{N_{\theta}}\bt_j$.
We define
\[
 \left\{
 \begin{aligned}
  &[A]_{mn}=\int_{\omega_i} k(x, \bar{\bt})\nabla \varphi_n\nabla \varphi_m=R_{\text{snap}}^T\bar{A}R_{\text{snap}},\\
  &[S]_{mn}=\int_{\omega_i} k(x, \bar{\bt}) \varphi_n\varphi_m=R_{\text{snap}}^T\bar{S}R_{\text{snap}},
\end{aligned}
\right.
\]
where $\bar{A}$ and $\bar{S}$ denote fine-scale matrices corresponding to the stiffness and mass matrices, respectively,  with the permeability $k(x, \bar{\bt})$.
We choose the smallest $M_i$ eigenvalues of the equation
  \[
  A\psi^i_k=\lambda_k S\psi^i_k
  \]
and take the corresponding eigenvectors in the snapshot  space  by setting $\psi_k^i=\sum_j \psi_{k,j}^i\varphi_j$, for $k=1,\cdots,M_i$, to form the reduced snapshot space, where $\psi_{k,j}^i$ are the coordinates of the vector $\psi^i_k$.

   Let $\{\chi_i\}_{i=1}^{N_H}$ be a set of partition of unity functions associated with the open cover $\{\omega_i\}_{i=1}^{N_H}$ of $\Omega$.
 Then we  multiply the partition of unity functions by the eigenfunctions to construct GMsFE space, i.e.,
  \[
  V_H=\text{span}\{\Psi_{il}: \Psi_{il}=\chi_i\psi^i_l: 1\leq i\leq N_H \quad \text{and}\quad 1\leq l\leq M_i\}.
  \]
We can  use  a single index  for  the multiscale basis function set  $\{\Psi_{il}\}$ and place them in  the following matrix
  \[
  R=[\Psi_1, \Psi_2,\cdots,\Psi_{M_v}],
  \]
where $M_v=\sum_{i=1}^{N_H} M_i$ denotes the total number of multiscale basis  functions. We note that once the matrix $R$ is constructed, it can be repeatedly used for simulation.

 Next we present the temporal discretization for the equation (\ref{examplelc}).   When $0<\gamma<1$,  the equation (\ref{examplelc}) is the subdiffusion equation. We use the  method in \cite{GLi} to discretize the fractional derivative and have
 \begin{equation*}
 \begin{split}
   \int_0^{t}(t-\tau)^{-\gamma}\frac{\partial u(x, \tau)}{\partial \tau}d\tau & = \sum_{k=0}^{n-1} \frac{u^{k+1}-u^k}{\Delta t} \int_{t_k}^{t_{k+1}} \frac{d\tau}{(t_{n}-\tau)^\gamma}+O(\Delta t) \\
   &=\frac{\Gamma(1-\gamma)}{s}\sum_{k=0}^{n-1}b_{k}(u^{k+1}-u^k)+O(\Delta t),
   \end{split}
  \end{equation*}
 where $u^{k}=u(x,t_k)$, $0=t_0<t_1\cdots<t_{n}=t$, $t_k=k\Delta t$ and
\begin{equation}\label{Gam1}
 \left\{
 \begin{aligned}
 &b_k:=(n-k)^{1-\gamma}-(n-k-1)^{1-\gamma}, \quad   k=0,1,\cdots,n-1,\\
 &s:=\Delta t^\gamma\Gamma(2-\gamma).\\
 \end{aligned}
   \right.
\end{equation}
When $1<\gamma<2$,  the $(\ref{examplelc})$ is the superdiffusion equation. We can use the method in \cite{Lc} to discretize the fractional  derivative and have
 \begin{equation*}
  \begin{split}
   \int_0^{t}(t-\tau)^{1-\gamma}\frac{\partial^2 u(x,\tau)}{\partial \tau^2}d\tau &=\sum_{k=0}^{n-1} \frac{u^{k+2}-2u^{k+1}+u^{k}}{\Delta t^{2}} \int_{t_k}^{t_{k+1}} \frac{d\tau}{(t_{n}-\tau)^{1-\gamma}}
   +O(\Delta t)^{2} \\
   &=\frac{{\Gamma(2-\gamma)}}{\tilde{s}}\sum_{k=0}^{n-1}\tilde{b}_{k}(u^{k+2}-2u^{k+1}+u^{k})+O(\Delta t)^{2},
  \end{split}
 \end{equation*}
 where
\begin{equation*}\label{Gam2}
 \left\{
 \begin{aligned}
 &\tilde{b}_k:=(n-k)^{2-\gamma}-(n-k-1)^{2-\gamma}, \ \ k=0,1,\cdots,n-1,\\
 &\tilde{s}:=\Delta t^\gamma\Gamma(3-\gamma).\\
 \end{aligned}
   \right.
\end{equation*}

Let $U^n$ be  the solution  at the $n-$th time level. Then we have the weak formulation for the subdiffusion equation $(0<\gamma<1)$,
\begin{equation}\label{diff1}
    \left\{
  \begin{aligned}
  \frac{1}{s}\tilde{a}\bigg(\sum_{k=0}^{n-1}(U^{k+1}-U^k)b_k, v\bigg)+a(U^{n}, v)&=(f(t_n), v),\quad \forall v\in V_H\\
  (U^0,v)&=(u(x,0), v),\quad \forall v\in V_H,
  \end{aligned}
     \right.
\end{equation}
where $s$ and $b_k$ are defined in (\ref{Gam1})  and
  \[
    a(u,v)=\int k(x,\bt)\nabla u\nabla vdx, \ \
    \tilde{a}(u,v)=\int q(x) uvdx.
  \]
 The weak formulation for the superdiffusion equation $(1<\gamma<2)$ reads
\begin{equation}\label{diff2}
    \left\{
  \begin{aligned}
  \frac{1}{\tilde{s}}\tilde{a}\bigg(\sum_{k=0}^{n-1}(U^{k+2}-2U^{k+1}+U^{k})\tilde{b}_k, v\bigg)+a(U^{n+1}, v)&=(f(t_{n+1}), v)\quad \forall v\in V_H,\\
  (U^0,v)&=(u(x,0), v) \quad \forall v\in V_H.
  \end{aligned}
     \right.
\end{equation}

For subdiffusion case $(0<\gamma<1)$,  we define
\begin{equation*}\label{c1}
c_k:=
\begin{cases}
b_0& k=0,\\
b_k-b_{k-1} & 1\leq k \leq n-1.
\end{cases}
\end{equation*}
For superdiffusion case $(1<\gamma<2)$, we define
\begin{equation*}\label{c2}
\tilde{c}_k=
\begin{cases}
n^{2-\gamma}-(n-1)^{2-\gamma}& k=0,\\
-2n^{2-\gamma}+3(n-1)^{2-\gamma}-(n-2)^{2-\gamma}& k=1,\\
(n+2-k)^{2-\gamma}-3(n+1-k)^{2-\gamma}+3(n-k)^{2-\gamma}-(n-k-1)^{2-\gamma}& 2\leq k \leq n-1,\\
-3+2^{2-\gamma}& k=n,\\
1& k=n+1.
\end{cases}
\end{equation*}

Then the weak formulation of $(\ref{diff1})$ can be rewritten as
  \[
  \tilde{a}(U^{n}, v)+sa(U^{n}, v)=\tilde{a}\big(\sum_{k=0}^{n-1}U^kc_k, v\big)+s\big((f(t_{n}), v)\big).
  \]
 We assume that $U^n$ has the expansion
  \[
  U^n=\sum_{j=1}^{M_v} u_{Hj}^n \Psi_j(x),
  \]
where $\{  \Psi_j(x) \}$ denote the GMsFEM basis functions. Let
  \[
  u_H^n=(u_{H1}^n,u_{H2}^n,\cdots,u_{HM_v}^n)^T.
  \]
Then for $k=1,\cdots,M_v$,
  \begin{eqnarray*}
  \label{c-eq}
  \begin{aligned}
  &\sum_{j=1}^{M_v} u_{Hj}^{n}\tilde a(\Psi_j, \Psi_k)+s\sum_{j=1}^{M_v}u_{Hj}^n a(\Psi_j, \Psi_k)\\
  &=\sum_{i=1}^n\sum_{j=1}^{M_v} c_i u_{Hj}^{i} \tilde a(\Psi_j, \Psi_k)+s(f^{n}, \Psi_k).
  \end{aligned}
  \end{eqnarray*}
 Let $B$, $K$ and $F$ be the weighted mass, stiffness matrices and load vector using FEM basis functions in fine grid, respectively. Then the equation gives the following algebraic system,
  \[
   R^TBRu_H^{n}+sR^TKRu_H^{n}=\sum_{i=1}^nc_iR^TBRu_H^{i}+sR^TF.
  \]
 If we define
 \[\tilde B=R^TBR,  \quad \quad \tilde K=R^TKR,\]
 then $u_H^n$ can be computed  by the iteration
  \begin{equation*}
  \label{iteration_c}
  u_H^{n}=\bigg(\tilde B+s\tilde K\bigg)^{-1}\bigg(\sum_{i=0}^{n-1}c_i\tilde{B}u_H^{i}+s R^TF\bigg).
  \end{equation*}
By using the multiscale basis functions, the solution of $(\ref{diff1})$ in fine grid can be  obtained by downscaling  through  the transformation $Ru_H^n$.

In a similar way to solving equation $(\ref{diff1})$, the GMsFEM solution of $(\ref{diff2})$ can be computed  by the iteration
  \begin{equation*}\label{iteration_c1}
  u_H^{n+1}=\bigg(\tilde B+\tilde{s}\tilde K\bigg)^{-1}\bigg(\sum_{i=0}^{n}\tilde{c}_i\tilde{B}u_H^{i}+\tilde{s} R^TF\bigg).
  \end{equation*}
We note that when GMsFEM is not applied, the full order model solution in fine grid is obtained by the iteration
  \begin{equation*}
  \label{iteration_f}
  u_h^{n}=\bigg(B+s K\bigg)^{-1}\bigg(\sum_{i=0}^{n-1}c_iBu_h^{i}+sF\bigg) \ \ \text{for} \ \ 0<\gamma<1,
  \end{equation*}
and
\begin{equation*}
  \label{iteration_f1}
  u_h^{n+1}=\bigg(B+\tilde{s}K\bigg)^{-1}\bigg(\sum_{i=0}^{n}\tilde{c}_iBu_h^{i}+\tilde{s} F\bigg)\ \ \text{for} \ \ 1<\gamma<2.
  \end{equation*}

By comparing  the GMsFEM model with the full order model, we see that  the size of $\tilde K$ and $\tilde B$ are $M_v \times M_v$, but the size of $K$ and $B$ are $N_h \times N_h$  ($M_v\ll N_h$).
Thus a much smaller system is solved in GMsFEM.  The matrix $R$ for multiscale basis functions is computed overhead and it can be repeatedly used for all the time levels. This significantly improves the efficiency for forward model simulations.

\begin{rem}
In the sequential data assimilation process,  ensemble members $\{\bt_j\}$ update and so $\bar{\bt}$ of ($\ref{MS-basis}$)  updates as well.
 We can update the GMsFEM basis matrix $R$  to improve the GMsFEM model.
\end{rem}

\subsection{Stochastic collation method using  $l_1$ regularized  least-squares}\label{SCM}
Stochastic collocation method is an efficient approach to approximate the solution of PDEs with random inputs. In this paper, we use stochastic collocation method (SCM) to obtain an expansion
for observation operator $\bb{H}(\bt)$    and efficiently evaluate the simulated observation ensemble.

  We use  generalized polynomial chaos (gPC) functions to represent $\bb{H}(\bt)$ by  $l_1$ penalized least-squares method. Let $i$ be a multi-index with $|i|=i_1+\cdots+i_{n_z}$ and $N_0$ be a nonnegative integer. The $N_0$th-degree gPC expansion of $\bb{H}(\bt)$ is then approximated  by a linear combination of gPC basis  $\{\Phi_i(\bt)\}_{i=1}^{P}$, i.e.,
\begin{equation}\label{o}
   \bb{H}(\bt)\approx \bb{H}^{N_0}(\bt):=\sum_{i}^{P}c_i\Phi_i(\bt), \quad \quad P=\frac{(N_0+n_z)!}{N_0!n_z!}.
\end{equation}
 The coefficients of expansion are obtained by choosing some collocation points and least-squares method. We first take $Q$ realizations $\{\bt^i\}_{i=1}^Q$ of $\bt$ in the support of prior distribution $p(\bt)$. Then for each $i=1,\cdots,Q$, we solve a deterministic problem at the node $\bt^i$ to obtain $\bb{H}(\bt^i)$.
   After we obtain all pairs $\big\{\bt^i, \bb{H}(\bt^i)\big\}$ ($i=1,\cdots,Q$), we are able to construct a approximation of $\bb{H} (\bt)$ such that $\bb{H}^{N_0}(\bt^i)=\bb{H}(\bt^i)$ for all $i=1,\cdots,Q$. Thus, $(\ref{o})$ can produce a system of linear equations
\begin{equation}\label{linea}
  \bf A\bf c=\bf b,
\end{equation}
where $\bff A\in \bb{R}^{Q\times P}$ is the matrix with the entries
  \[
  \bff A_{ij}=\Phi_j(\bt^{i}), \quad i=1,\cdots,Q,\quad j=1,\cdots,P.
  \]
  and the right term $\bf b$ satisfing
  \[
   b_i={\bb H}(\bt^i), \quad i=1,\cdots,Q.
  \]
  If we solve the system in the ordinary least-squares method, the system $(\ref{linea})$ should be overdetermined, i.e.,  $Q$ should be much larger than $P$.
  This means that we need solve the forward model for a large number of samples. To reduce the computation burden, we can take much fewer samples  ($Q\leq P$) and use $l_1$ regularized least-squares, i.e.,
  \[
  \min \|\bff c\|_{1} \quad \quad \text{s.t.} \quad \quad \bff {Ac=b},
  \]
which is equivalent to the optimization  problem
  \begin{equation}
  \label{l1-problem}
   \min_{\bff c} \|\bff {Ac-b}\|_2+\alpha\|\bff c\|_1 .
  \end{equation}
 We use the lagged diffusivity fixed point method \cite{CR} for the $l_1$ penalized least-squares problem.
 Due to the nondifferentiability of the $l_1$ norm, we take an approximation to the penalty $\|{\bf x}\|_1$ such as  $\sqrt{|{\bf x}|^2+\beta^2}$, where $\bf x$ is a scalar and $\beta$ is a small positive parameter. We denote the approximated penalty by
 \[
 \call { J(\bff c)}=\sum_{i=1}^P \psi(|c_i|^2),
 \]
where
 \[
 \psi(t)=2\sqrt{t+\beta^2}.
 \]
For any $\bff v\in \bb{R}^P$,
\begin{eqnarray*}
  \frac{d}{d\tau}\call { J(\bff{c}+\tau \bff v )} &=& \sum_{i=1}^P \psi'(|c_i|^2)c_iv_i = \langle \text{diag}(\psi'(\bff c))\bff c, \bff v\rangle,
\end{eqnarray*}
where diag $(\psi'(\mathbf{c}))$ denotes the $n\times n$ diagonal matrix whose $i$th diagonal entry is $\psi'(|c_i|^2)$, and $\langle \cdot ,\cdot \rangle$ denotes the Euclidean inner product on $\bb{R}^{P}$. From this we obtain the gradient
\begin{equation*}
  \bf grad\mathcal{J}(\mathbf{c})=\mathcal{L}(\mathbf{c})\mathbf{c} ,
\end{equation*}
where $\mathcal{L}(\mathbf{c})=\text{diag}(\psi'(\bff c))$ and  is positive semidefinite.
 For convenience, we present the lagged diffusivity fixed point method for the $l_1$ regularized least-squares problem (\ref{l1-problem})  in Algorithm $\ref{lagged}$.

\begin{algorithm}
\caption{Lagged diffusivity fixed point method for the $l_1$ regularized least-squares problem (\ref{l1-problem})}
  \textbf{Input}:
         $\nu$:=0, $\mathbf{c}_{0}$ := \text{initial guess}, $\alpha$\\
  \textbf{Output}: $\bf{c}$\\
  \text{begin fixed point iterations}  \\
  $\mathcal{L}_{\nu} :=\mathcal{L}(\mathbf{c}_{\nu})$;\\
  $\mathbf{g}_{\nu} :=A^\mathrm{T}(A\mathbf{c}_{\nu}-\mathbf{d})+\alpha\mathcal{L}_{\nu}\mathbf{c}_{\nu}$; \\
  $H=A^\mathrm{T}A+\alpha\mathcal{L}_\nu$; \\
  $\mathbf{s}_{\nu+1} := -H^{-1}\mathbf{g}_{\nu}$;\\
  $\mathbf{c}_{\nu+1} := \mathbf{c}_{\nu}+\mathbf{s}_{\nu}$; \\
   $\nu :=\nu+1$;\\
         \label{lagged}
\end{algorithm}

 We note that the forward model is  solved $Q$ times to obtain the sampling vector $\bff b$.  The accuracy of surrogate model can be ensured  using much fewer samples
 (i.e., $Q<P$)  when the $l_1$ regularized least-squares method is used. As the number of  unknown parameters increases, the simulation times for the forward model will increase significantly.
 To treat the challenge,  we use  GMsFEM to  solve the forward model on a coarse grid for each sample to  improve  the simulation efficiency.

\section{Two-stage ensemble Kalman filter using GMsFEM coarse model}
\label{Sur}

In this section, we present a two-stage EnKF using GMsFEM coarse model to accelerate posterior exploration and improve the sequential assimilation performance.
The goal of inverse problems is to identify an appropriate solution which can minimize the misfit  between the forward model and measurements.
 Here it is equivalent to solving the minimization  problem (\ref{cos}). In this paper, we use EnKF methods to solve the minimization problem.
 EnKF is a sampling  method and avoids expensive  gradient computation for solving  minimization problem.

 Although EnKF method can  avoid  lineralization and  repeated  sampling to  explore  posterior density,  it requires  to compute the forward model many times
 in each forecast step described  in algorithm $\ref{Seq}$.  When the number of unknown  parameters $\bt$ is large, we need a large number of ensemble members to estimate $\bt$.
   This implies that the computation of the forecast is very expensive.
    To improve the computation efficiency,  we construct a surrogate model based on sparse gPC  and GMsFEM to approximately represent the full order model.
    However, EnKF has an inherent  constraint assumption  for the prior distribution, which must be a Gaussian distribution.
    As we know, the support of Gaussian distribution is $\bb R$, but the posterior is often concentrated in a small portion of the entire prior support in many inference problems. Thus, it may be much more efficient to build a surrogate only over the important region of a posterior than the entire prior support. Inspired by the idea,  we proposed a two-stage EnKF method.
    In the first stage, we build a new prior by a very coarse GMsFEM model, where we exclude the unimportant region of the posterior, and the initial ensemble members are drawn from the new prior to enter the EnKF assimilation process.     The second stage is the surrogate model based EnKF.  We update the surrogate model dynamically when a new analysis is obtained.

The objective of the  surrogate model is to construct a representation that quantifies the primary features of the high-fidelity model while providing the computational efficiency required for uncertainty quantification. In stationary Bayesian inference, we may need to build surrogate model only once. However,  EnKF method integrates new measurement data in each   assimilation step.
Thus, the surrogate model need to be updated sequentially.  To construct the current  surrogate model,  we exclude the unimportant region by the previous analysis.

In the first stage, we can apply standard  EnKF method based on the very first few levels of measurement data to construct the new prior.  In the second stage,
we use only a few ensemble members from the previous analysis and $l_1$ regularized least-squares  to build the current surrogate model, which allows fast forward model evaluations to generate observation ensemble.  The outline of two-stage EnKF is presented in Algorithm $\ref{IN}$.
\begin{algorithm}
\caption{ Two-stage EnKF with unknown measurement noise algorithm}
  \textbf{Input}:  $n_1$ and $n_2$ ($n_1<n_2$), $M_1$ and $M_2$ ($M_1< M_2$),  $I_1$ and $I_2$ ($I_1< I_2$)\\
  \textbf{Output}: 
   final ensemble $\Theta_{I_2}$\\
  $\bf1$. As the GMsFEM  described in section $\ref{GMs}$, we build  the multiscale \\
  $~~~~~$basis functions matrix $R$.\\
  \textbf{First stage}:\\
  $\bf2$. Take $M_1$ basis functions to construct a very  coarse model\\
  $\bf3$. Set the initial ensemble $\Theta_0$ with $n_1$ samples\\
   $~~~~~$\textbf{for} $k=1, \cdots, I_1$\\
   $~~~~$ $~~~~~$Run algorithm $\ref{Seq}$, where $\bb{H}$ is the GMsFEM model based on  $M_1$ basis functions\\
   $~~~~~$\textbf{end for}\\
   \textbf{Second stage}:\\
  $\bf4$. The step $2$ and $3$ is to obtain $\Theta_{I_1}^{a}$ based on  the first $I_1$ levels of data information. The \\
  $~~~~~$support of $\Theta_{I_1}^{a}$ can be much smaller than the original prior support. We calculate the \\
   $~~~~~$mean and covariance of $\Theta_{I_1}^{a}$ obtained in step $3$. Take new initial ensemble $\Theta_0^{\text{new}}$ with\\
   $~~~~~$$n_2$ samples from the new prior. \\
  $\bf5$. \textbf{for} $k=I_1+1,\cdots,I_2$\\
   $~~~~$ $(a)$ Let $\Theta_{I_1}^{a}=\Theta_0^{\text{new}}$. Construct the surrogate model using multiscale  basis functions based \\
   $~~~~$ $~~~~$ on $\Theta_{k-1}^{a}$, where the forward model is solved by GMsFEM with $M_2$ basis functions.\\
   $~~~~$ $(b)$  Generate the observation ensemble by substituting ensemble members\\
   $~~~~$ $~~~~$  into the surrogate model based on  algorithm $\ref{lagged}$.\\
   $~~~~$ $(c)$  Run algorithm $\ref{Seq}$  and obtain the new analysis $\Theta^{a}_k$\\
  $~~~~$ \textbf{end for}\\
  \label{IN}
\end{algorithm}

The algorithm \ref{IN} can render an  accurate posterior for Gaussian models.  However, for non-Gaussian distributed parameter, we need to construct the surrogate model based on the support of non-Gaussian distribution, while the EnKF only works well for Gaussian prior. For this situation,   we use normal-score EnKF to perform non-Gaussian models.
 We note that the Legendre orthogonal polynomials can be used in surrogate model construction  if the support of unknown parameters is bounded.
\begin{algorithm}
\caption{Two-stage EnKF by normal-score transformation }
  \textbf{Input}: $n_1$ and $n_2$ ($n_1<n_2$), $M_1$ and $M_2$ ($M_1< M_2$),  $I_1$ and $I_2$ ($I_1< I_2$)\\
  \textbf{Output}: 
   final ensemble $\Theta_{I_2}$\\
   The step $\bf1$ and $\bf2$ is the same as $\bf1$ and $\bf2$ in algorithm $\ref{IN}$\\
   \textbf{First stage}:\\
  $\bf3$. Set the initial ensemble $\Theta_0$ with $n_1$ samples, $\Xi_0={\bb F}_0(\Theta_0)$\\
   $~~~~$\textbf{for} $k=1, \cdots, I_1$\\
   $~~~~$ $~~~~$ Run algorithm $\ref{non}$,  where $\bb{H}$ is the GMsFEM model based on  $M_1$ basis functions\\
   $~~~~$\textbf{end for}\\
   \textbf{Second stage}:\\
  $\bf4$. The step $2$ and $3$ is to obtain $\Theta_{I_1}^{a}$ and $\Xi_{I_1}^a$, and we calculate the mean and\\
   $~~~~$covariance of $\Theta_{I_1}^{a}$ obtained in step $3$. Take new initial ensemble $\Theta_0^{\text{new}}$ with $n_2$ \\
   $~~~~$samples from the new prior. \\
  $\bf5$. \textbf{for} $k=I_1+1,\cdots,I_2$\\
   $~~~~$  The step $(a)$ and $(b)$ are the same as in algorithm $\ref{IN}$\\
   $~~~~$ $(c)$ Run algorithm $\ref{non}$ and obtain the new analysis $\Theta_k^{a}$ and $\Xi_{k}^a={\bb F}_{k}(\Theta_k^a)$\\
  $~~~~$\textbf{end for}\\
  \label{Xl}
\end{algorithm}
\begin{rem}
In algorithm \ref{IN} and \ref{Xl},  $I_2$ is the total number of data assimilation steps. In step $3$, we can use  ES method to build the new prior  when some of the very first measurement data is uninformative.
\end{rem}

\section{Numerical examples}\label{exam}

In this section, we consider the time fractional PDE model (\ref{examplelc}) and estimate the model's unknown parameters and structures using the proposed two-stage EnKF. A few numerical results
will be presented for the estimation for different unknown sources of the dynamic model.
In Subsection \ref{channel}, we recover   a channel structure in permeability $k(x)$ when the fractional derivative $\gamma$ is known. In Subsection \ref{sourcese}, we estimate  the source locations
 when the diffusion type of equation (\ref{examplelc}) is  unknown. In Subsection \ref{coefsc}, we will recover a permeability pattern  when measurement noise is  unknown and is treated as a hyperparameter.

For the numerical examples,  we consider a dimensionless  square domain $\Omega=[0,1]\times[0,1]$ for spatial variable and $(0, T]$ for time, and we set the initial condition as $u(x,0)=0$ and in addition, for the superdiffusion equation, we set
\[
\frac{\partial u(x;t)}{\partial t}\bigg|_{t=0}=0.
\]
Measurement data are generated synthetically by using FEM in a fine grid with time step $\Delta t=0.001$, and the measurement noise is set to be $\sigma=0.01$. For any given realizations of $\bt$, we solve the time fractional diffusion-wave equation using GMsFEM with time step $\Delta t=0.002$. The regularization parameter $\alpha$ is set as $0.01$ in the lagged diffusivity fixed point algorithm. For all numerical examples, the number of samplers in constructing the gPC   surrogate is set to be the number of the gPC basis functions for approximation. We parameterize the Gaussian random fields by
 Karhunen-Lo$\grave{e}$ve expansion (KLE) with a given covariance function $C(x,y)$  and truncate the KLE to approximately represent the random fields. For the random filed $h(x,\omega)$ by the first $N$ terms can be represented by
\[
h(x,\omega)={\bb E}[h(x,\omega)]+\sum_{i=1}^N \sqrt\lambda_i\theta_i(\omega)\varphi_i(x),
\]
where $\bb{E}[\cdot]$ is the expectation operator, $\theta_i(\omega)\sim {\bf N}(0,1)$ and $(\lambda_i, \varphi_i)$ are the eigenpairs  of the eigenvalue problem
\begin{equation}\label{KLE}
\int_\Omega C(x,y)\varphi_i(y)dy=\lambda_i \varphi_i(x).
\end{equation}
 We sort the eigenvalues in ascending order, i.e., $\lambda_1\geq\lambda_2\cdots$, and their corresponding eigenfunctions are also sorted accordingly and $\{\theta_i\}_{i=1}^N$ are uncorrelated random variables. We will compare the estimation results obtained by using standard EnKF method with the proposed  two-stage EnKF.


\subsection{Recover  a channel structure in permeability field}\label{channel}

\begin{figure}[htbp]
  \centering
  \includegraphics[width=2.5in, height=1.7in]{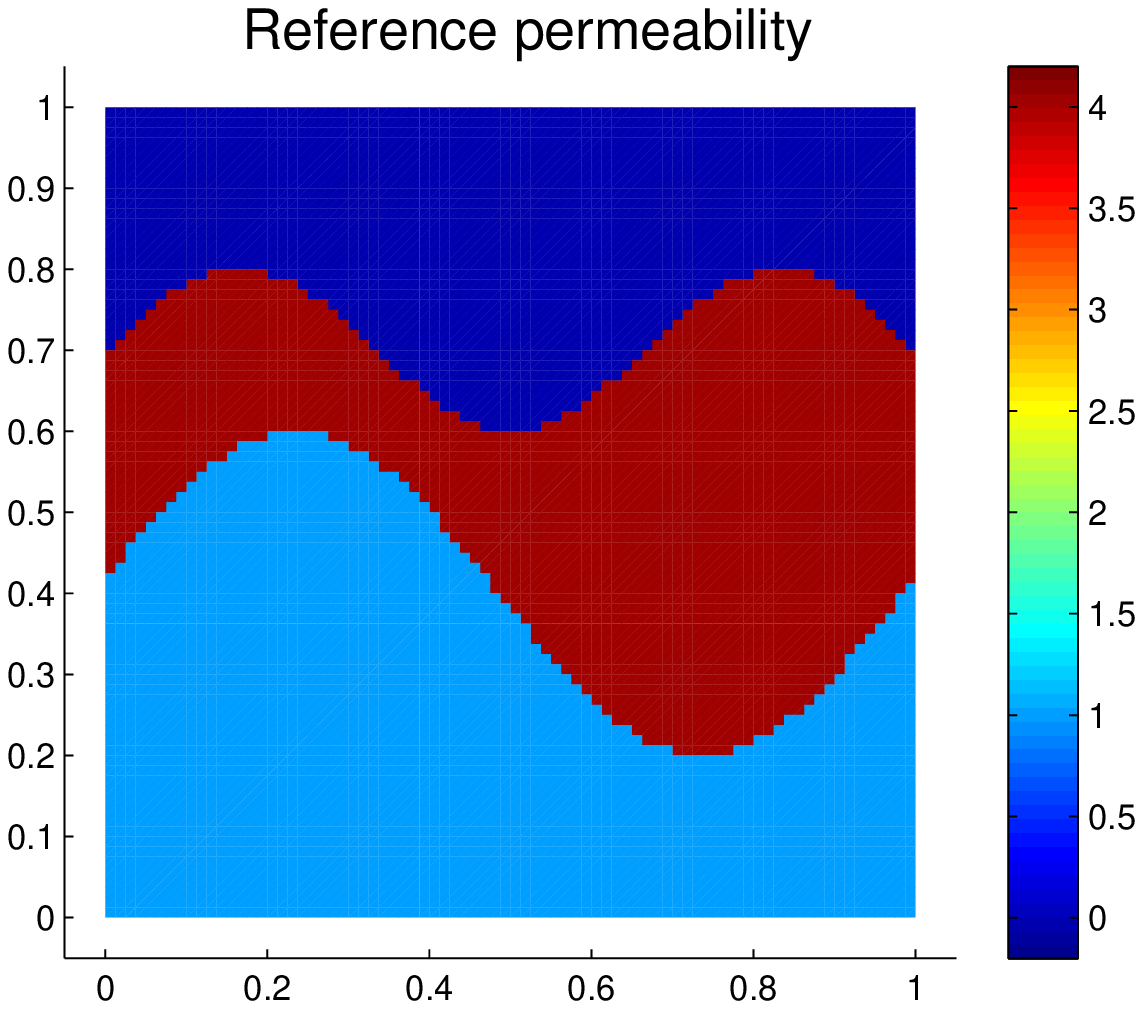}
  \includegraphics[width=2.5in, height=1.7in]{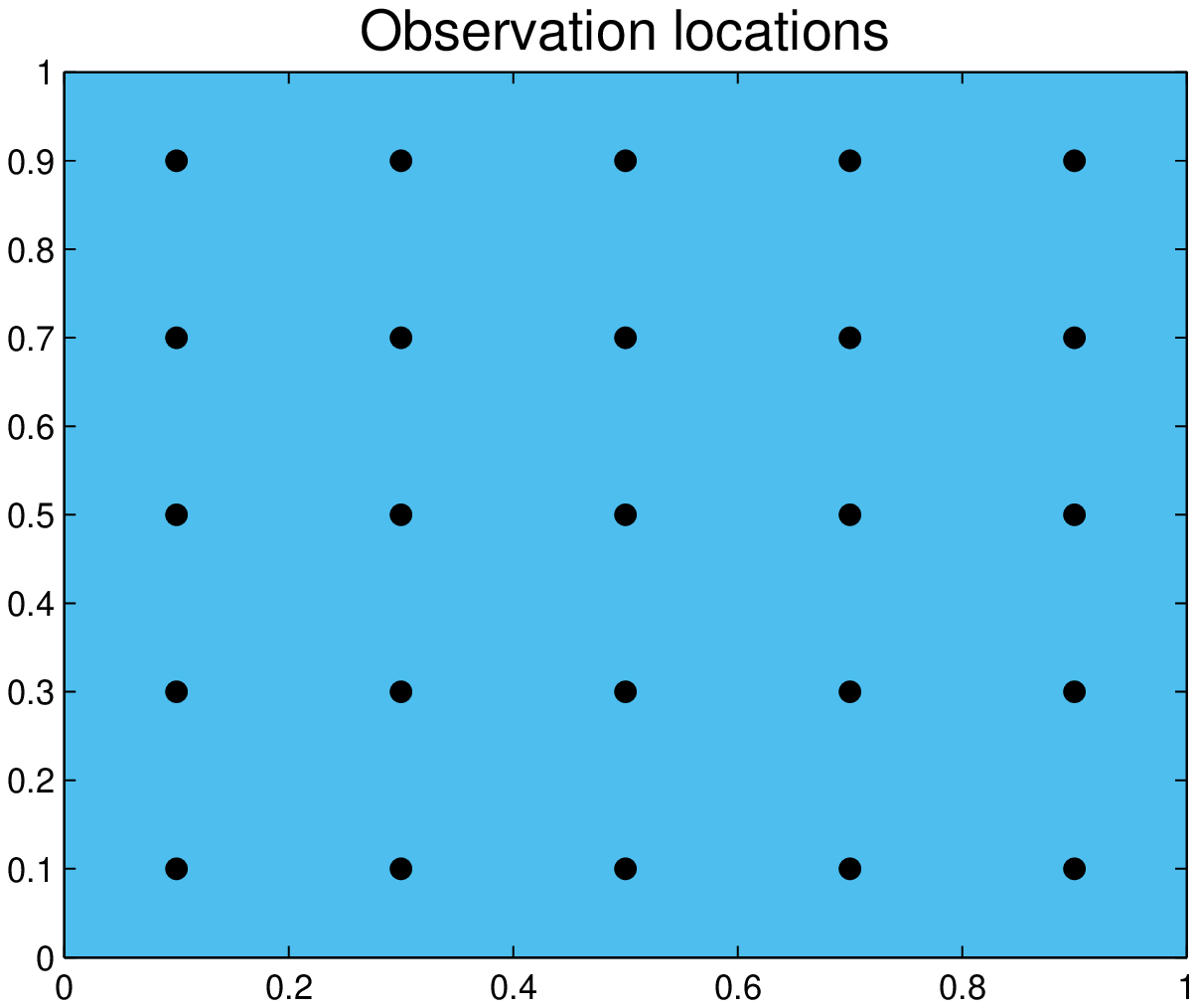}
  \caption{True profile of $\log k(x)$ (left)  and  observation locations (right).}\label{chh}
\end{figure}

In this subsection, we consider the  subdiffusion model (\ref{examplelc}) with mixed boundary condition, where Dirichlet boundary conditions is
\begin{equation}
u(0,y;t)=1 ,\quad \quad \quad u(1,y;t)=0,
\end{equation}
and there is no flow on the other boundaries. The source term is set as $f=10$, the end time is $T=0.11$, and the fractional derivative is given by $\gamma=0.5$. The permeability field is unknown here,
and we only  have the prior information of the permeability field, which is structured with a channel that lies between $y=0$ and $y=1$. The spatial domain is divided into $3$ parts by this channel and the permeability is a constant at each subregion. Thus, we can describe the boundaries of the channel by two curves $\Gamma_{1}$ and $\Gamma_{2}$, which can be expressed as

\[
  \Gamma_{1}(x)=\sum_{i=1}^{m_{1}} w_{i}^{1} \phi_{i}(x), \quad \Gamma_{2}(x)=\sum_{i=1}^{m_{2}} w_{i}^{2} \phi_{i}(x),
\]
where $\phi_{i}(x)$ are interpolation basis functions. In order to reduce dimension of unknown parameters, $\phi_{i}(x)$ are set to be Karhunen-Lo$\grave{e}$ve expansion (KLE) basis functions, i.e. $\phi_{i}(x)=\sqrt\lambda_i\varphi_i(x)$, and ${\bb E}[\Gamma_i(x)]=0\quad(i=1,2)$. The covariance function is set as
\[
C(x_1,x_2)=\xi^2\exp(-\frac{\|x_1-x_2\|^2}{2l^2})
\]
to solve the eigenvalue problem $(\ref{KLE})$, where $\xi=1$ and $l^2=0.1$. A  criterion which is adopted for the choice of truncated term $N$ is
\[
\frac{\sum_{j=1}^{N}\lambda_j}{\sum_{j=1}^{\infty}\lambda_j}>99.99\%.
\]
We truncate the KLE expansion for functions $\Gamma_i(x),\quad(i=1,2)$ by this criterion, and then we get $m_1=5$ and $m_2=5$.
In order to constrain the curves in the unit square domain, we make a bijective transformation
 \[
  \tilde\Gamma_{i}(x)=\frac{1}{2}+\frac{1}{\pi}\arctan(\Gamma_{i}(x)), \quad i=1, 2.
\]
Thus, $\tilde\Gamma_{1}(x)\in[0,1]$ and  $\tilde\Gamma_{2}(x)\in[0,1]$ can be guaranteed, i.e., the curves lie in the physical domain. Then we can construct the level set functions corresponding to the two curves as
\[
  L_{i}(x,y)=\text{Heaviside}(y>\tilde\Gamma_{i}(x)), \quad i=1,2.
\]
Hence, the random field can be parameterized as
\[
\log k(x,y)=c_{1}L_1L_2+c_{2}(1-L_1)L_2+c_3(1-L_1)(1-L_2),
\]
where $\{c_i\}_{i=1}^{i=3}$ are unknowns.
Thus, we have the unknown parameter vector
\[
\bt=(c_1,c_2,c_3,w_1^1,\cdots,w_{m_1}^1,w_1^2,\cdots,w_{m_2}^2).
\]
In this example, the reference permeability is generated by setting
\begin{eqnarray*}
  \Gamma_1(x) &=& 0.7+0.1\sin(3\pi x), \\
  \Gamma_2(x) &=& 0.4+0.2\sin(2\pi x+0.1),
\end{eqnarray*}
and the value of each subregion is set as $(c_1,c_2,c_3)=(0,4,1)$, which is shown in Figure \ref{chh}(left). In the example, the number of  artificial time steps for data assimilation is  $I_2=9$.
 Measurements are taken at time instances $0.012+0.01I:0.002:0.018+0.01I$ in each data assimilation step, where $I\in\{1,2,\cdots,I_2\}$ and the measurement locations are distributed on the uniform $5\times 5$ grid of the domain $[0.1, 0.9]\times [0.1, 0.9]$ as shown in Figure \ref{chh} (right).

The forward model is defined  on a uniform $80\times80$ fine grid, and we set the coarse grid $5\times5$ for GMsFEM simulation. We
construct the local snapshot space with dimension  $M^i_{\text{snap}}=20$ and  select $M_c=10$ multiscale basis functions at each block to construct the coarsen reduced order model and solve the optimization problem (\ref{cos}).

To construct a new prior, we use $5\times 10^3$ ensemble members by the standard EnKF method in the first stage. Then the new prior is constructed  by incorporating data information from the first three data assimilation steps. Seven local multiscale basis functions ($M_i=7$) are selected in constructing the basis  matrix $R$ to construct coarse model. Then, we construct the surrogate model after obtaining the new prior ensemble in the second stage. Eight local multiscale basis functions ( $M_i=8$) are selected  in constructing the matrix $R$ to construct the gPC surrogate model in the second stage. When the order of gPC is set as $N_0=3$, the number of random samplers is 560 for computing vectors $\bf b$ and $\bf A$ in Section \ref{SCM}. In this stage, the number of ensemble members is set as $10^4$. For a comparison, we also implement the standard EnKF, which directly uses the previous analysis as a prior to the current moment in data assimilation process.
In the standard EnKF, we also use GMsFEM with $8$ local GMsFE basis functions to solve forward model for $10^4$ ensemble members. But no new prior and no gPC surrogate model are build in the
standard EnKF method.  To complete the assimilation process,  the standard EnKF method takes  about $4.7$ hours for CPU time, while the two-stage EnKF just take about  $1$ hour CPU time.
Thus, the two-stage EnKF is much more efficient than the standard EnKF.

\begin{figure}[htbp]
  \centering
  \includegraphics[width=2.1in, height=1.7in]{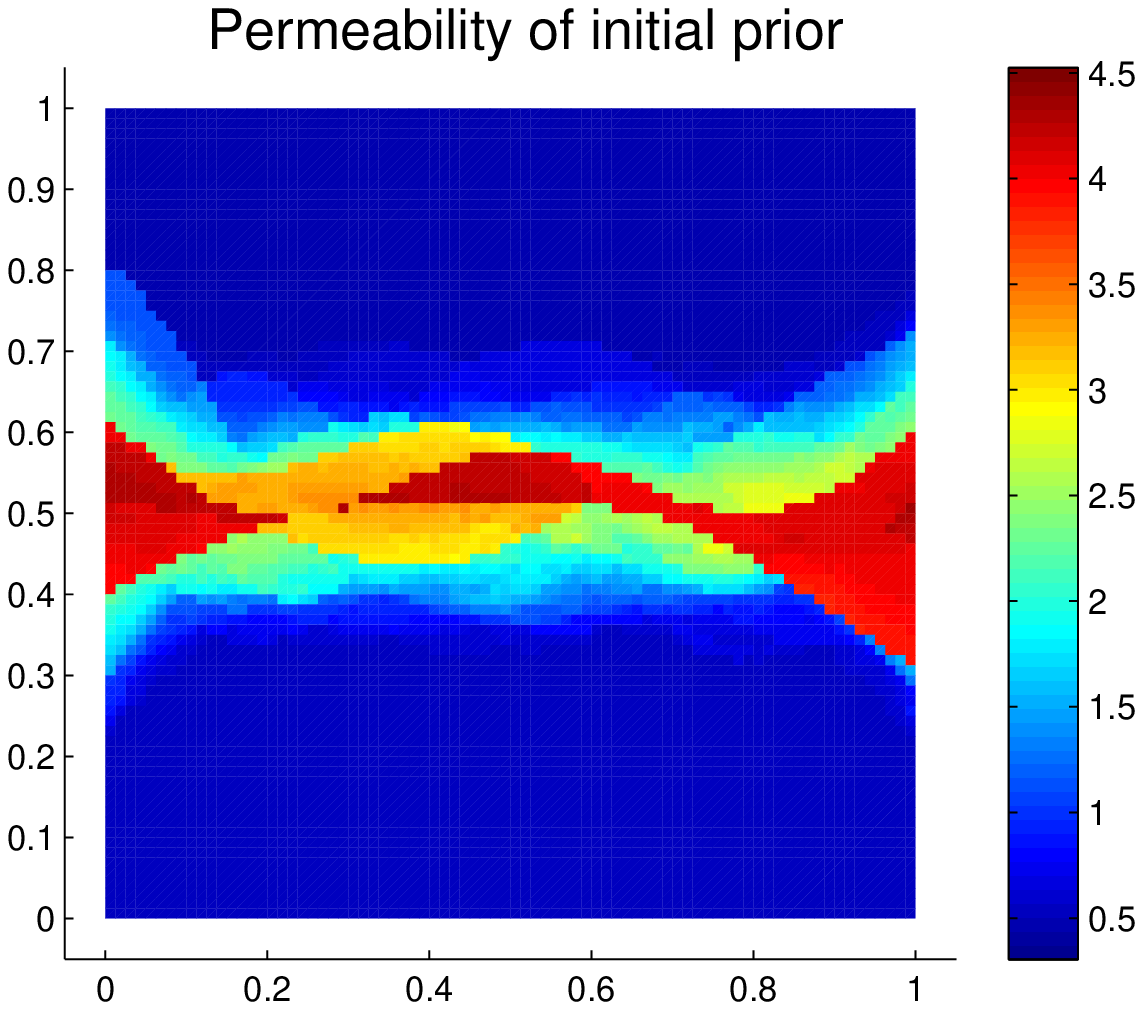}
  \includegraphics[width=2.1in, height=1.7in]{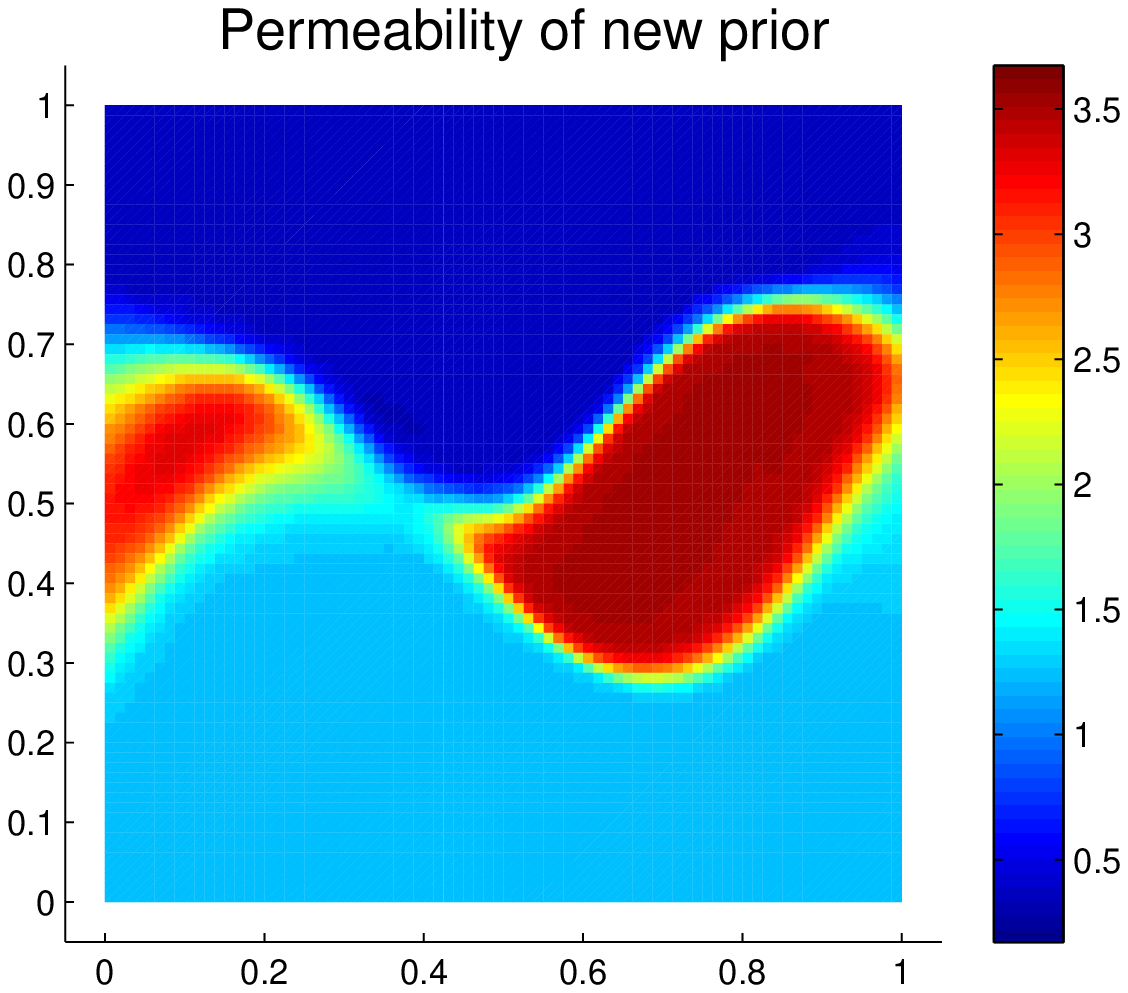}
  \caption{Log $k(x,\omega)$ of initial prior for two different EnKF methods (left), log $k(x,\omega)$ of new initial prior for the two-stage EnKF (right).}\label{ini}
\end{figure}

\begin{figure}[htbp]
  \centering
  \includegraphics[width=2.5in, height=1.7in]{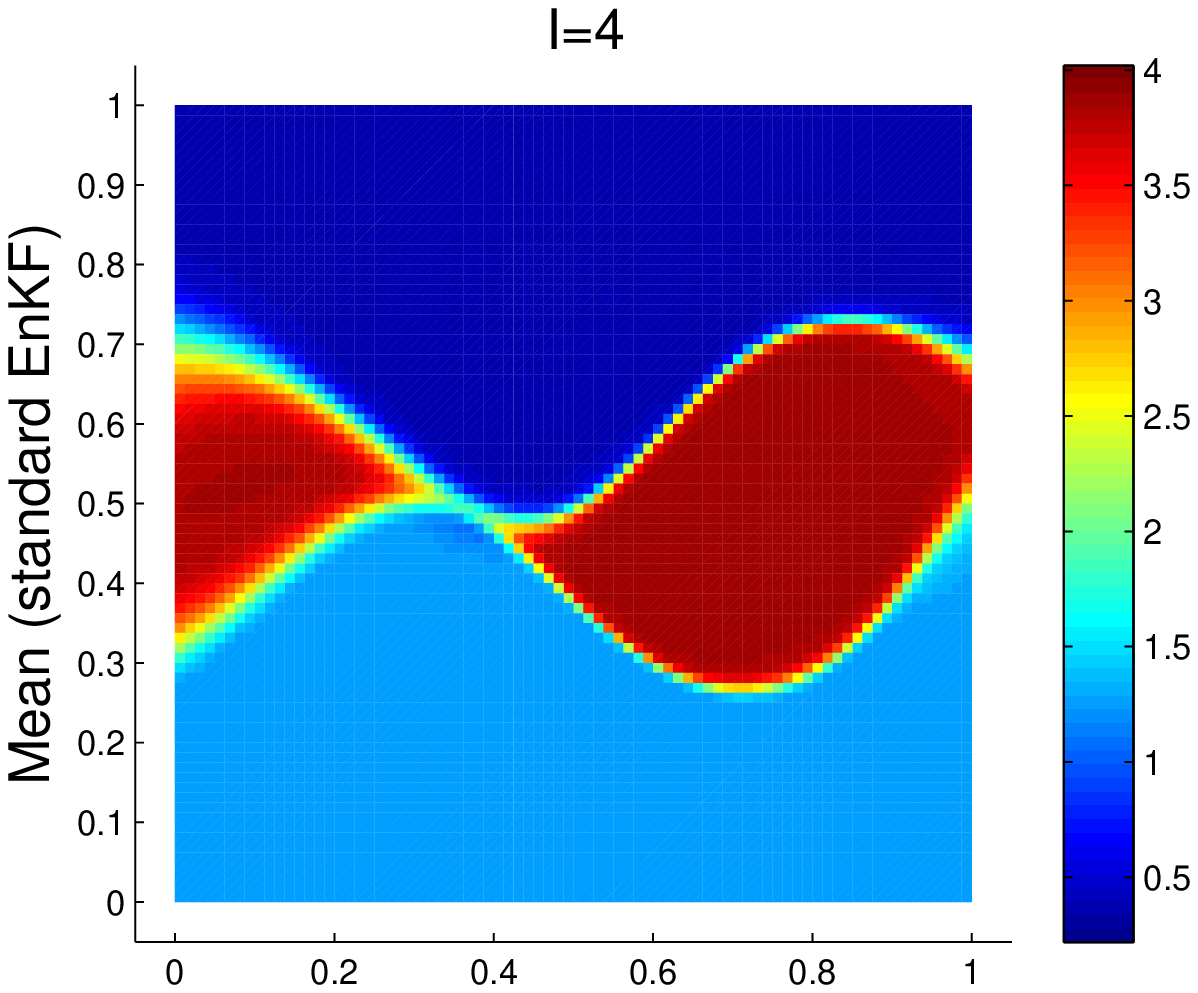}
  \includegraphics[width=2.5in, height=1.7in]{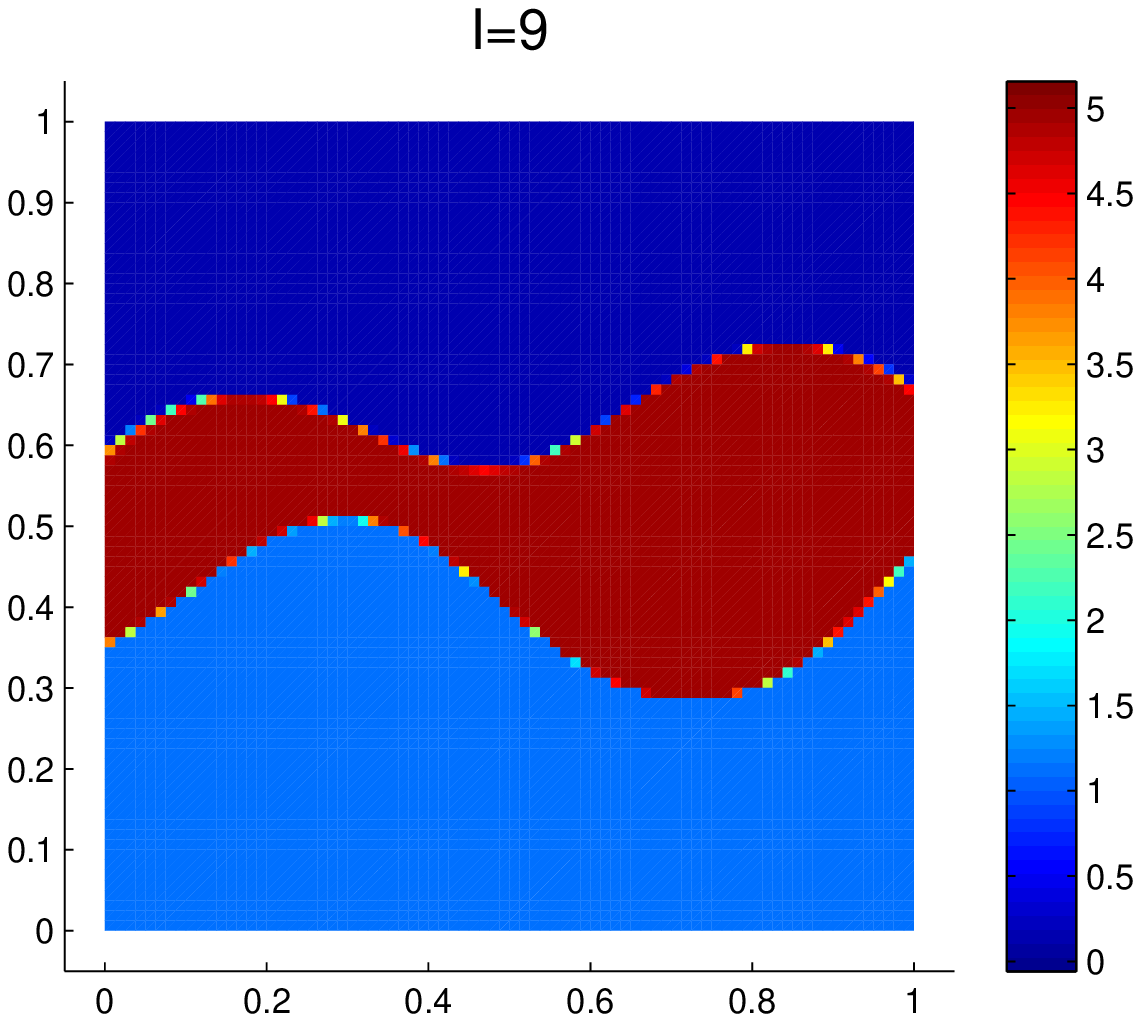}
  \includegraphics[width=2.5in, height=1.7in]{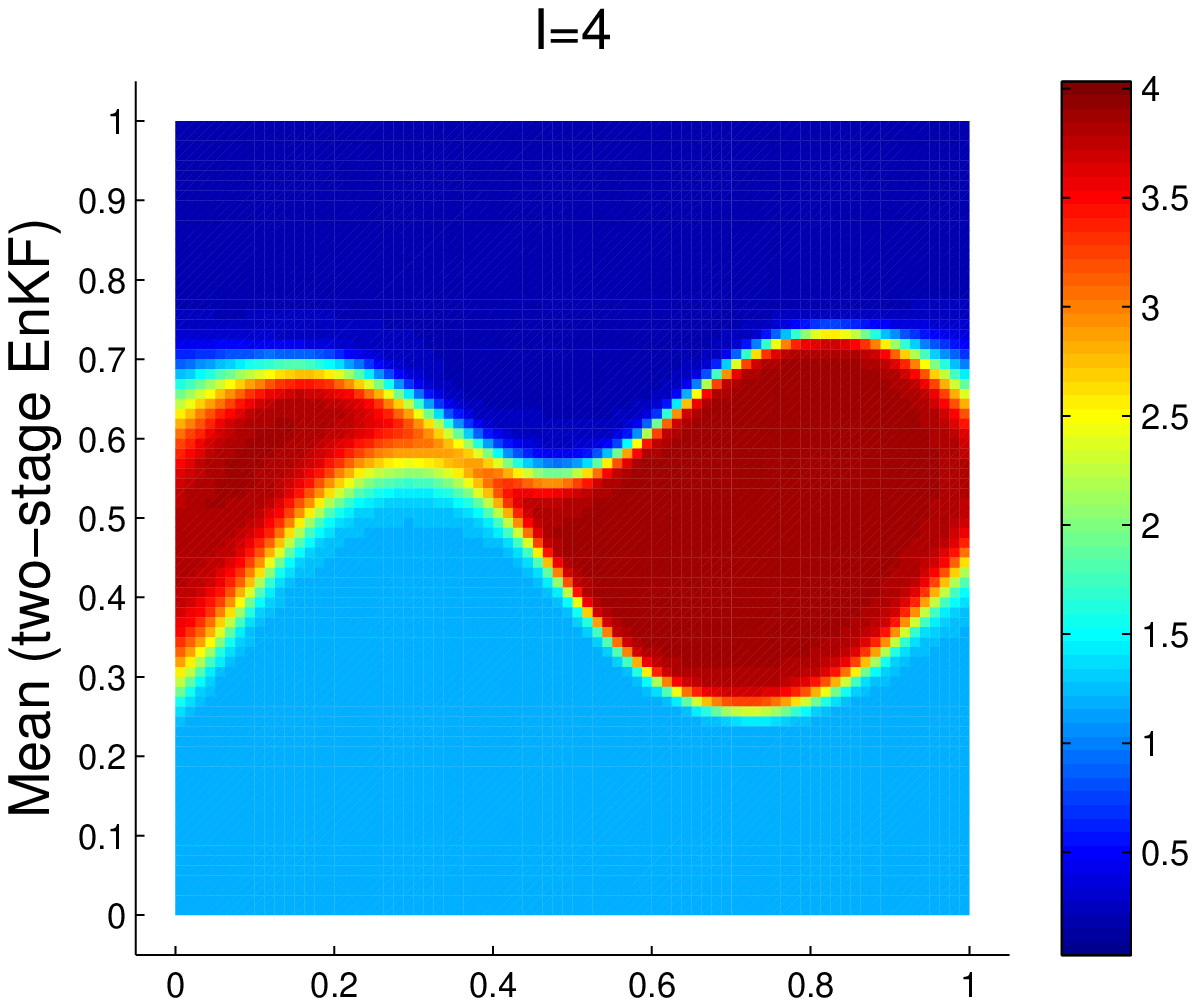}
  \includegraphics[width=2.5in, height=1.7in]{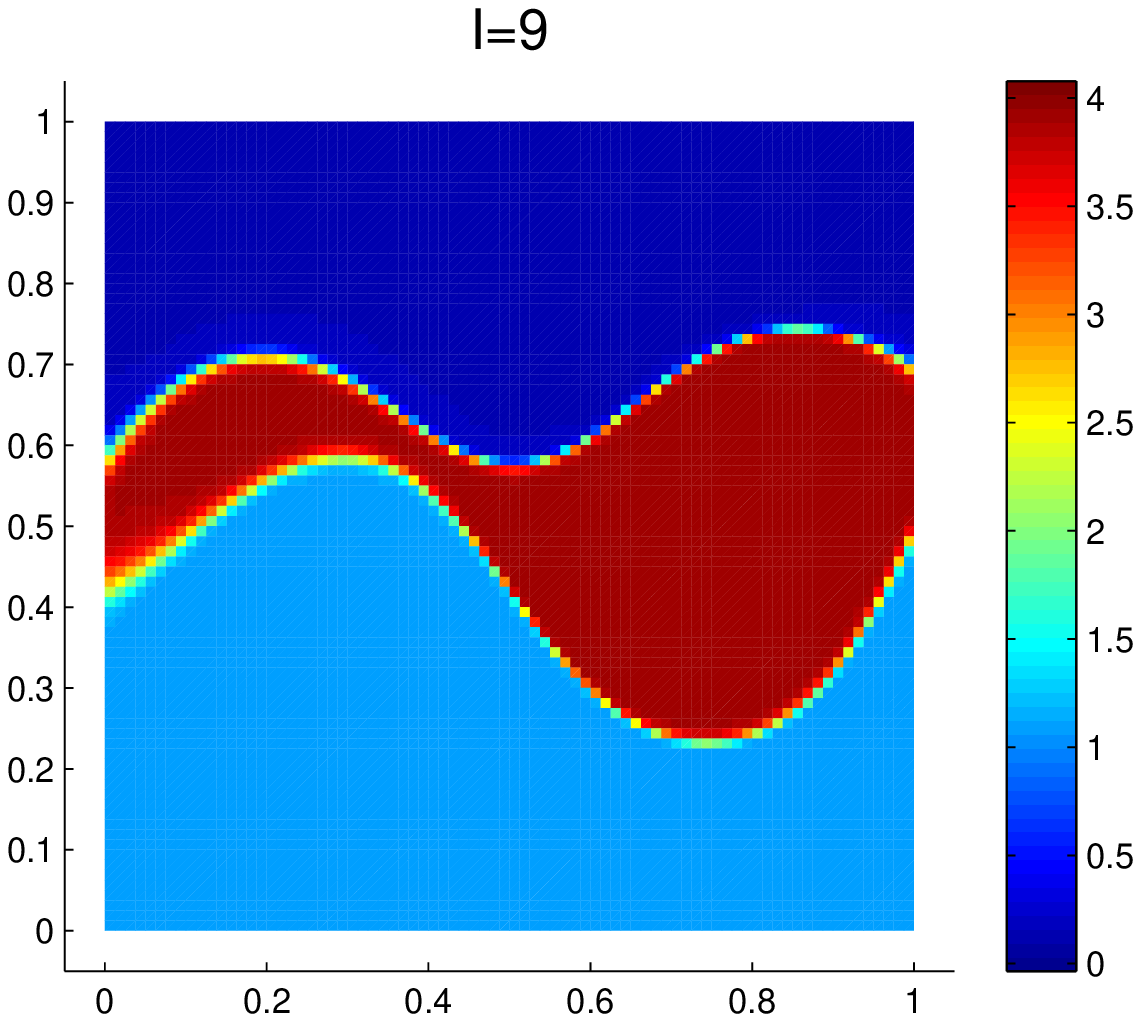}
  \includegraphics[width=2.5in, height=1.7in]{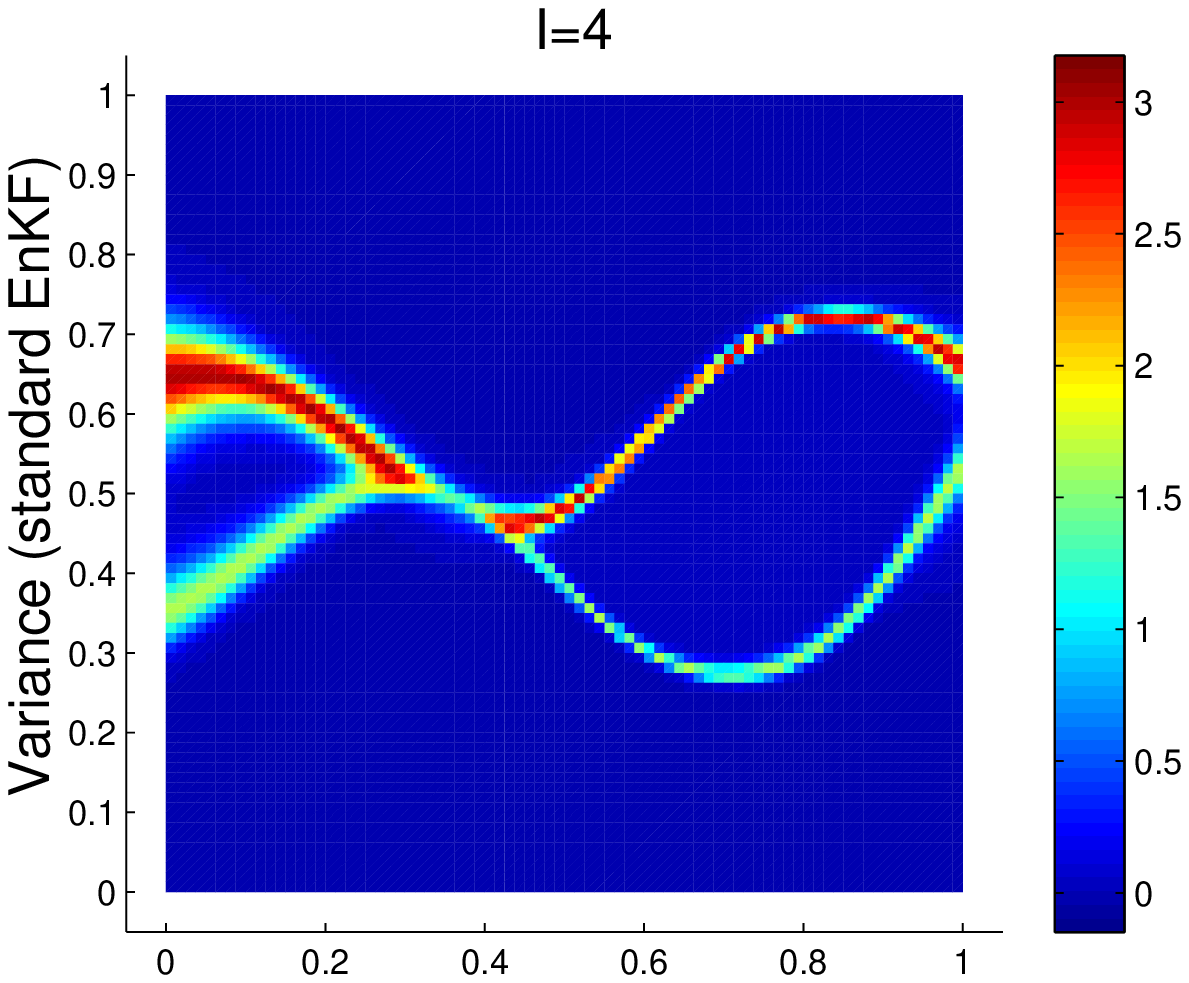}
  \includegraphics[width=2.5in, height=1.7in]{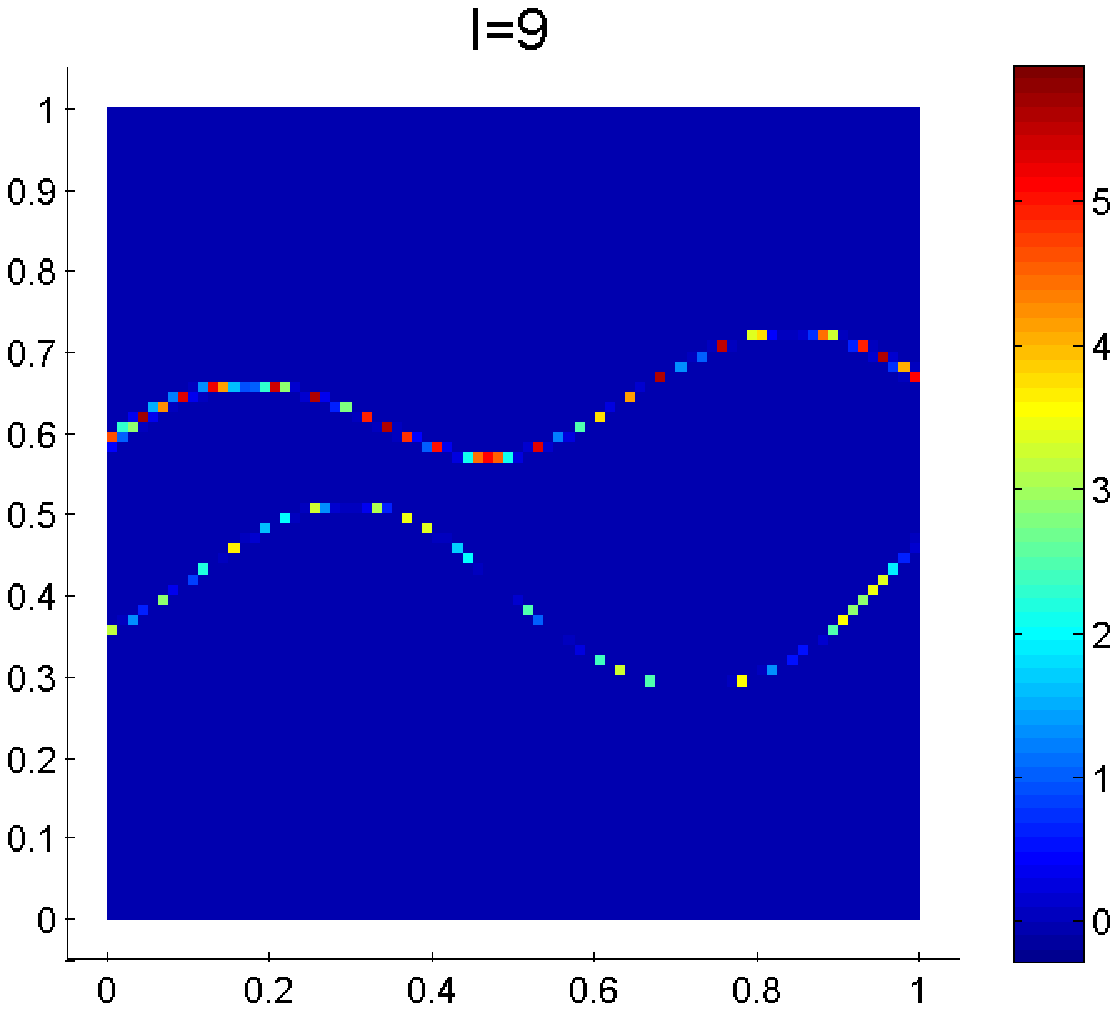}
  \includegraphics[width=2.5in, height=1.7in]{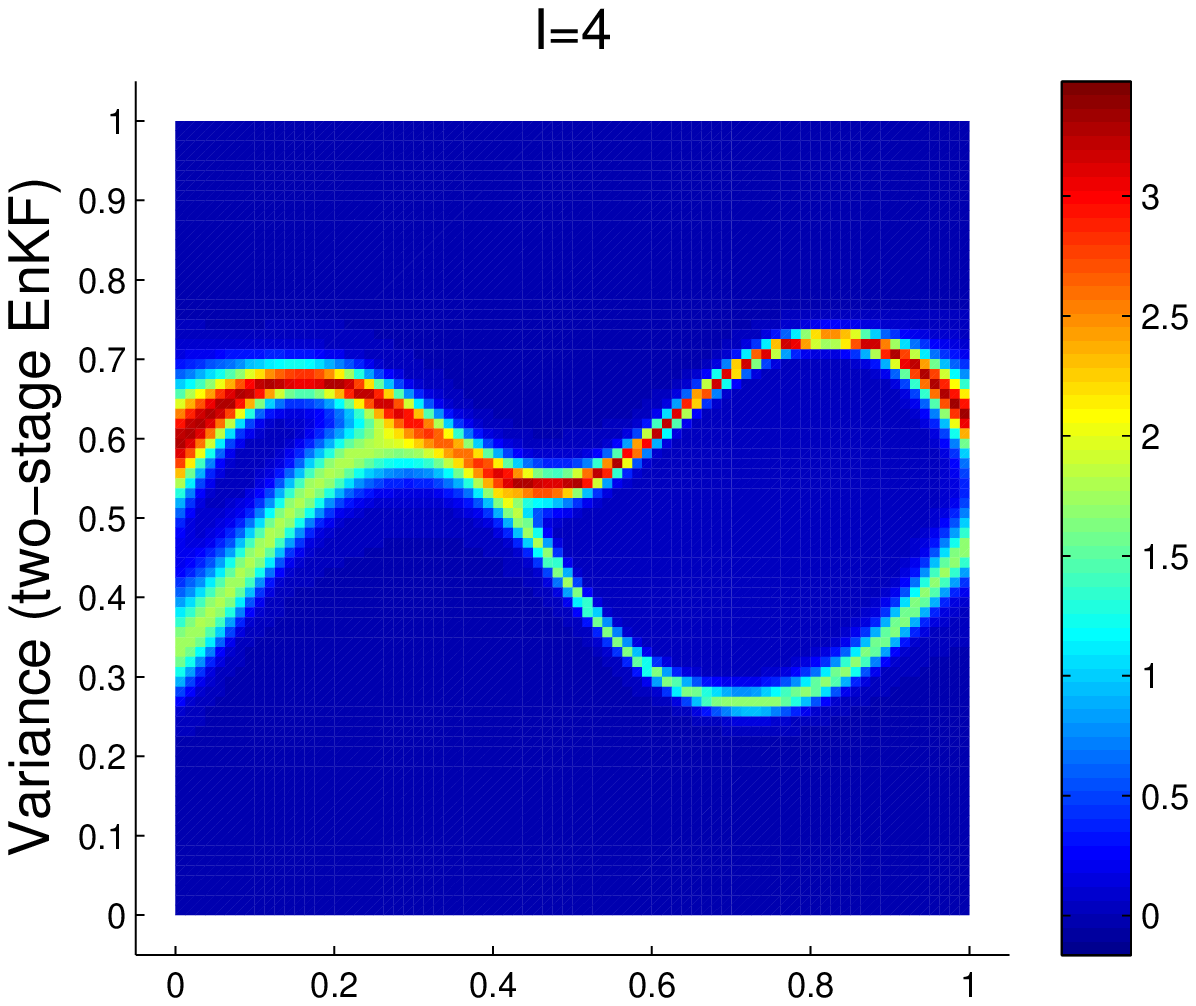}
  \includegraphics[width=2.5in, height=1.7in]{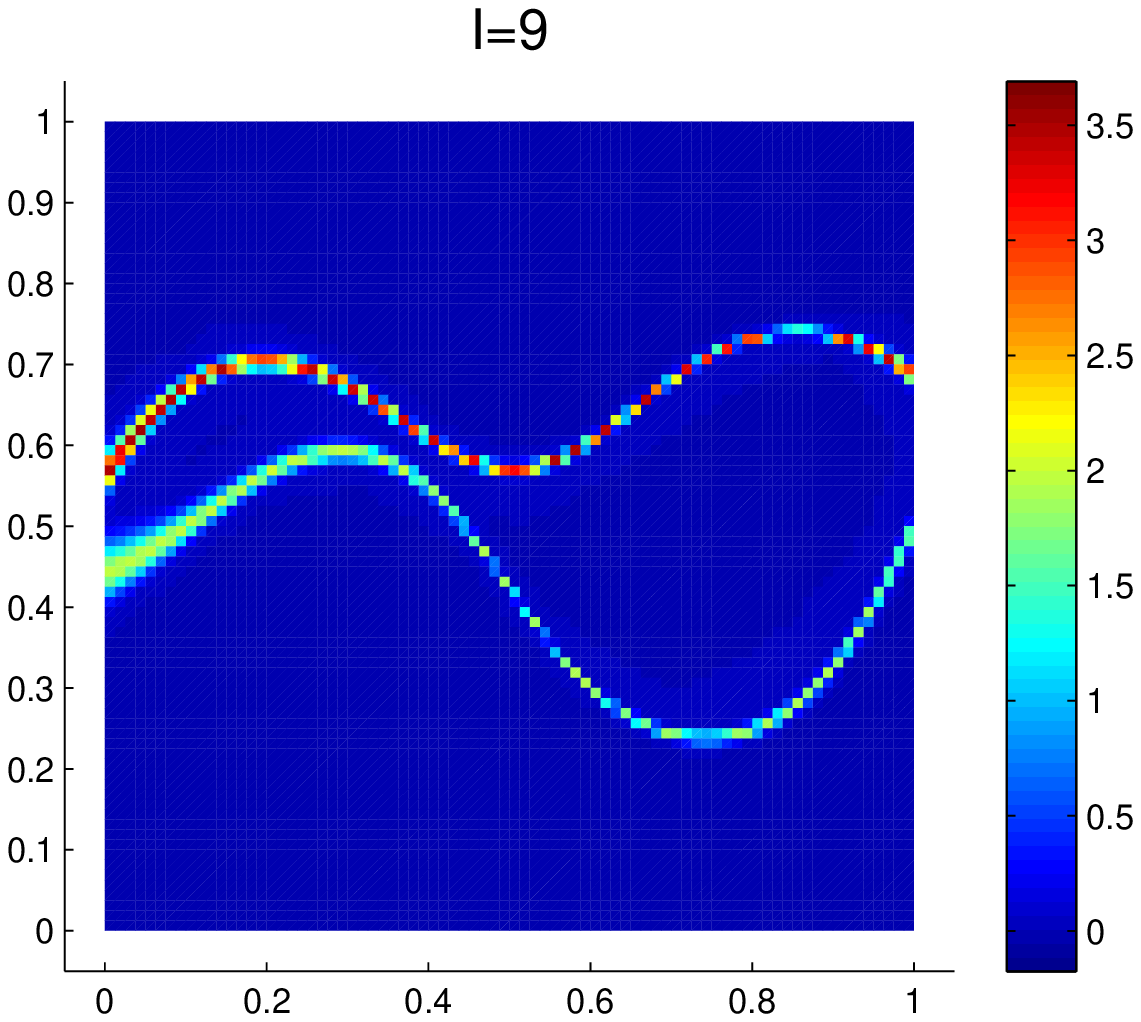}
  \caption{Mean and variance of log $k(x,\omega)$ by standdard EnKF (the first row and the third  row) and two-stage EnKF (the second row and the fourth row) at the assimilation step $4$ and $9$.}\label{chhigh1}
 \end{figure}

Figure \ref{ini} shows  the  initial prior ensemble (mean) and new prior ensemble (mean) for the  two-stage EnKF method in logarithmic scale.
We can see that the initial prior gives a very rough structure for  the channel, and  the new prior improves the channel pattern.  Then we compare the mean estimate by standard EnKF with the one by two-stage EnKF.  The posterior mean and variance for the two methods are presented   in Figure \ref{chhigh1}, from which we find the the two-stage EnKF  offers  more accurate estimation than standard EnKF.
The large  variance is concentrated around  the two boundaries of the channel because  two boundaries are unknown.

In order to measure the estimate accuracy in the update process, we define the relative errors $\varepsilon_{\Gamma_1}$ and $\varepsilon_{\Gamma_2}$ corresponding to the two curves after each update by
\begin{equation}
  \varepsilon_{\Gamma_1} := \frac{\|\tilde\Gamma_{1}(x)-\Gamma_{1}(x)\|}{\|\Gamma_{1}(x)\|}, \quad
  \varepsilon_{\Gamma_2} := \frac{\|\tilde\Gamma_{2}(x)-\Gamma_{2}(x)\|}{\|\Gamma_{2}(x)\|},
\end{equation}
where $\Gamma_1(x)$ and $\Gamma_2(x)$ are the reference boundary functions. We note that the first three data assimilation steps are used to produce the new prior for the two-stage EnKF.
We list  the relative errors  in table $\ref{chre}$.
As expected, the relative error gradually decreases when more measurement data is used in the inference.  The table shows that  the two-stage EnKF gives  more accurate estimates  than standard EnKF.

\begin{table}[tp] \caption {The relative error of $\tilde\Gamma_1$ and $\tilde\Gamma_2$}
\centering
\begin{threeparttable}
\begin{tabular}{ccccccc}\toprule
  \multirow{2}{*}{Update}    &\multicolumn{2}{c}{Two-stage EnKF} &\multicolumn{2}{c}{Standard EnKF} \\ \cmidrule(lr){2-3} \cmidrule(lr){4-5}
                & $\varepsilon_{\tilde\Gamma_1}$  &$\varepsilon_{\tilde\Gamma_2}$ & $\varepsilon_{\tilde\Gamma_1}$  &$\varepsilon_{\tilde\Gamma_2}$ \\\midrule
  initial prior     & 0.318     & 0.407 & 0.318    & 0.407 \\
  I=4    & 0.120     &0.129  & 0.183    & 0.181 \\
  I=5    & 0.111     &0.102  & 0.144   & 0.180\\
  I=6    & 0.109     &0.095  & 0.118    & 0.176 \\
  I=7    & 0.098     &0.086  & 0.120    & 0.179 \\
  I=8    & 0.099     &0.088  & 0.120    & 0.179\\
  I=9    & 0.096     &0.078  & 0.120    & 0.178 \\\bottomrule\label{chre}
\end{tabular}
\end{threeparttable}
\end{table}

\begin{figure}[tbp]
  \centering
  \includegraphics[width=1.5in, height=1.4in]{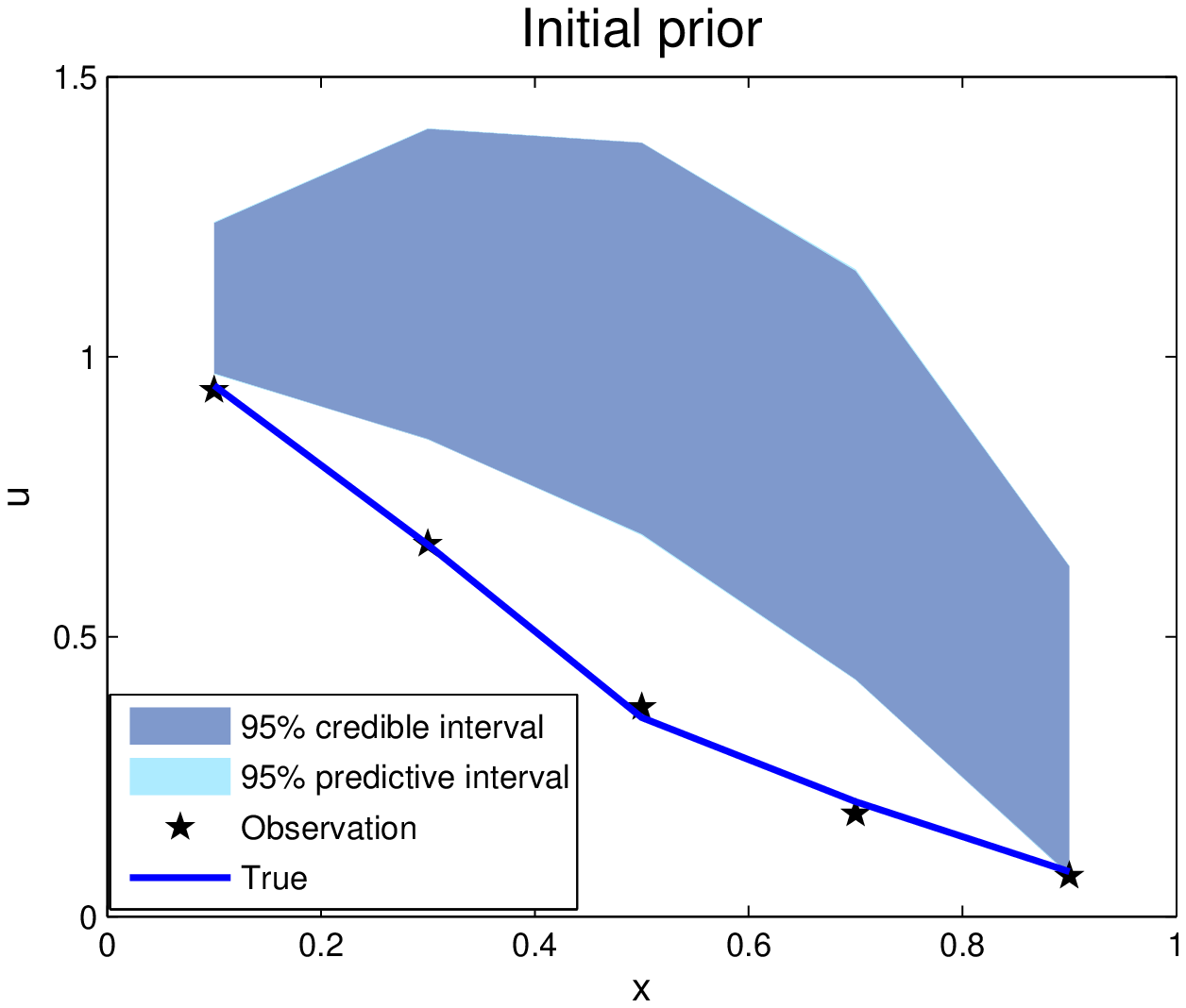}
  \includegraphics[width=1.5in, height=1.4in]{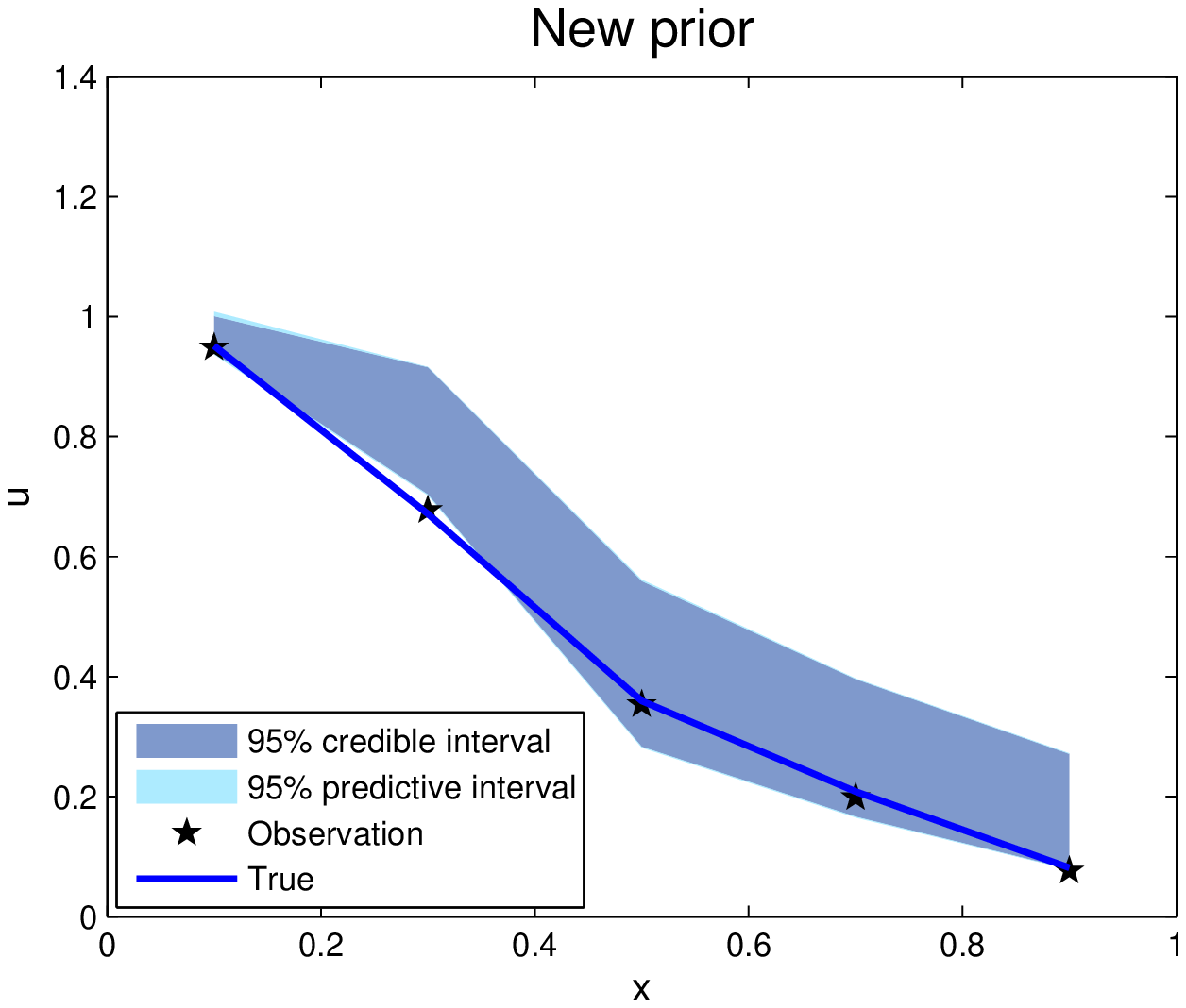}
  \includegraphics[width=1.5in, height=1.4in]{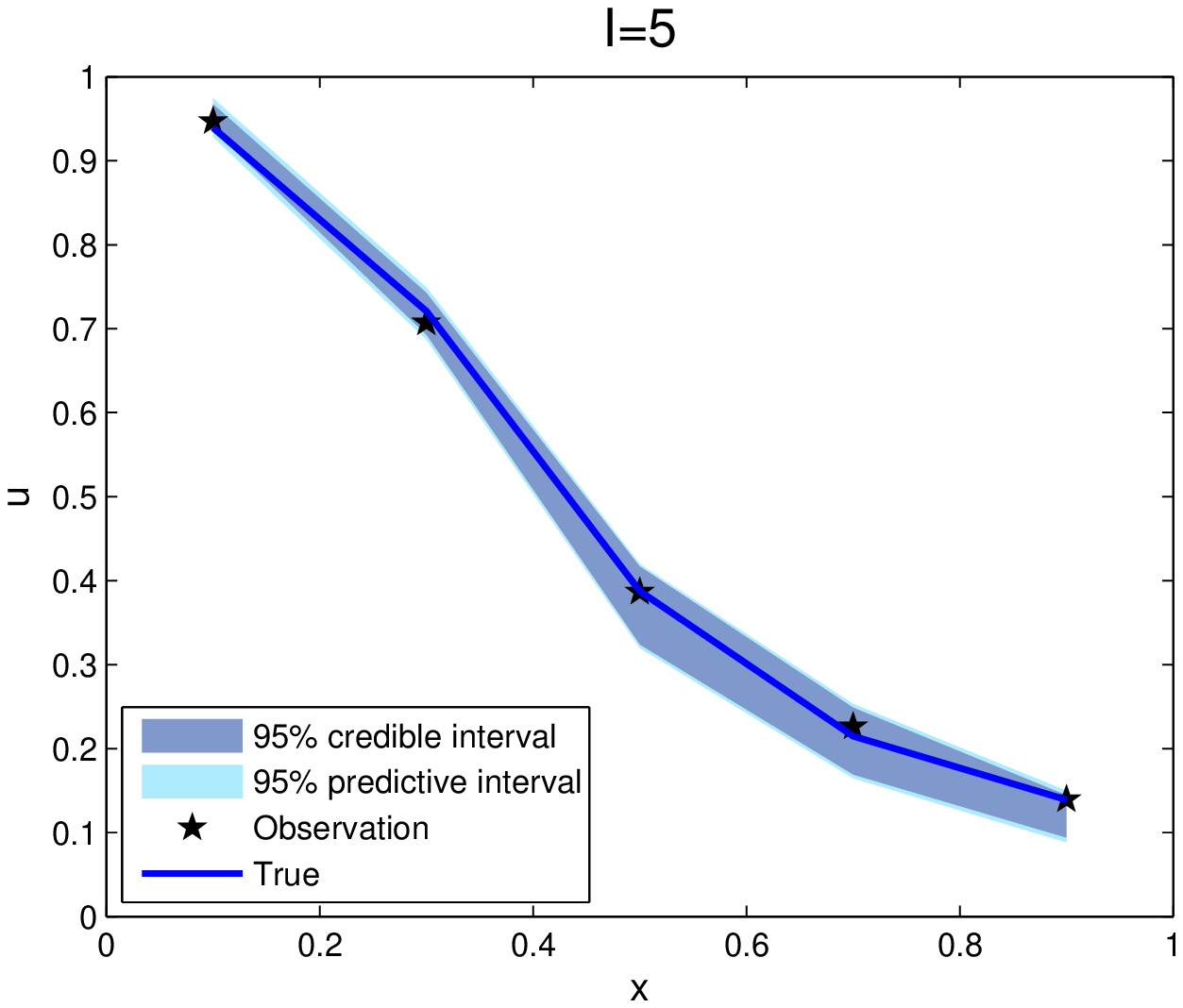}
  \includegraphics[width=1.5in, height=1.4in]{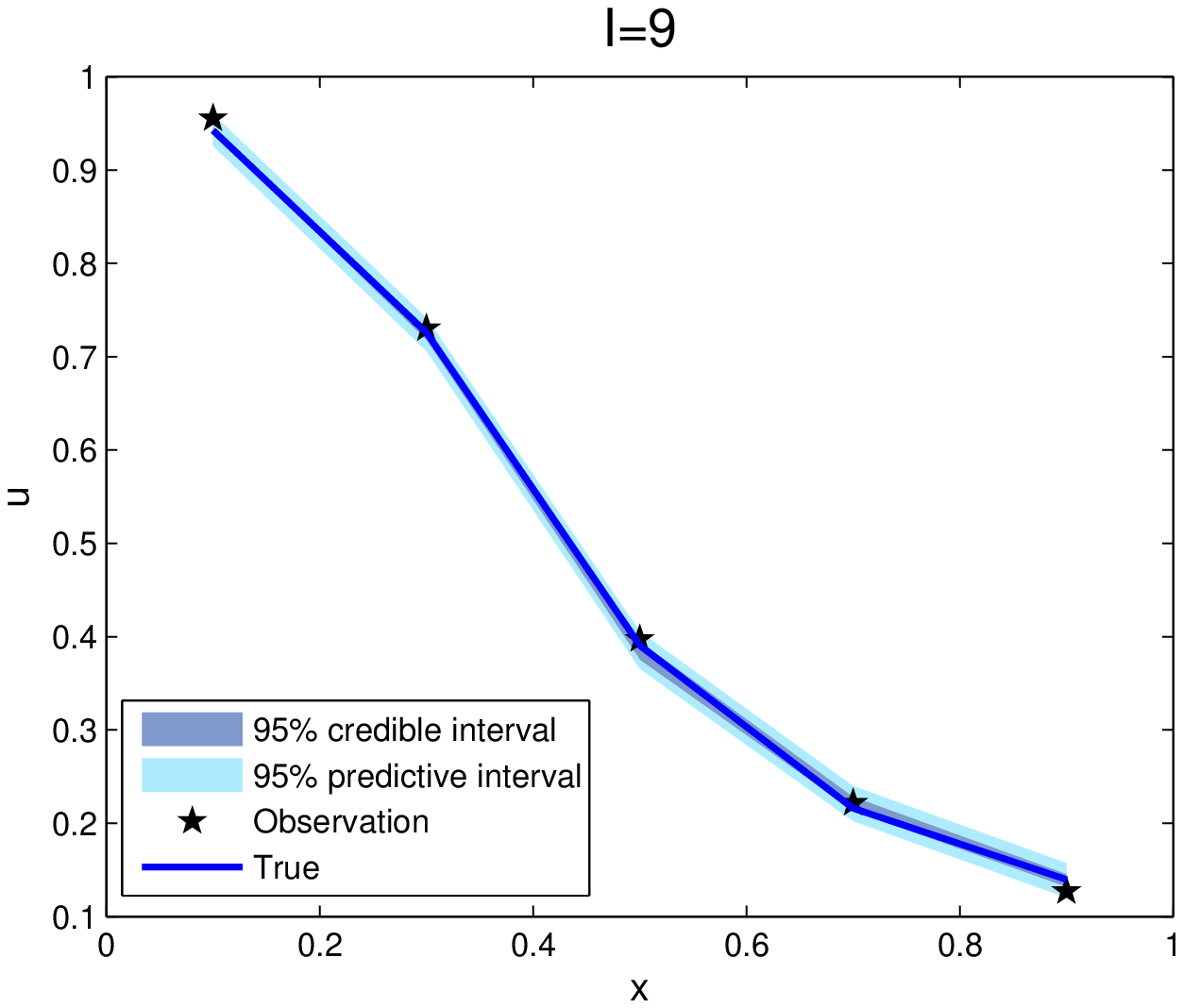}
  \caption{95\% predictive interval, 95\% credible interval, observation and true value by two-stage EnKF for $u(x, 0.5)$ at different assimilation steps.}\label{chpredict1}
 \end{figure}

\begin{figure}[tbp]
  \centering
  \includegraphics[width=1.5in, height=1.4in]{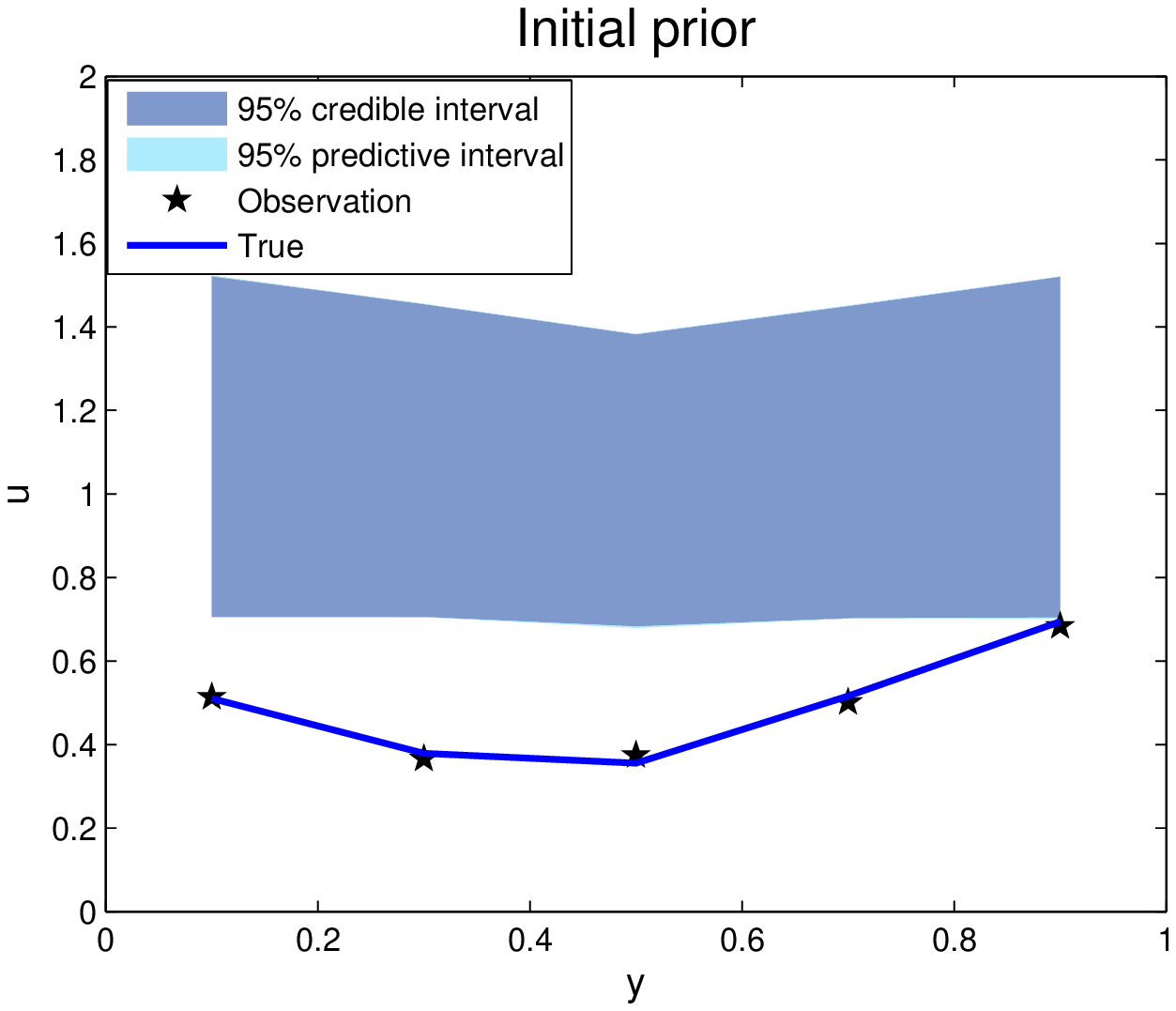}
  \includegraphics[width=1.5in, height=1.4in]{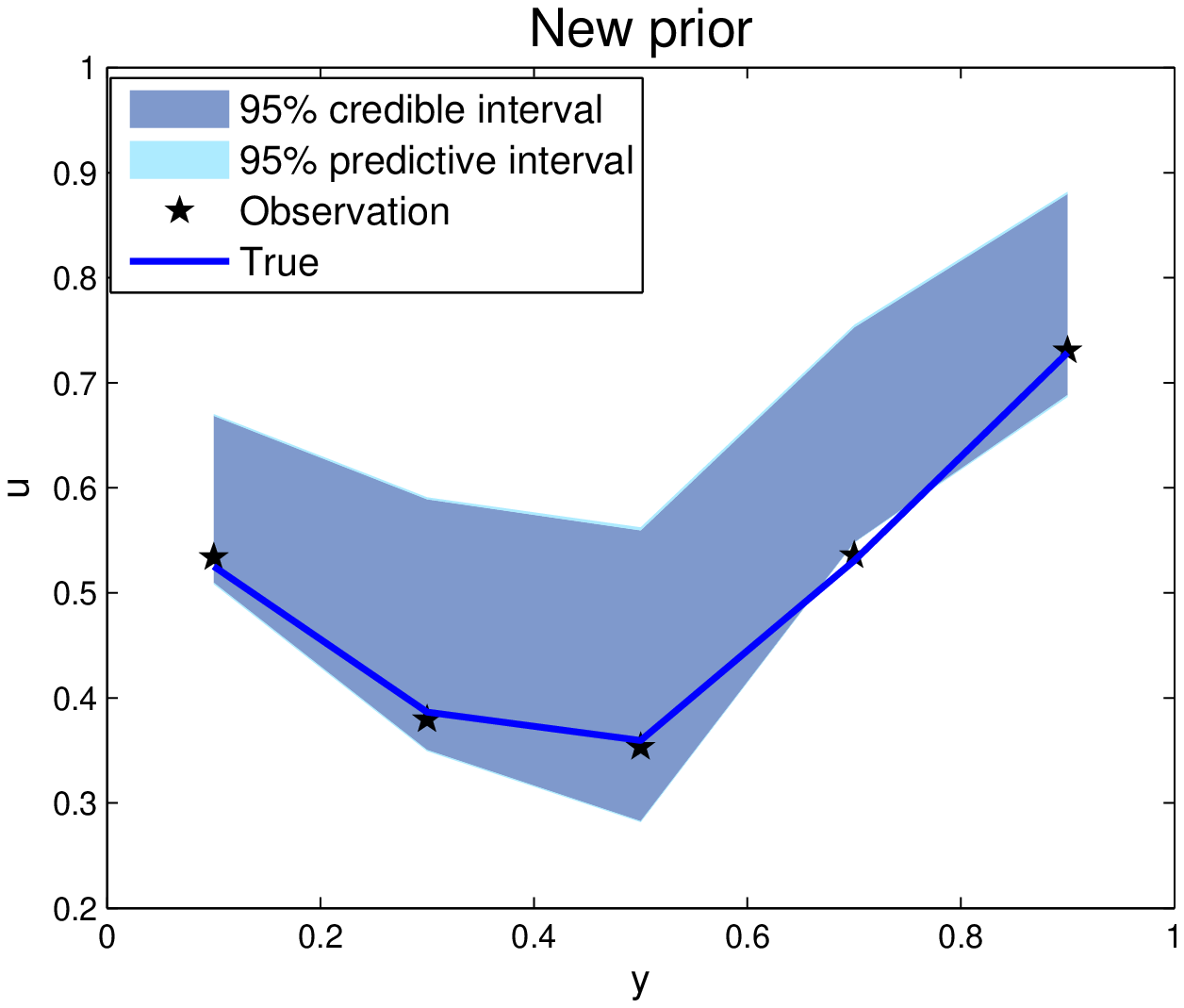}
  \includegraphics[width=1.5in, height=1.4in]{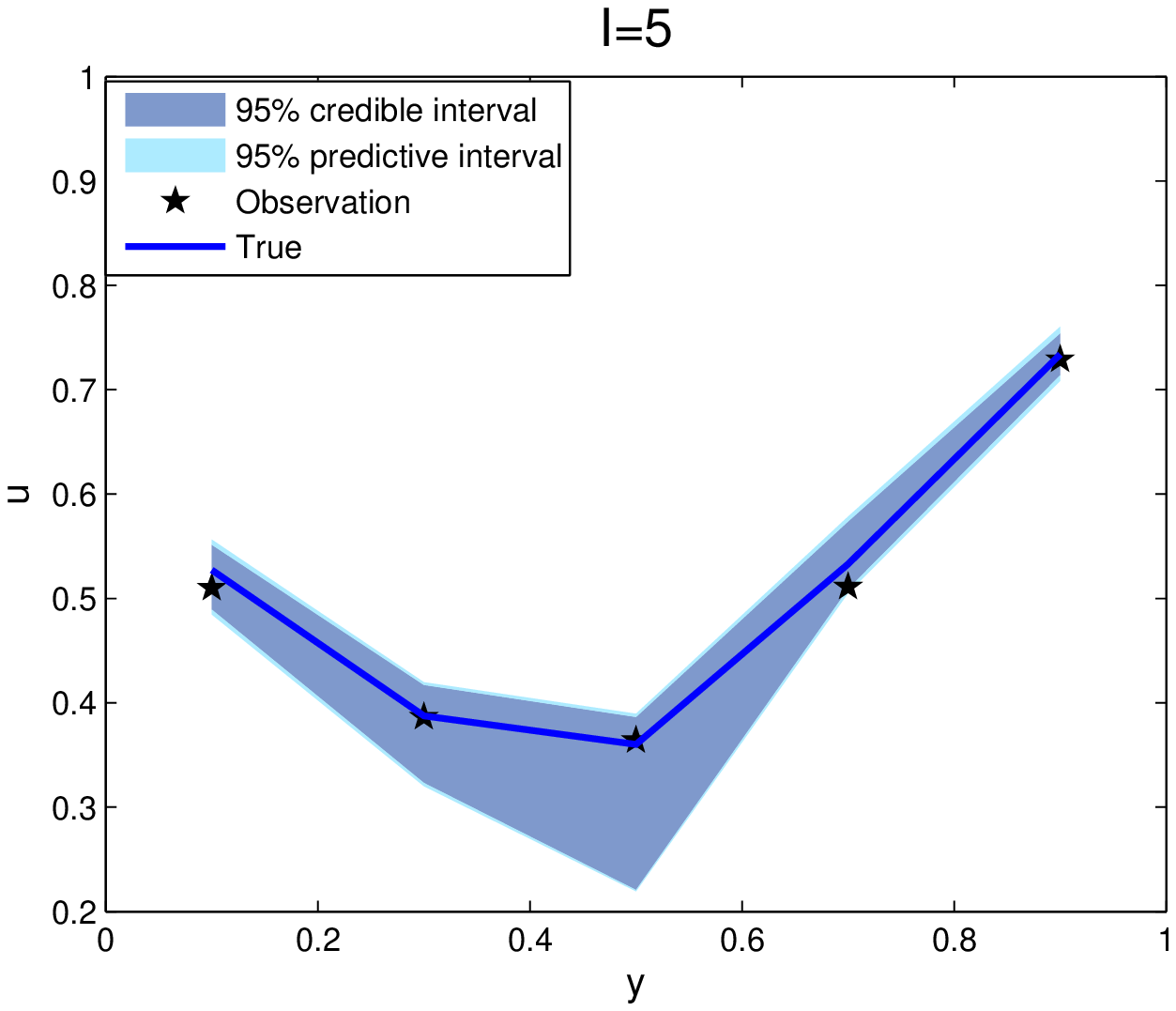}
  \includegraphics[width=1.5in, height=1.4in]{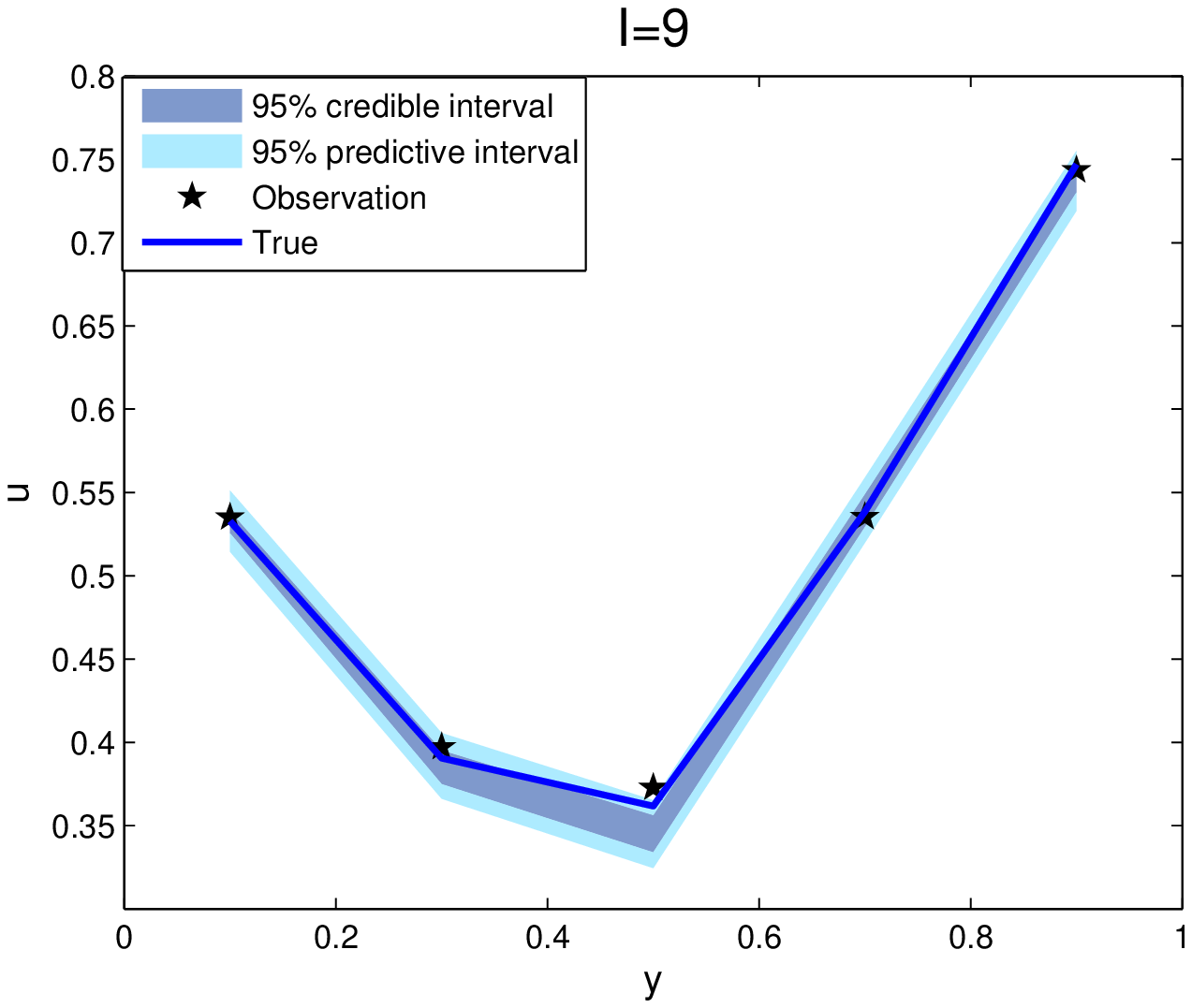}
\caption{95\% predictive interval, 95\% credible interval, observation and true value by two-stage EnKF for $u(0.5, y)$ at different assimilation steps.}\label{chpredict2}
 \end{figure}

To access the prediction using the posterior model,  we compute  the $95\%$ credible and predictive intervals for model response at $u\big((x, 0.5); t\big)$ and $u\big((0.5, y); t\big)$
for the two EnKF methods.  We note that  the realizations of model response of initial prior are  constructed by GMsFE model. This coarsen model error  is added to the estimated error variance to construct prediction intervals. As illustrated in Figure \ref{chpredict1} and \ref{chpredict2},  both the credible interval and predictive interval become  narrower as  assimilation moves on. This means that as the uncertainty from input $\bt$ decreases with respect to assimilation step and the uncertainty associated with the model fit and predictions decreases.  The observation/measurment data are almost contained in the predictive intervals when sufficient data is used in the posterior model.

\begin{figure}[tbp]
  \centering
  \includegraphics[width=3.3in, height=2.5in]{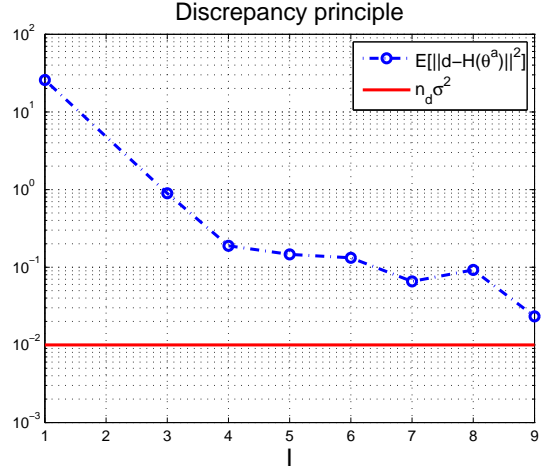}
  \caption{$\bb{E}[\|\bm{d}-\bb{H}(\bt^{a})\|^{2}]$ and   $n_d\sigma^{2}$ via  the assimilation time. We note $I=1$ denotes the initial prior time and $I=3$ the new prior time.}\label{chUnba}
 \end{figure}

For this example, the measurement noise is the additional Gaussian type. Let $\bt^*$ denotes the truth  parameter, the discrepancy between the data and forward model is defined by
\[
[\bm{d}-\bb{H}(\bt^*)]_i\sim {\bf N}(0,\sigma^2), \quad \quad i=1,\cdots,n_d,
\]
where $n_d$ is the number of measurements and its expectation is given by
\[
\bb{E}[\|\bm{d-\bb{H}(\theta^*)}\|^2]=n_d\sigma^2.
\]
We plot $\bb{E}[\|\bm{d}-\bb{H}(\bt^a)\|^2]$ against data assimilation step in Figure $\ref{chUnba}$. Then we can see the expectation of the discrepancy tends to $n_d\sigma^2$ as assimilation time moves on.
 This implies that  the final ensemble mean is an accurate estimate.

\begin{figure}[tbp]
  \centering
  \includegraphics[width=2.1in, height=1.7in]{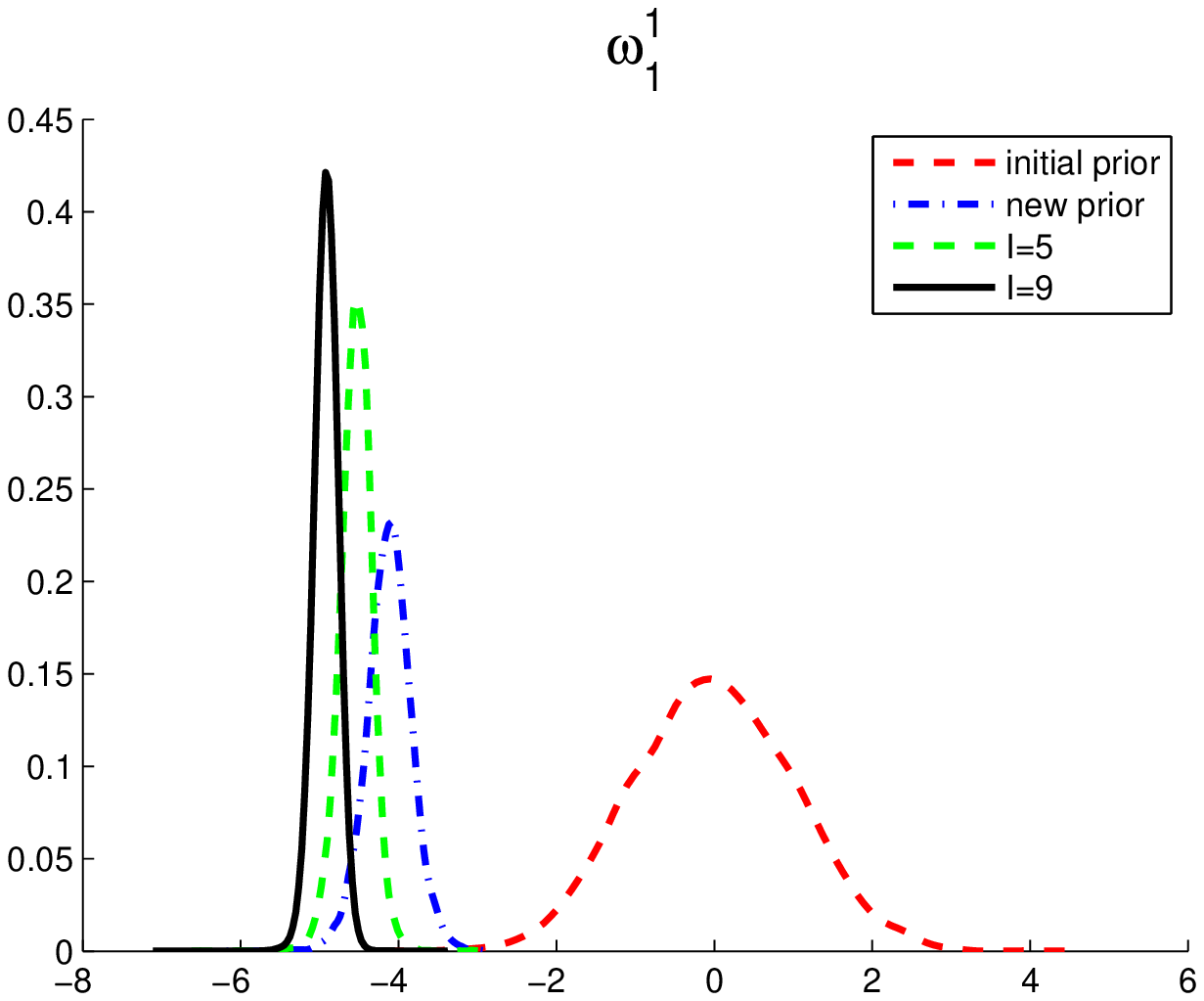}
  \includegraphics[width=2.1in, height=1.7in]{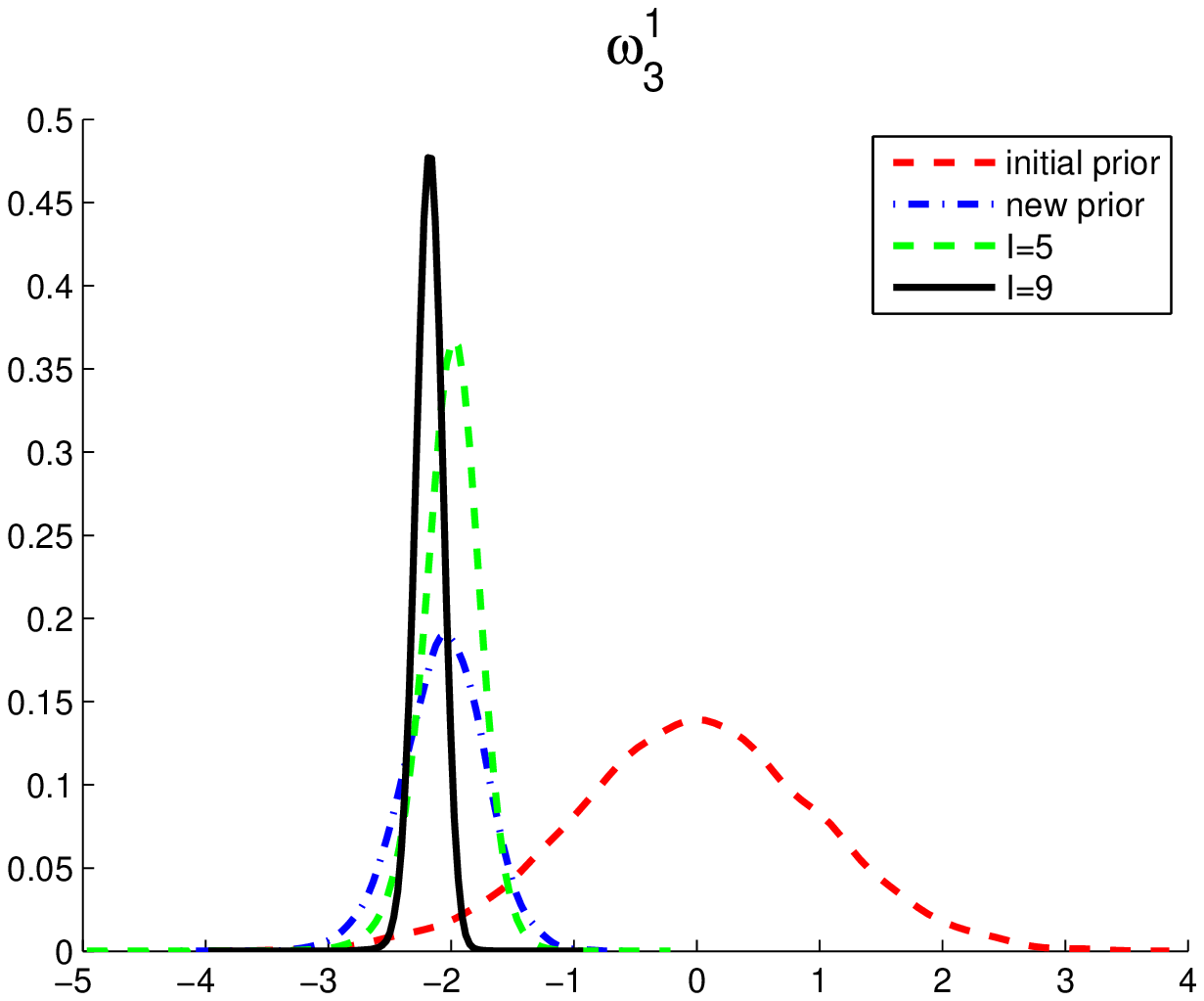}
  \includegraphics[width=2.1in, height=1.7in]{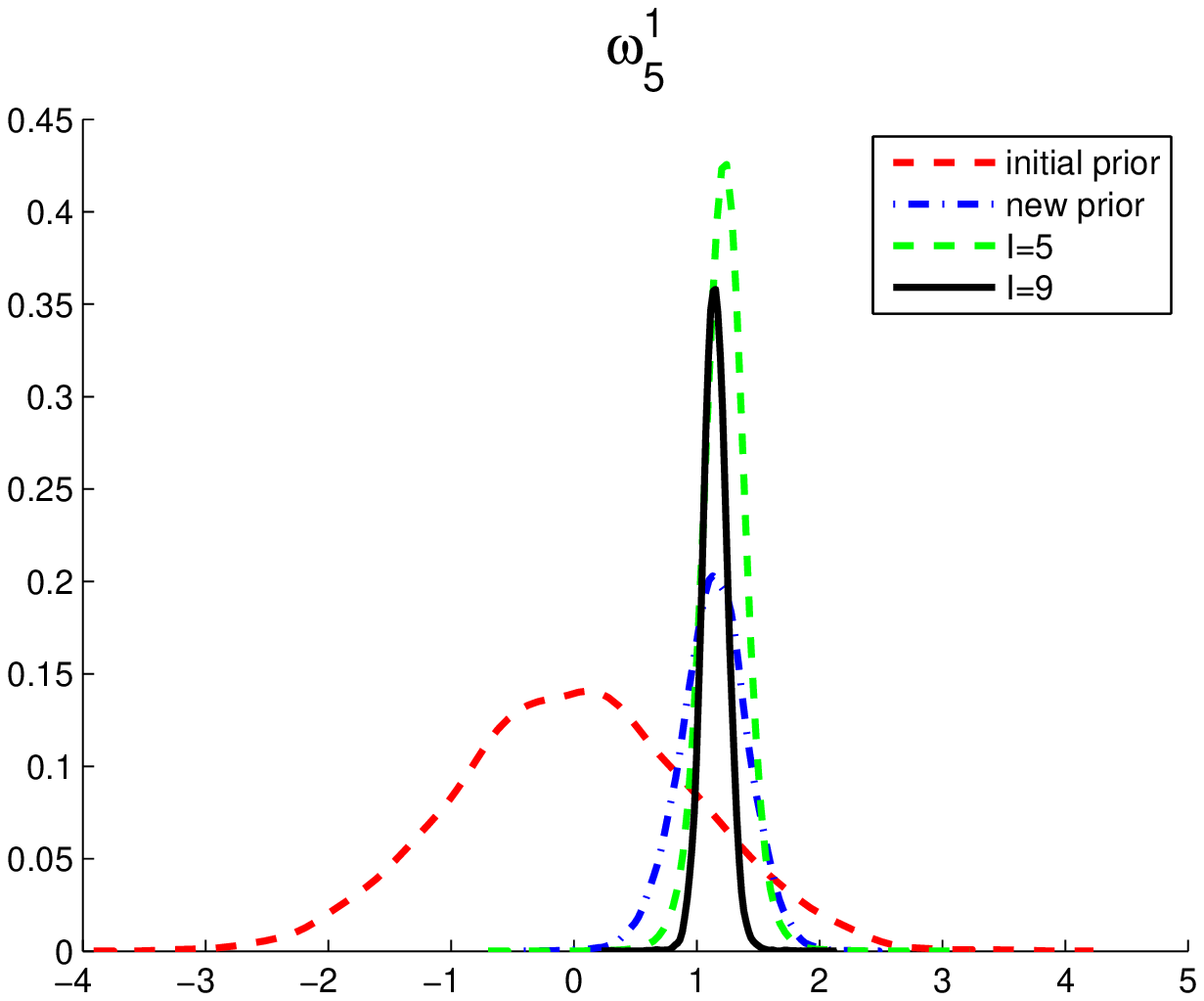}
  \includegraphics[width=2.1in, height=1.7in]{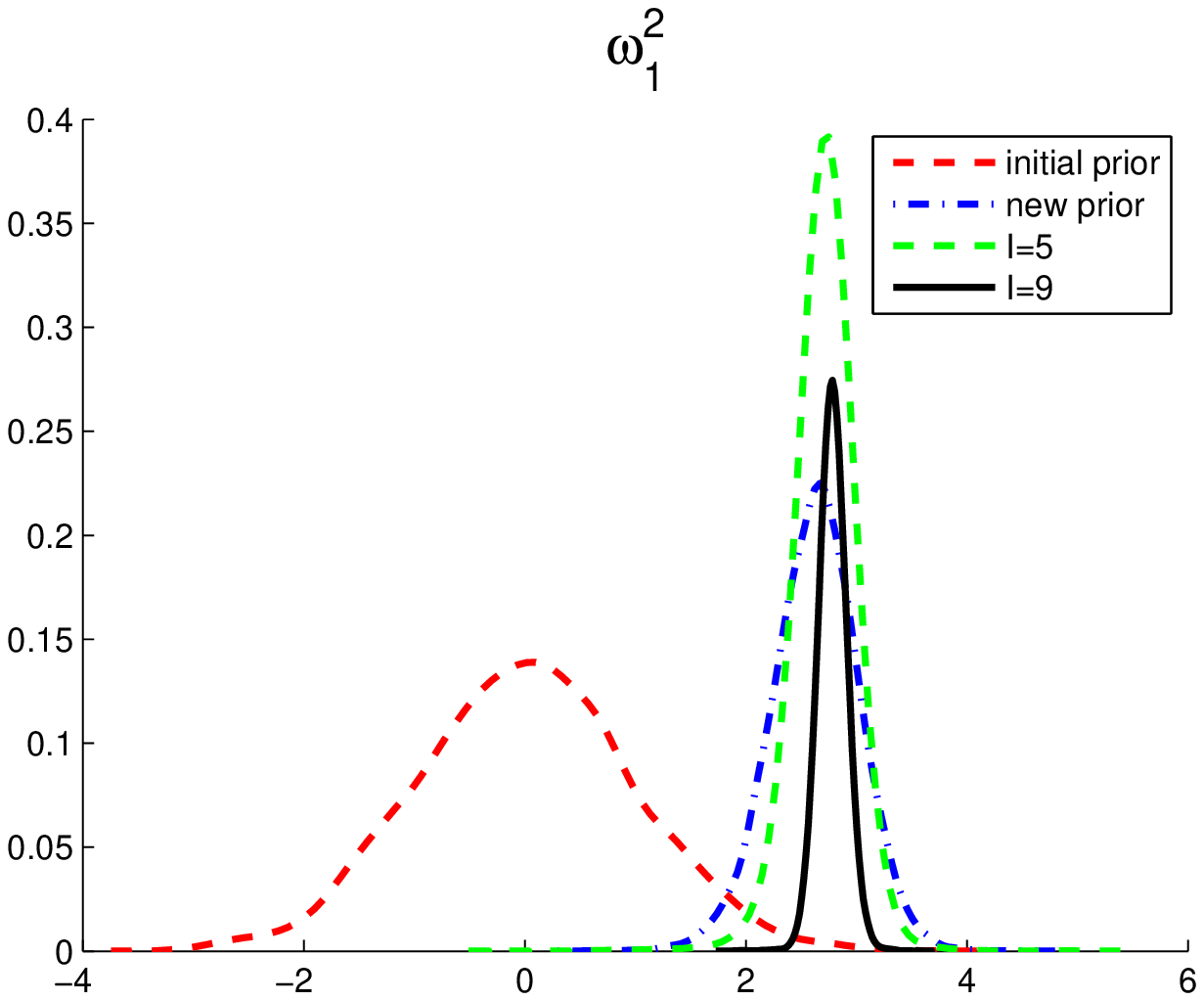}
  \includegraphics[width=2.1in, height=1.7in]{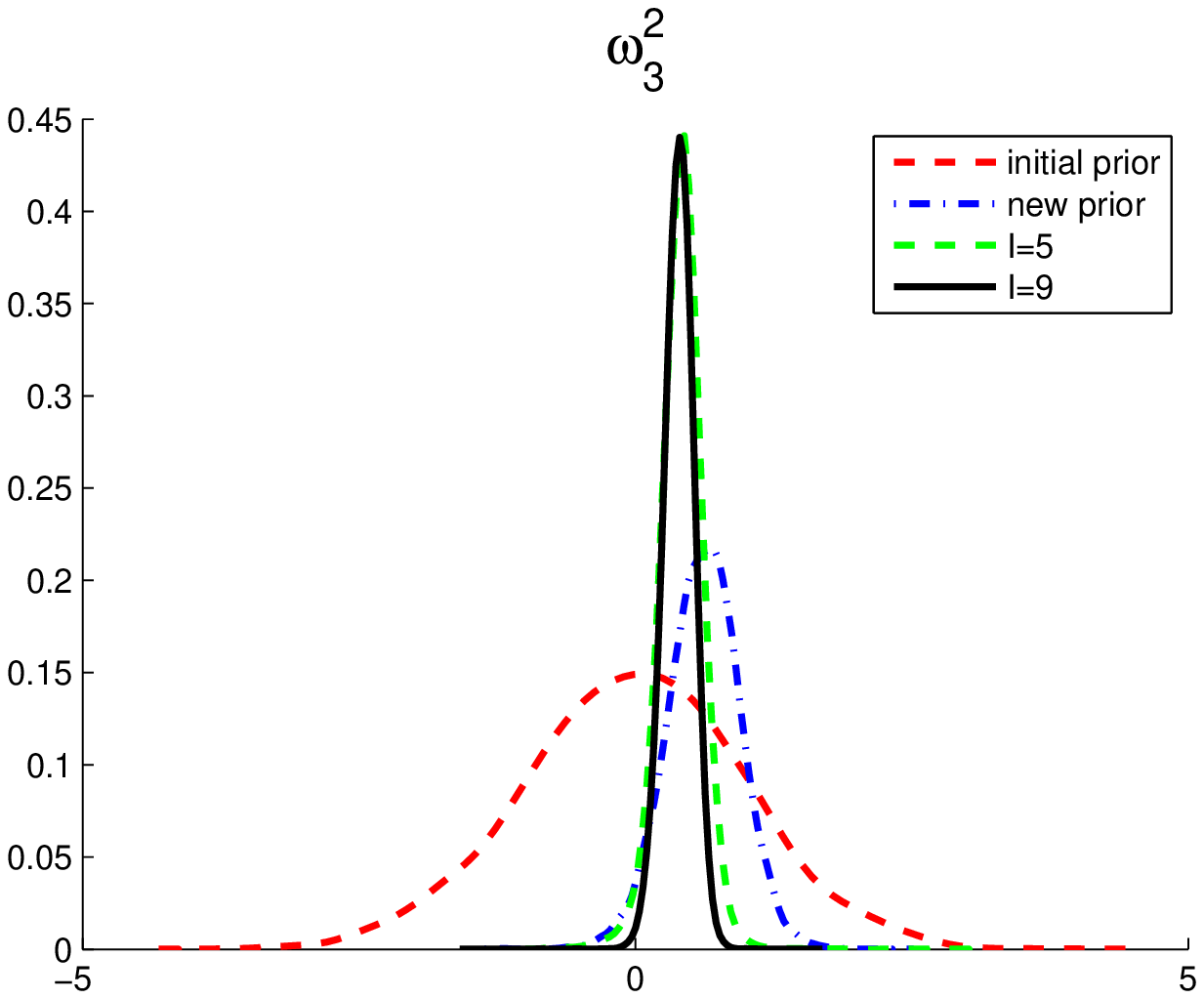}
  \includegraphics[width=2.1in, height=1.7in]{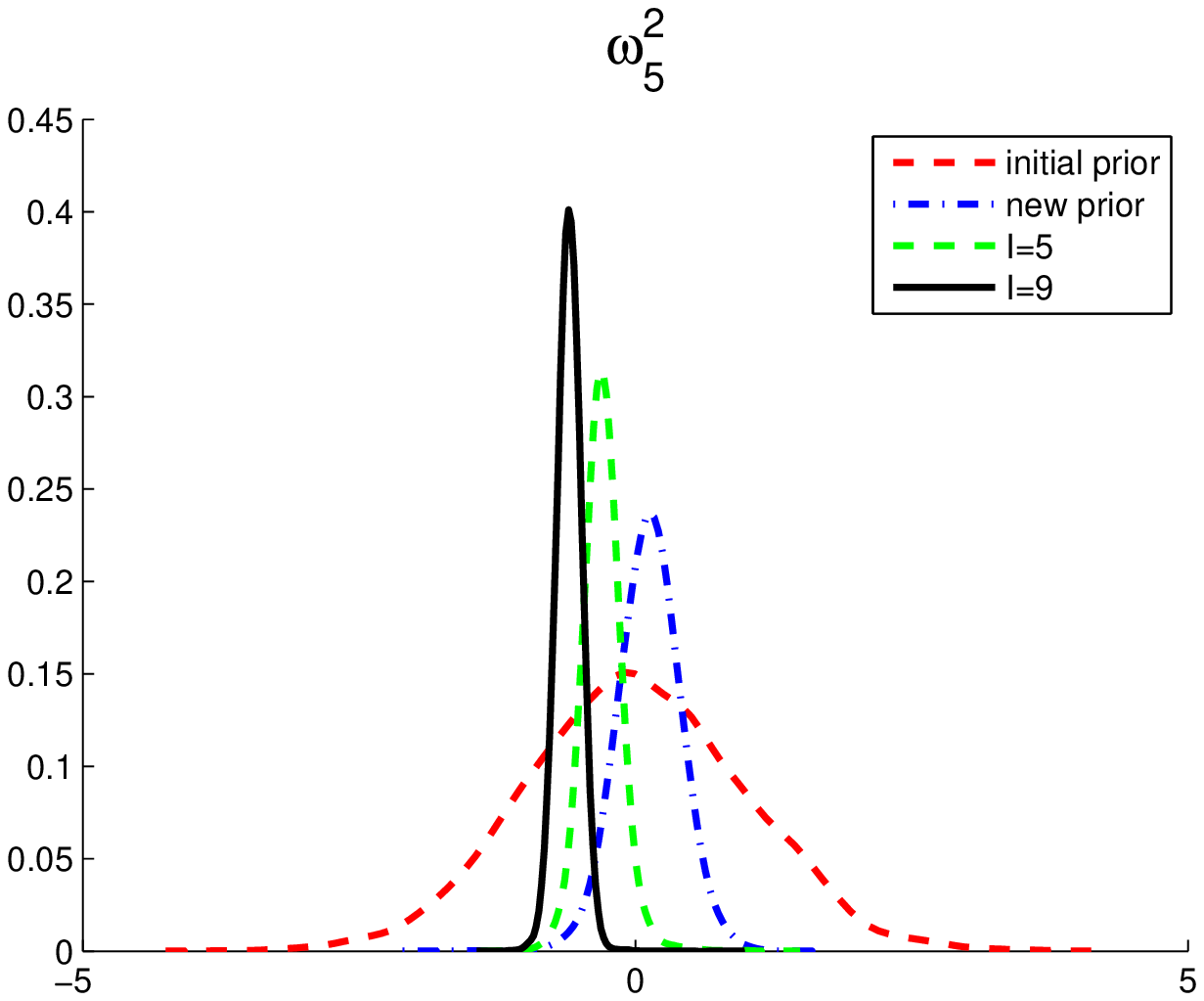}
  \caption{Marginal posterior density estimation for $\omega_j^i$ ($i=1,2$; $j=1,3,5$)  in  different data assimilation steps.}\label{chde}
 \end{figure}

Figure \ref{chde} shows the marginal posterior density against the data assimilation step in the two-stage EnKF. Although there are no truth for these parameters  $\{\omega_i^1\}_{i=1}^5$ and $\{\omega_i^2\}_{i=1}^5$, the data is sufficiently informative to identify a small range of values for the unknown parameters. The important region of marginal densities becomes narrower as data information gains.

\subsection{Estimate  the source locations and fractional derivative}\label{sourcese}
\label{source}
In this section, we consider a diffusion-wave problem, i.e., the diffusion type is unknown, where the boundary conditions are the same as in Subsection \ref{channel}. The end time is set as $T=0.1$, and
the source term is given by
\begin{eqnarray*}
  f(x, t)&=\frac{s_1}{2\pi\tau_1^2}\exp\{-\frac{\|\chi_1-x\|^2}{2\tau_1^2}\}[1-H(t-T_m)] \nonumber\\
   &+\sum_{i=2}^m \frac{s_i}{2\pi\tau_i^2}\exp\{-\frac{\|\chi_i-x\|^2}{2\tau_i^2}\}H(t-T_{i-1}),
\end{eqnarray*}
where $H(t)$ is the heaviside function, $s_i$ is the strength, $\tau_i$ is the width and $\chi_i$ is the source location. In this example,
both $\chi_i$ and $\gamma$ are unknown. We just have the prior $\chi_i=(\chi_i^1, \chi_i^2)\in{[0,1]\times[0,1]}$, and $\gamma \in (0,2)$ and $\gamma\neq 1$. The prior is bounded, which can be seen as uniform distribution and $\gamma$  lies  in the  interval $(0,2)$ for diffusion-wave equation. But its values are  uncontrollable during the data assimilation steps, the samples may run out of the  interval. Thus, the two-stage EnKF by normal-score method can be used to avoid the issue. To this end, we use a bijective map $\bb{F}:E \rightarrow{\bb{R}}^5$, where
\[
E=(0,2)\times (0,1)^4.
\]
Let $\bt=(\gamma,\chi_1,\chi_2)^{\rm T}$ and $\bb{F}({\bt})=\bm{q}$, the map $\bb{F}$ is set as
\[
\begin{aligned}
&q_i=\bb{F}_i(\theta_i),\quad i=1,\cdots,5.
\end{aligned}
\]
The dynamic system correspondingly becomes
\[
 \left\{
 \begin{aligned}
 \bm{q}_{k-1} &={\bb F}_{k-1}(\bt_{k-1}), \\
 \bt_k &=\bt_{k-1},\\
  \bm{q}_k &=\bm{q}_{k-1}, \\
  \bm{y}_k &=\bb{H}({\bt}_{k-1})+\bm{\varepsilon}_k,
  \end{aligned}
\right.
\]
where $\bm{q}_k=(q^k_1,\cdots,q^k_5)^{\rm T}$ and $k$ is the artificial time for data assimilation, and the Kalman gain is approximated by $\bm{K}_k=\text{Cov}(\bm{q}_k,\bm{y}_k)\text{Cov}(\bm{y}_k,\bm{y}_k)^{-1}$.

\begin{figure}[htbp]
  \centering
  \includegraphics[width=2.5in, height=1.7in]{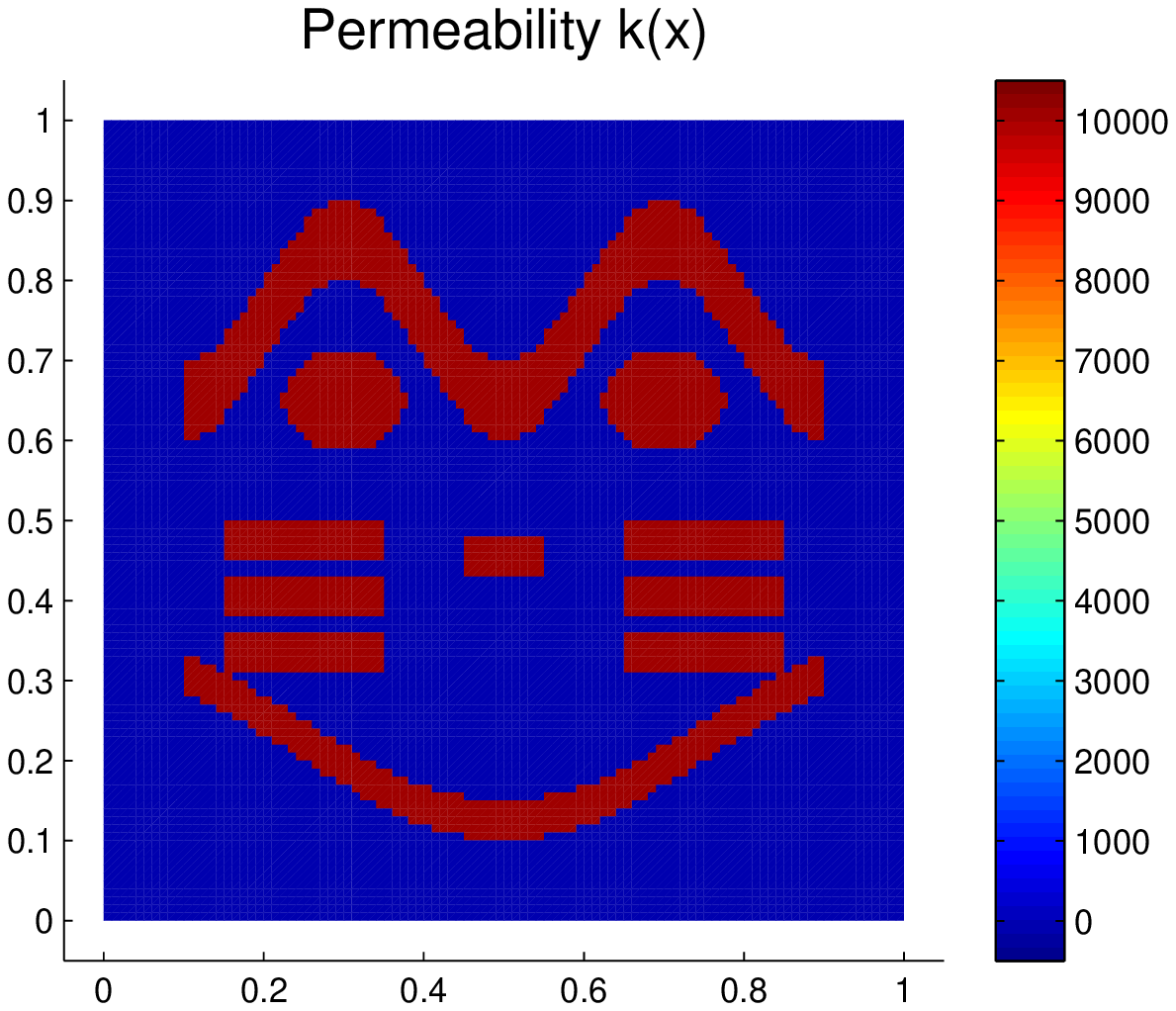}
  \includegraphics[width=2.5in, height=1.7in]{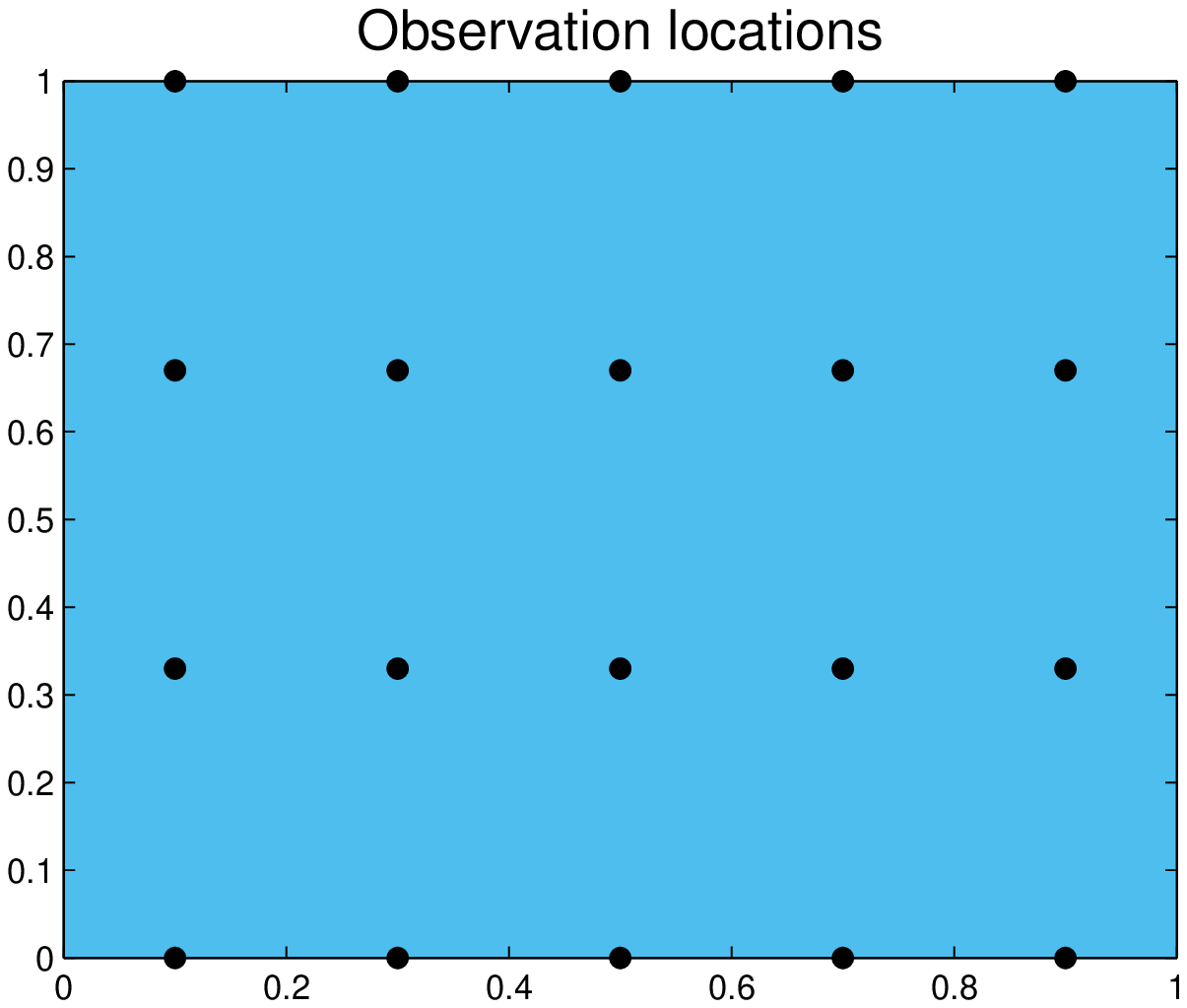}
  \caption{Spatial distribution of the permeability $k(x)$ (left) and observation locations (right).}\label{sourcehigh}
 \end{figure}

For the source term, we take $m=2$, $T_1=0.05$, $T_2=T$, and $\tau_i=0.1$ ($i=1,2$), $s_1=3$, $s_2=1$, and the truth source locations are set as $\chi^{tr}_1=(0.2,0.6)$, $\chi^{tr}_2=(0.5,0.3)$ and $\gamma^{tr}=0.5$. Here, the total artificial time steps for data assimilation is set as $I_2=8$. Measurements are taken from the subdiffusion equation at time instances $0.012+0.01I:0.006:0.018+0.01I$ in each data assimilation step, where $I\in\{1,2,\cdots,8\}$ and the locations are distributed on the uniform $5\times 4$ grid of the domain $[0.1, 0.9]\times [0, 1]$ as shown in Figure \ref{sourcehigh} (right). The high-contrast permeability is known and has the spatial distribution shown in  Figure \ref{sourcehigh} (left).

The forward model is defined on $100\times100$ uniform fine grid, and GMsFEM is implemented on $5\times5$ coarse grid.
 We choose $M_{\text{snap}}=10$ eigenfunctions for the local snapshot space  and use  GMsFEM with $9$ local multiscale basis functions  to approximate the forward model in solving the optimization problem (\ref{cos}).

In the subsection, we focus on the two-stage EnKF for the inverse problem.  To this end,  we use $3\times 10^3$  ensemble members to construct a new prior by standard EnKF in the first stage. The new prior is constructed by incorporating data information from the  first data assimilation step. $9$ local multiscale basis functions ($M_i=9$) are chosen  in constructing the matrix $R$ to construct coarse model. Then, we construct the surrogate model after obtaining the new prior ensemble in the second stage. The data assimilation process begins from $I=1$ in second stage. The same local multiscale basis functions are used  in constructing the matrix $R$ to construct the gPC surrogate model. When the order of gPC is set as $N_0=7$, the number of random samplers is 792 when computing vectors $\bf b$ and $\bf A$ in Section \ref{SCM}. In the second  stage of the proposed EnKF, the number of ensemble members is also set as $3\times 10^3$.

\begin{figure}[htbp]
  \centering
  \includegraphics[width=3in, height=2.5in]{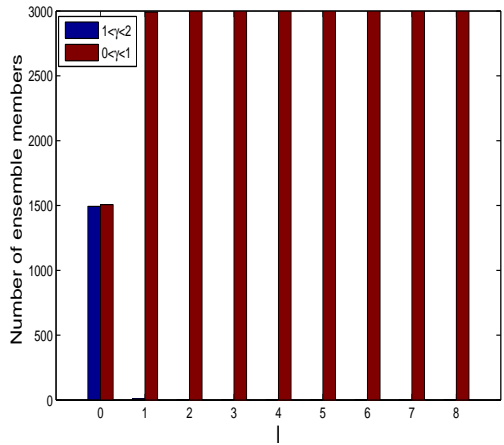}
  \caption{ Number of ensemble members applied to different forward models against data assimilation steps,  where $I=0$ is the process of generating a new prior.}\label{sourceModel}
\end{figure}

\begin{figure}[htbp]
  \centering
  \includegraphics[width=4in, height=3in]{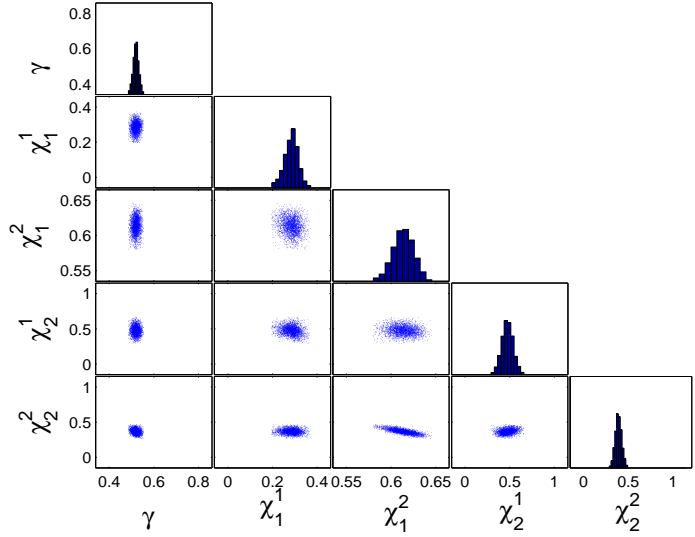}
  \caption{1-D and 2-D posterior marginals of $\bt$}\label{sourcecor}
 \end{figure}

We only know the problem is a diffusion-wave equation, which may be  the subdiffusion or  superdiffusion. Thus, it is important  to identify  the forward model before estimating   the source locations.
The result of identifying model is shown in Figure \ref{sourceModel}, where we count the number of ensemble members applied to the two diffusion models.
By the figure, we can find that the new prior has actually identified  the forward model, where $I=0$ denotes the process of generating a new prior. In the second stage of two-stage EnKF, the superdiffusion model is almost not used. We define the relative error $\varepsilon$ corresponding to the posterior distribution by
\begin{equation}\label{relio}
  \varepsilon := \frac{\|\bb{E}(\bt)-\bt^{tr}\|}{\|\bt^{tr}\|},
\end{equation}
where $\bt^{tr}=(\gamma^{tr}, \chi_1^{tr}, \chi_2^{tr})$. Then we compute  the relative error of the final ensemble mean and get  $\varepsilon=0.08$.
Figure \ref{sourcecor} shows all of the one and two-dimensional posterior marginal of $\bt$,  where there exist some correlation  between the source locations and $\gamma$, such as $\chi_1^2$ and $\chi_2^2$, but the others appear  uncorrelated and mutually independent based on the shape of their 2-D marginal.

\begin{figure}[htbp]
 \centering
  \includegraphics[width=1.5in, height=1.4in]{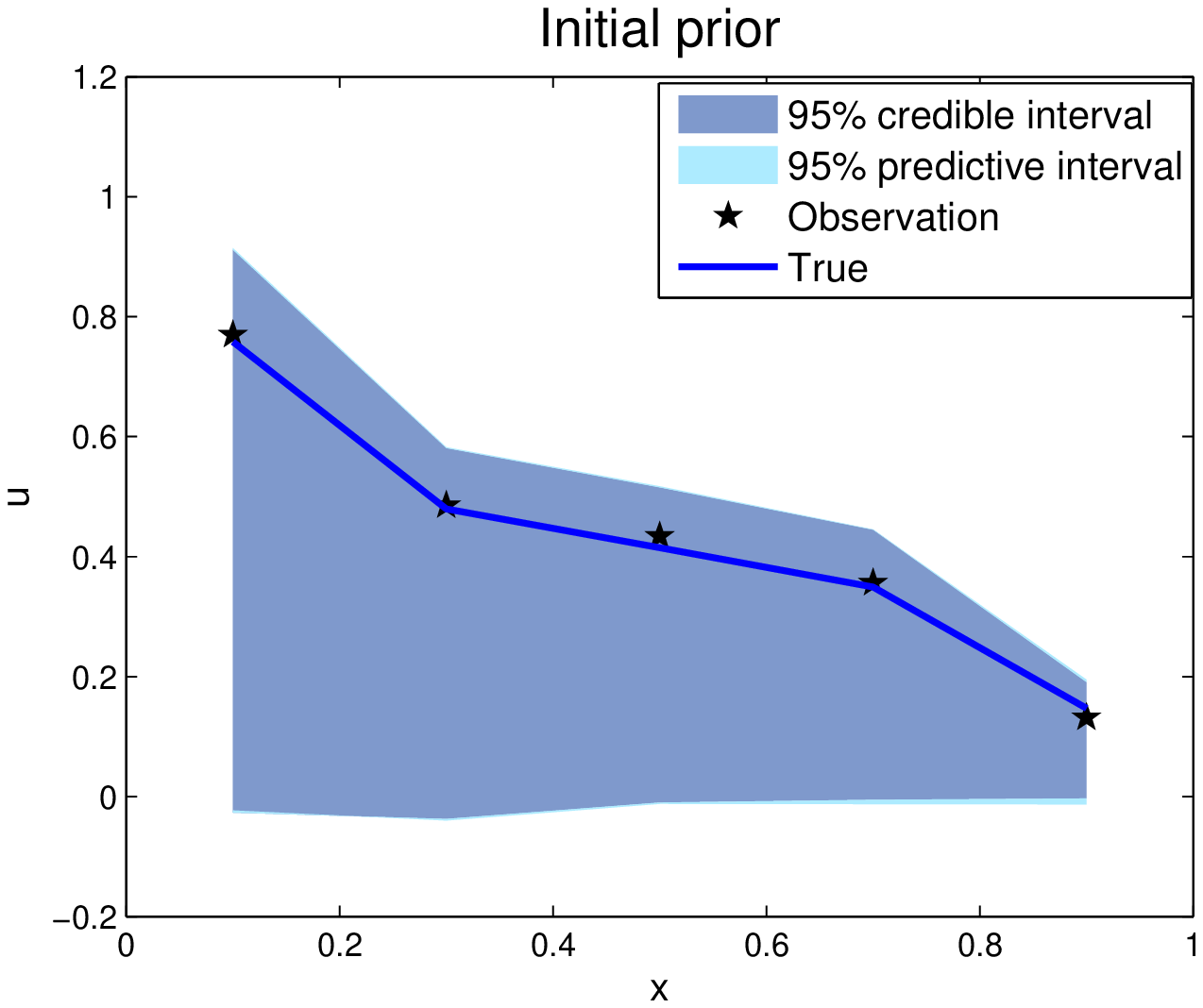}
  \includegraphics[width=1.5in, height=1.4in]{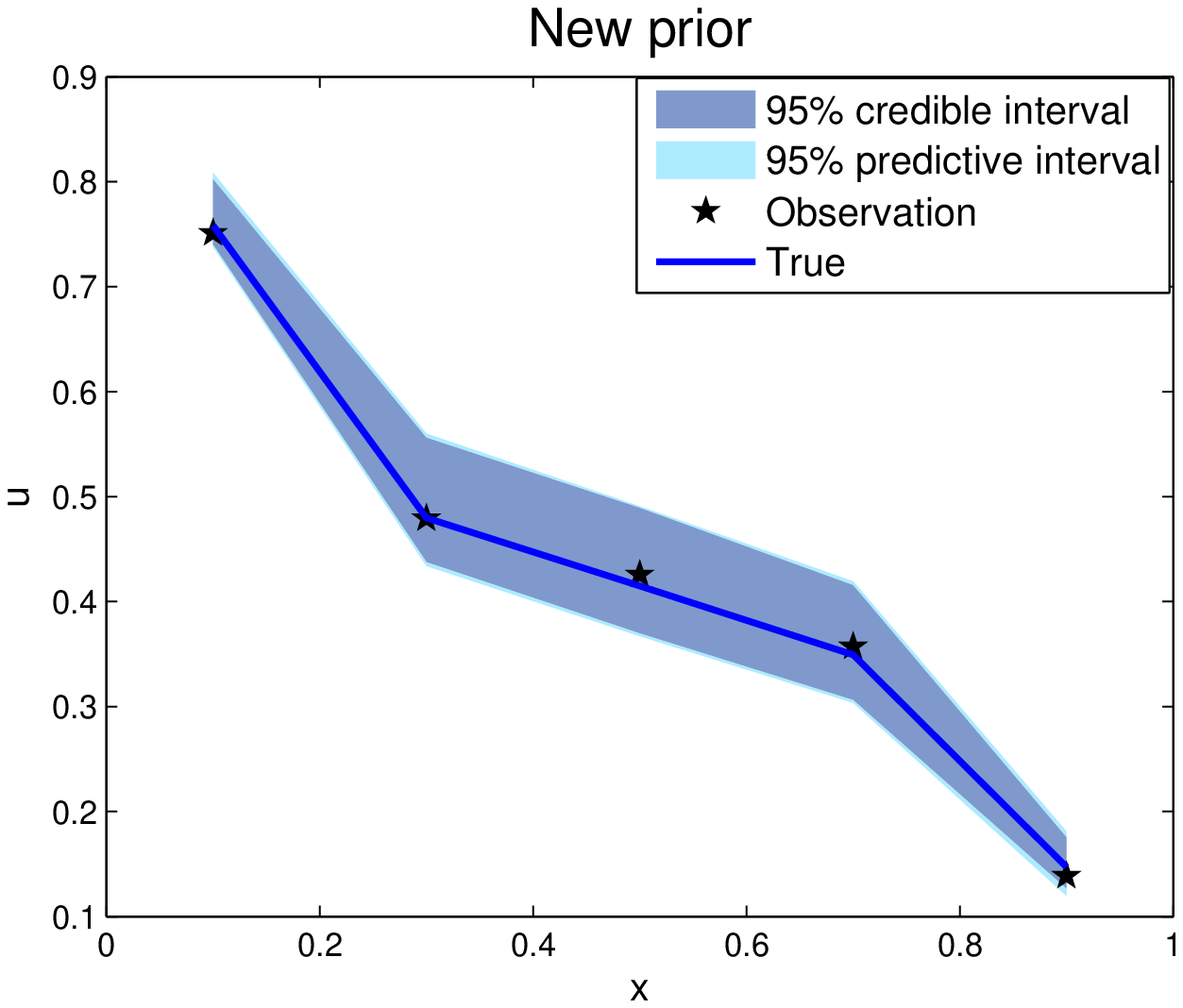}
  \includegraphics[width=1.5in, height=1.4in]{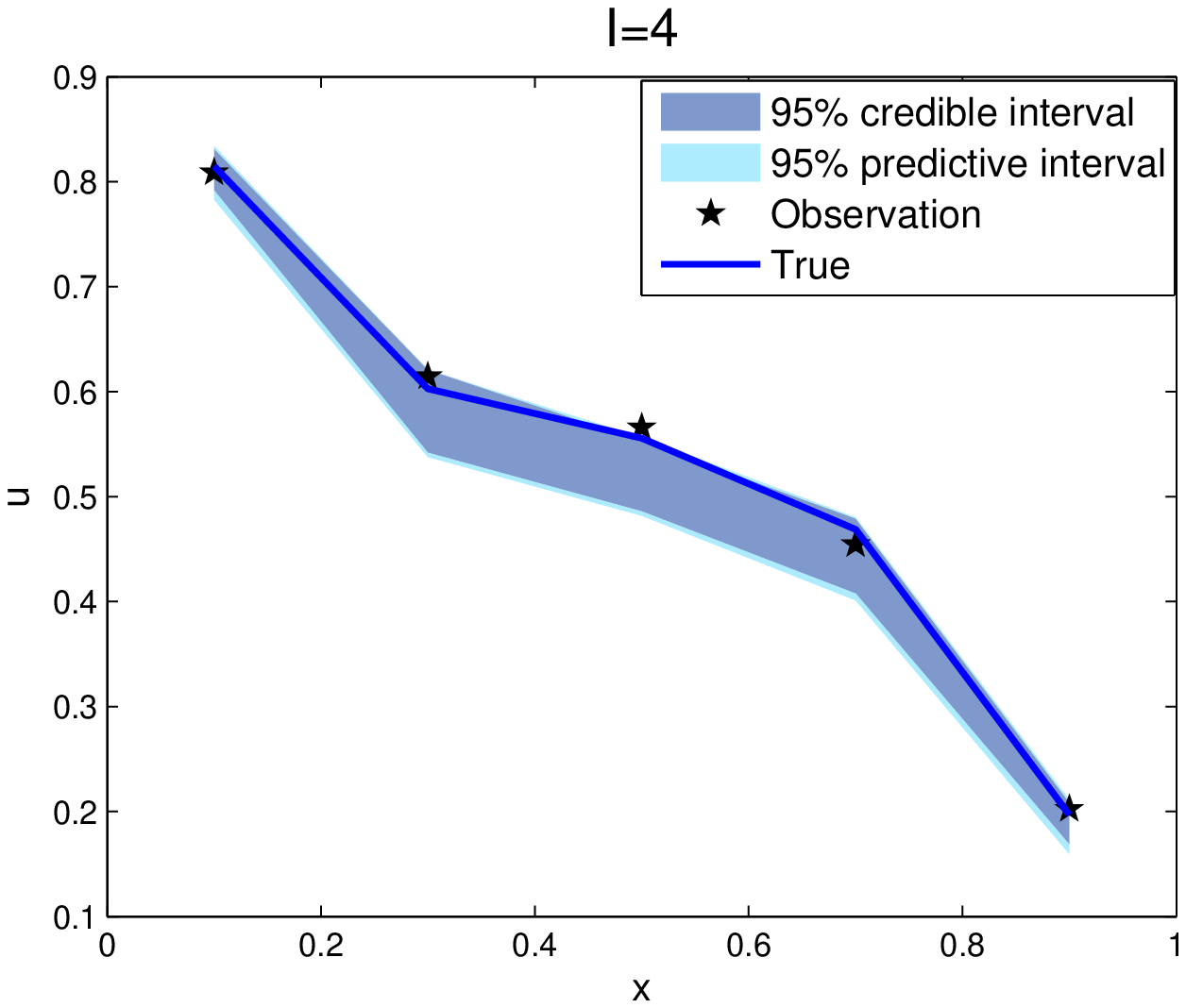}
  \includegraphics[width=1.5in, height=1.4in]{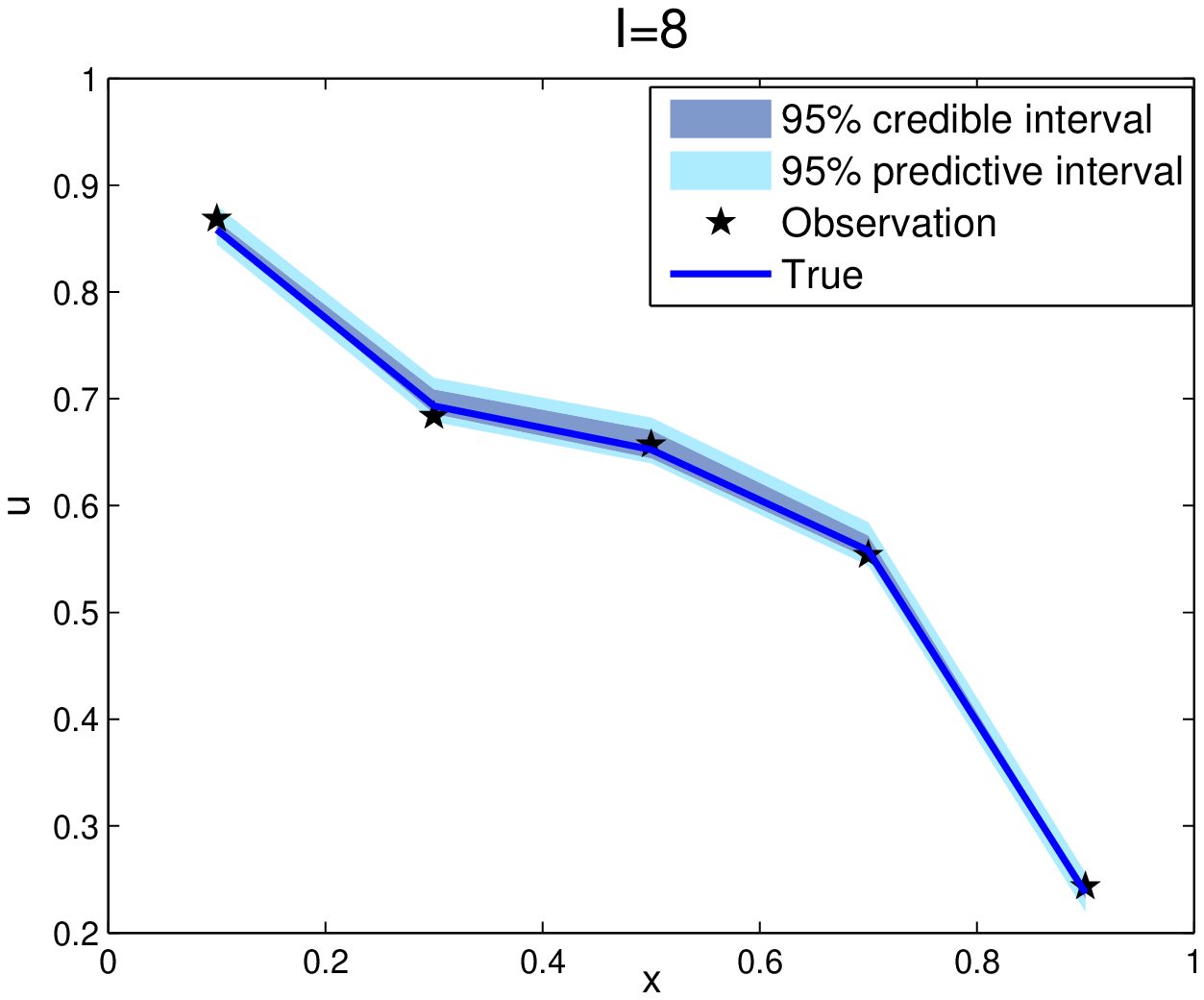}
  \caption{95\% predictive interval, 95\% credible interval, observation and true value for $u((x, 0); t)$.}\label{sopredict1}
\end{figure}

\begin{figure}[htbp]
  \centering
  \includegraphics[width=1.5in, height=1.4in]{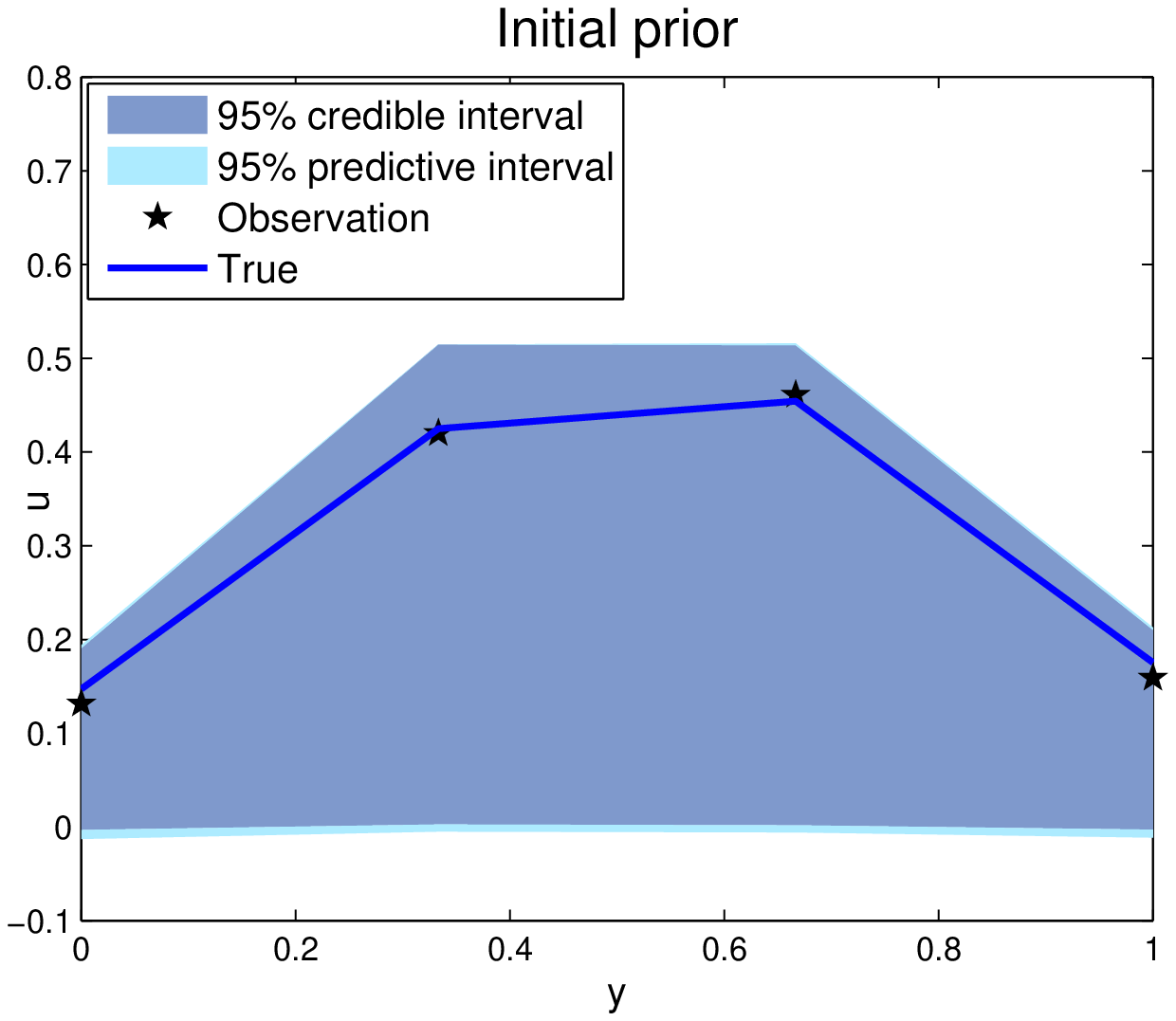}
  \includegraphics[width=1.5in, height=1.4in]{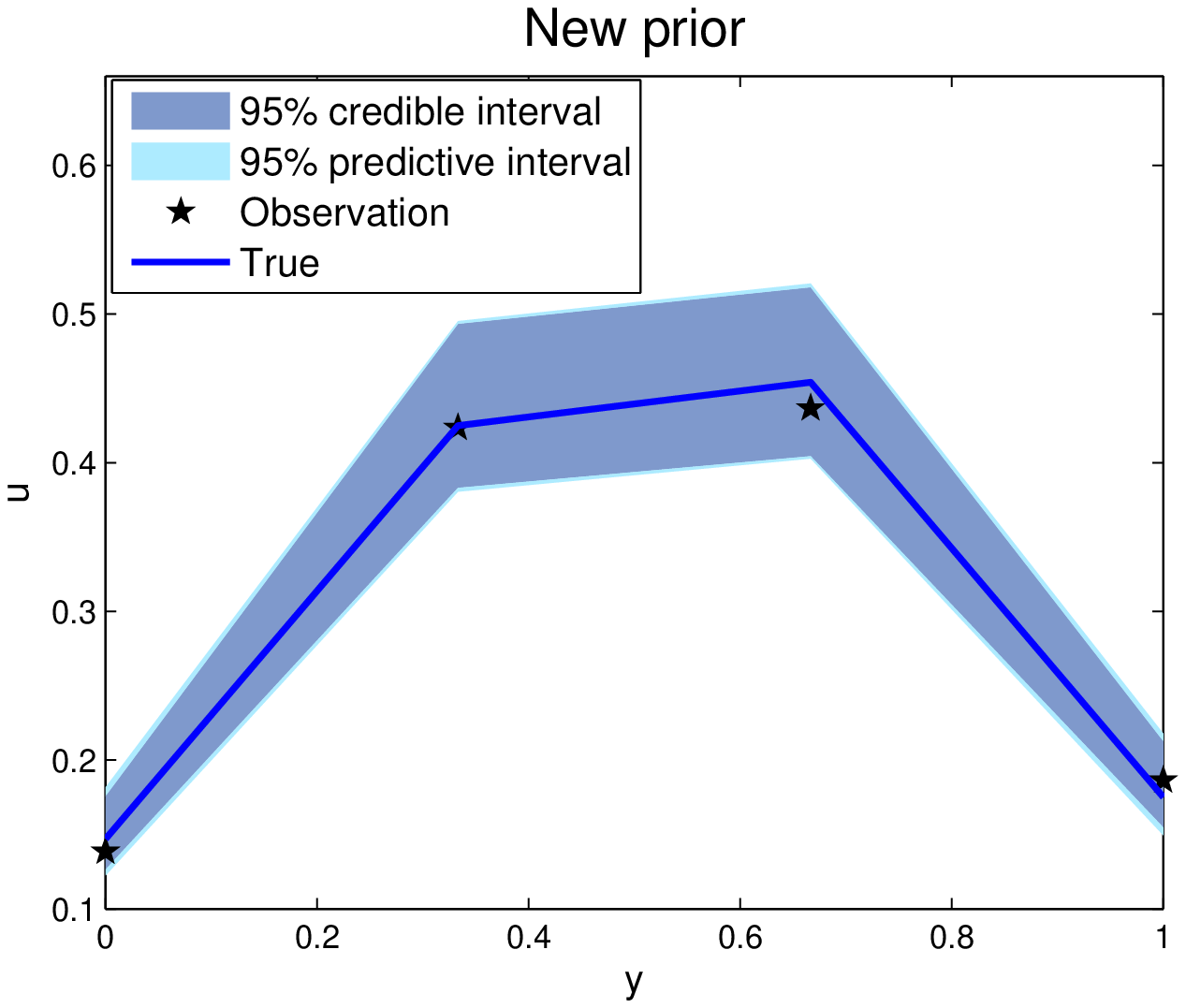}
  \includegraphics[width=1.5in, height=1.4in]{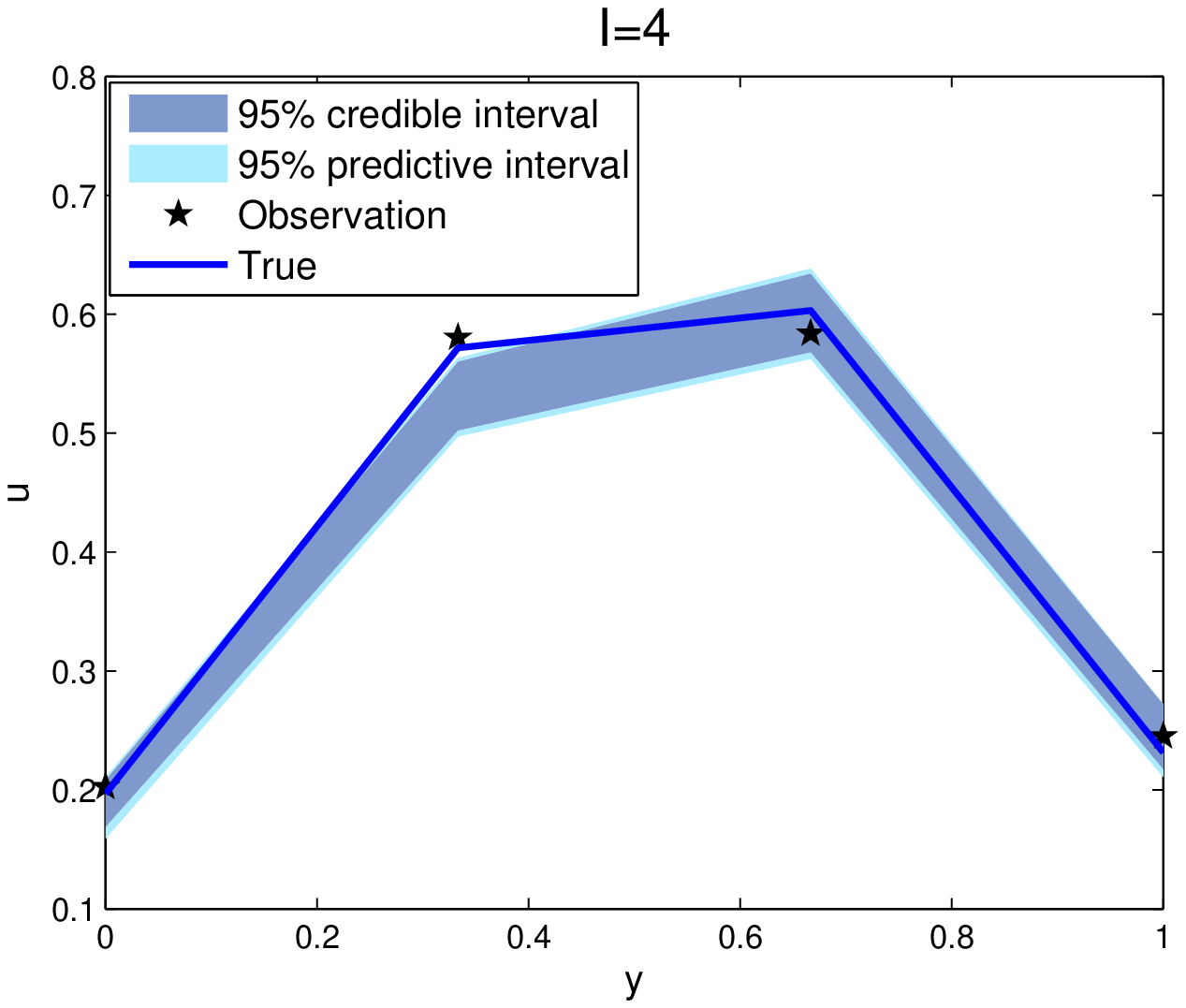}
  \includegraphics[width=1.5in, height=1.4in]{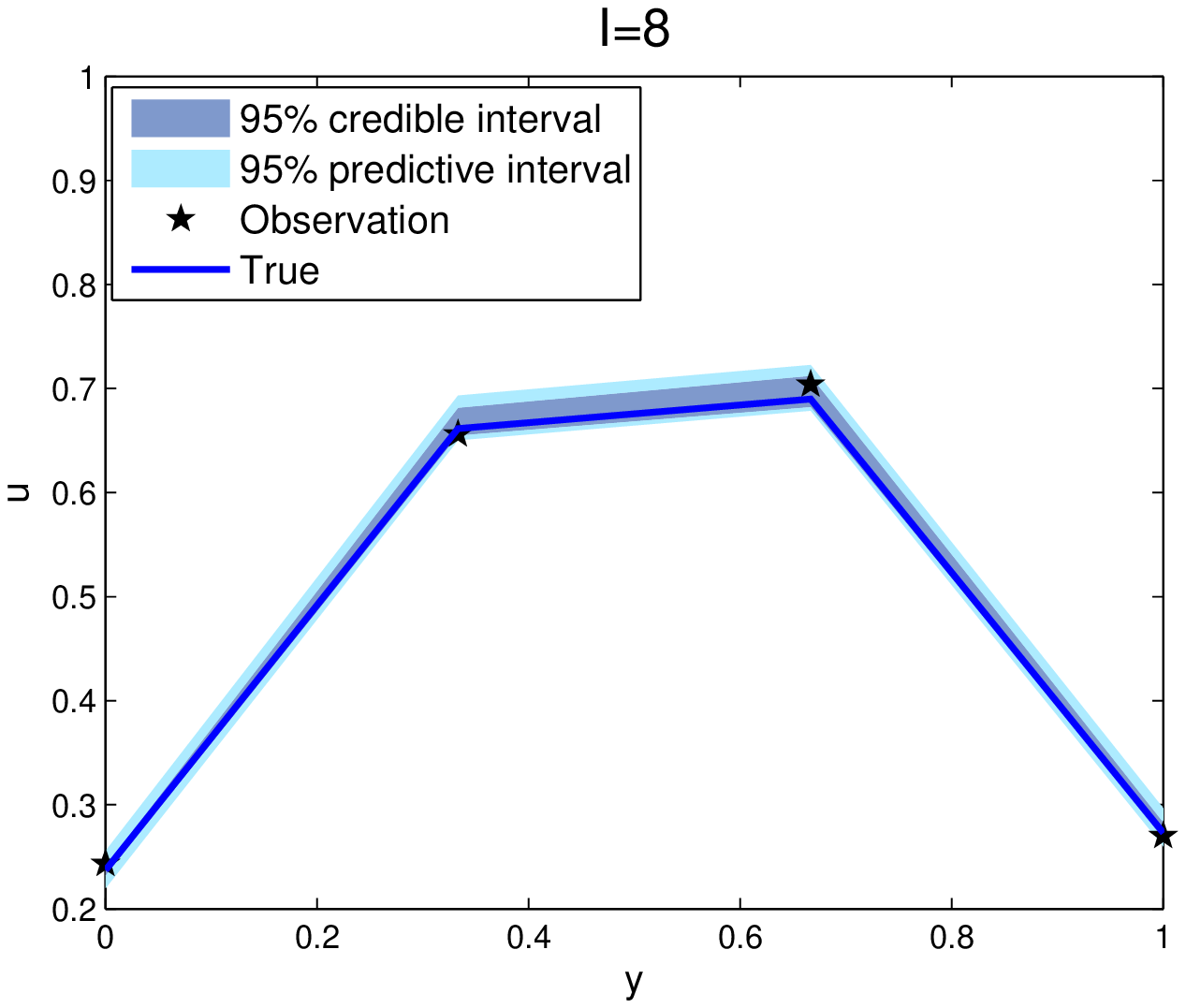}
  \caption{95\% predictive interval, 95\% credible interval, observation and true value for $u((0.9, y); t)$.}\label{sopredict2}
\end{figure}

\begin{figure}[htbp]
  \centering
  \includegraphics[width=2.1in, height=1.7in]{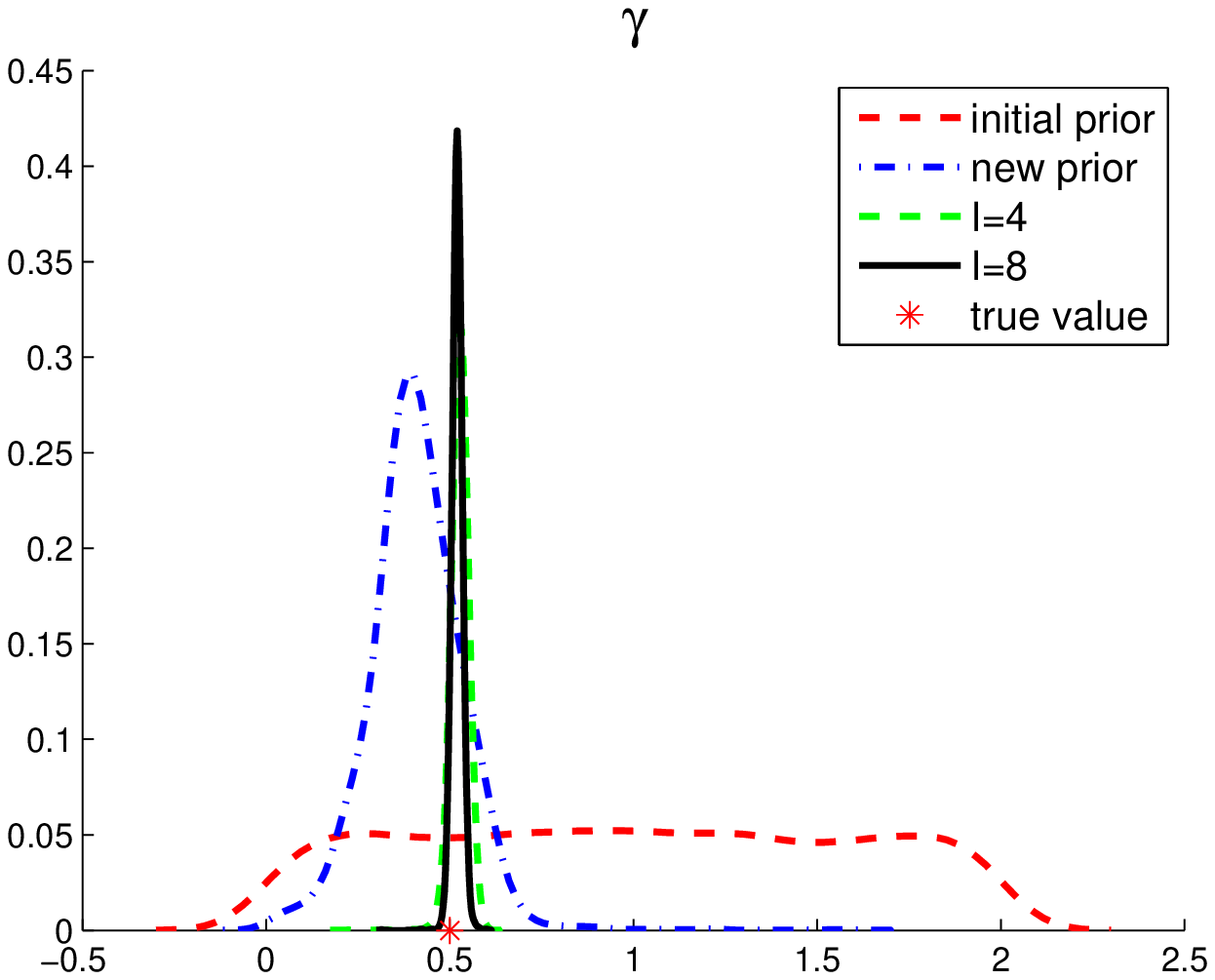}
  \includegraphics[width=2.1in, height=1.7in]{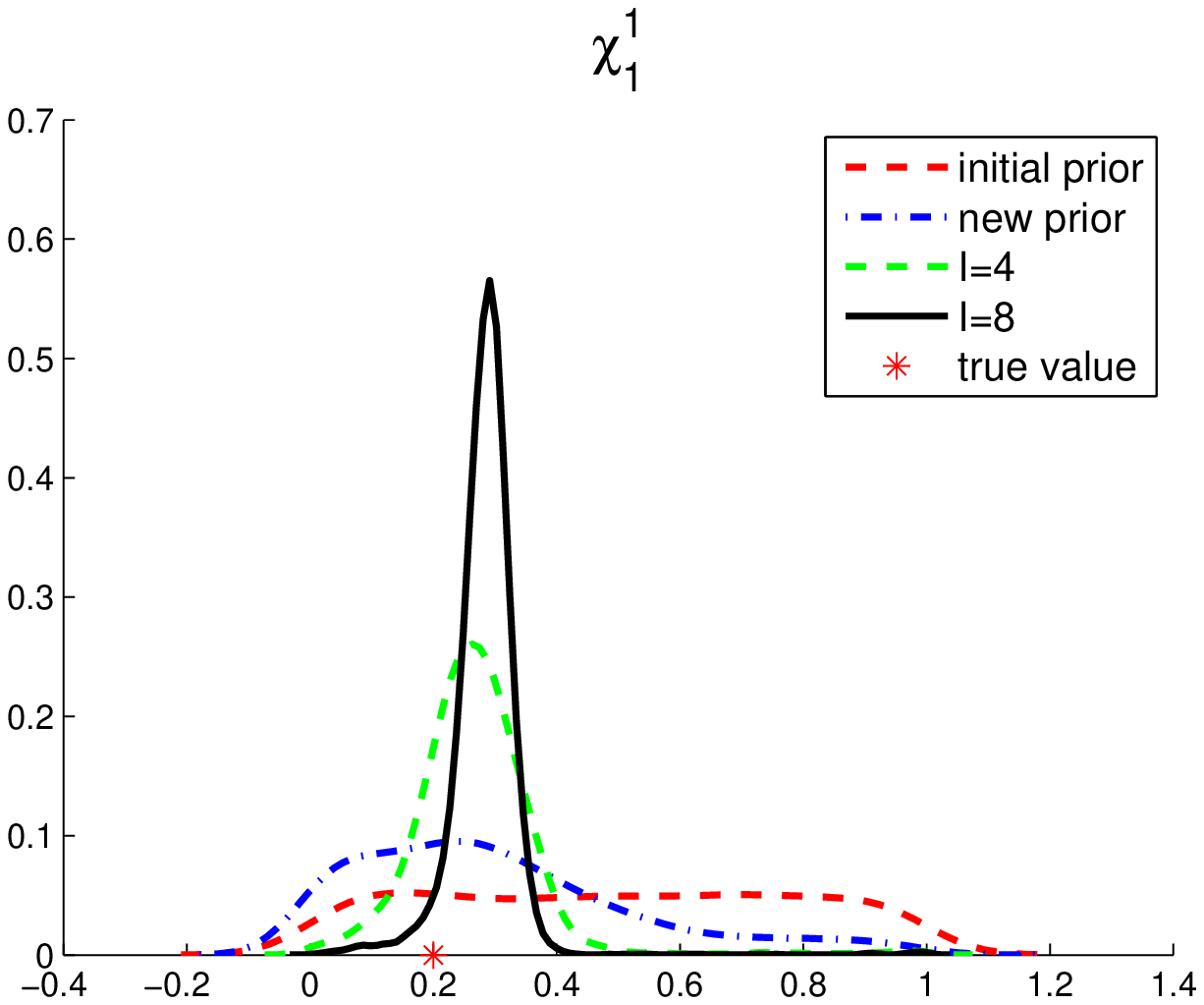}
  \includegraphics[width=2.1in, height=1.7in]{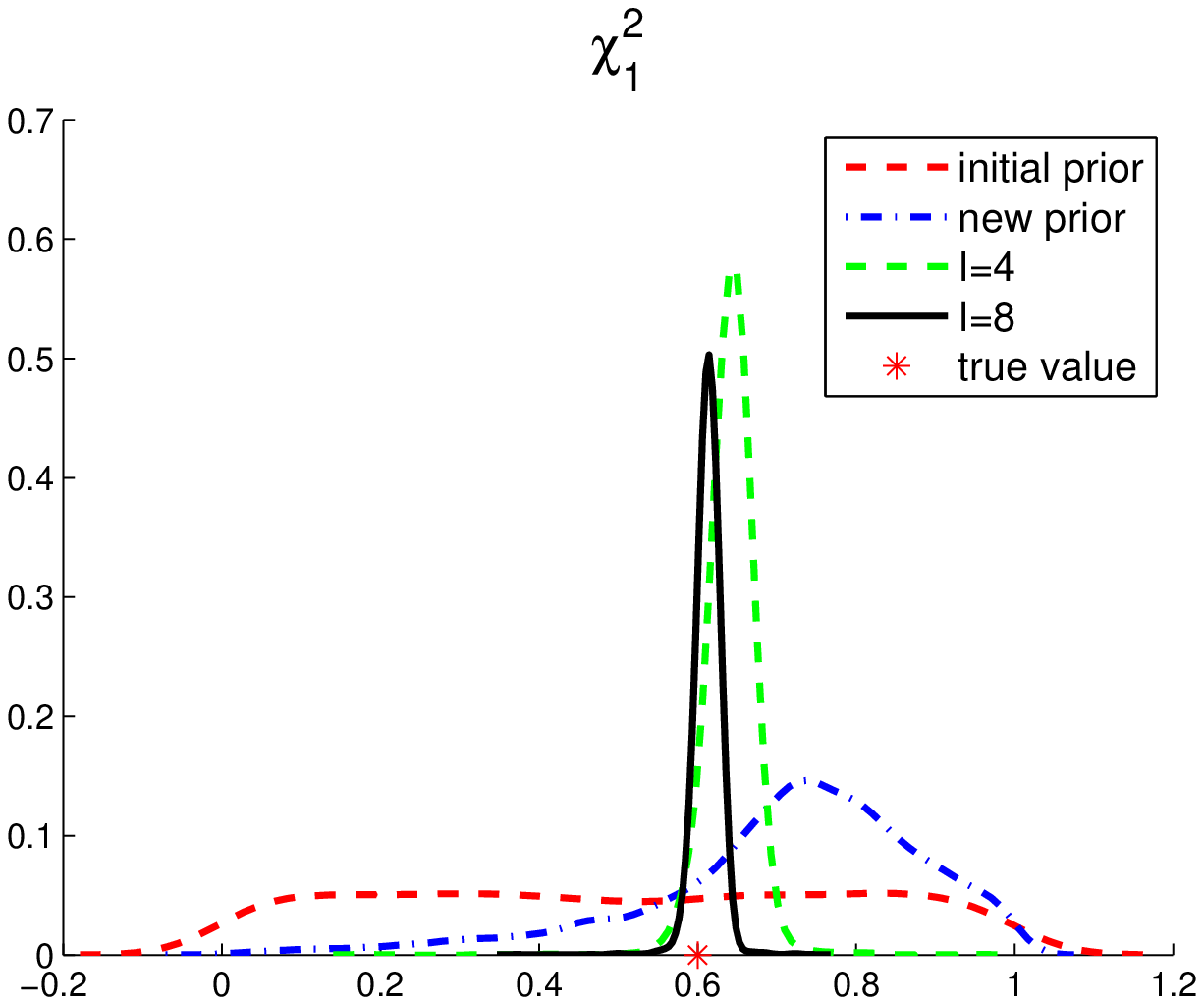}
  \includegraphics[width=2.1in, height=1.7in]{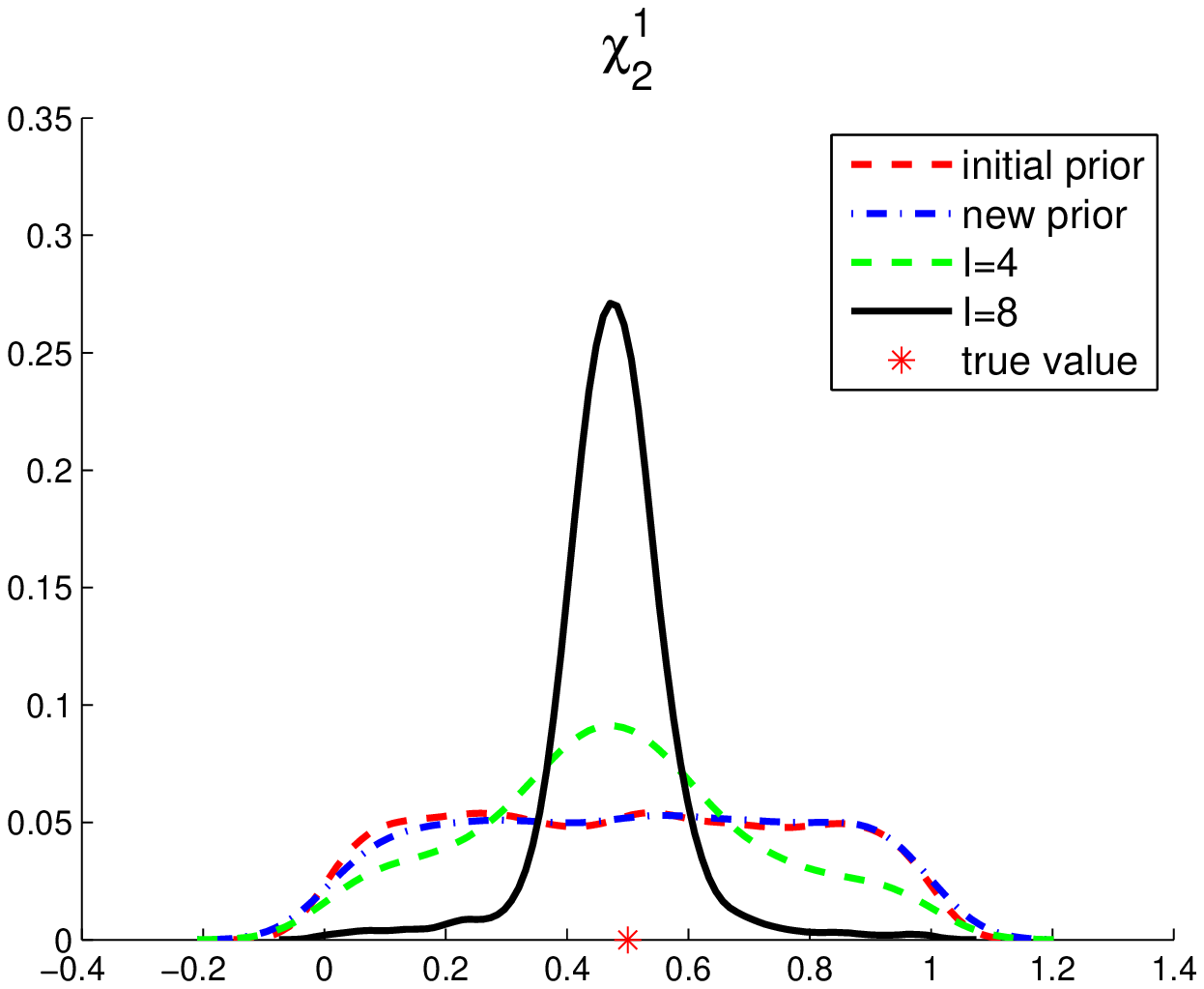}
  \includegraphics[width=2.1in, height=1.7in]{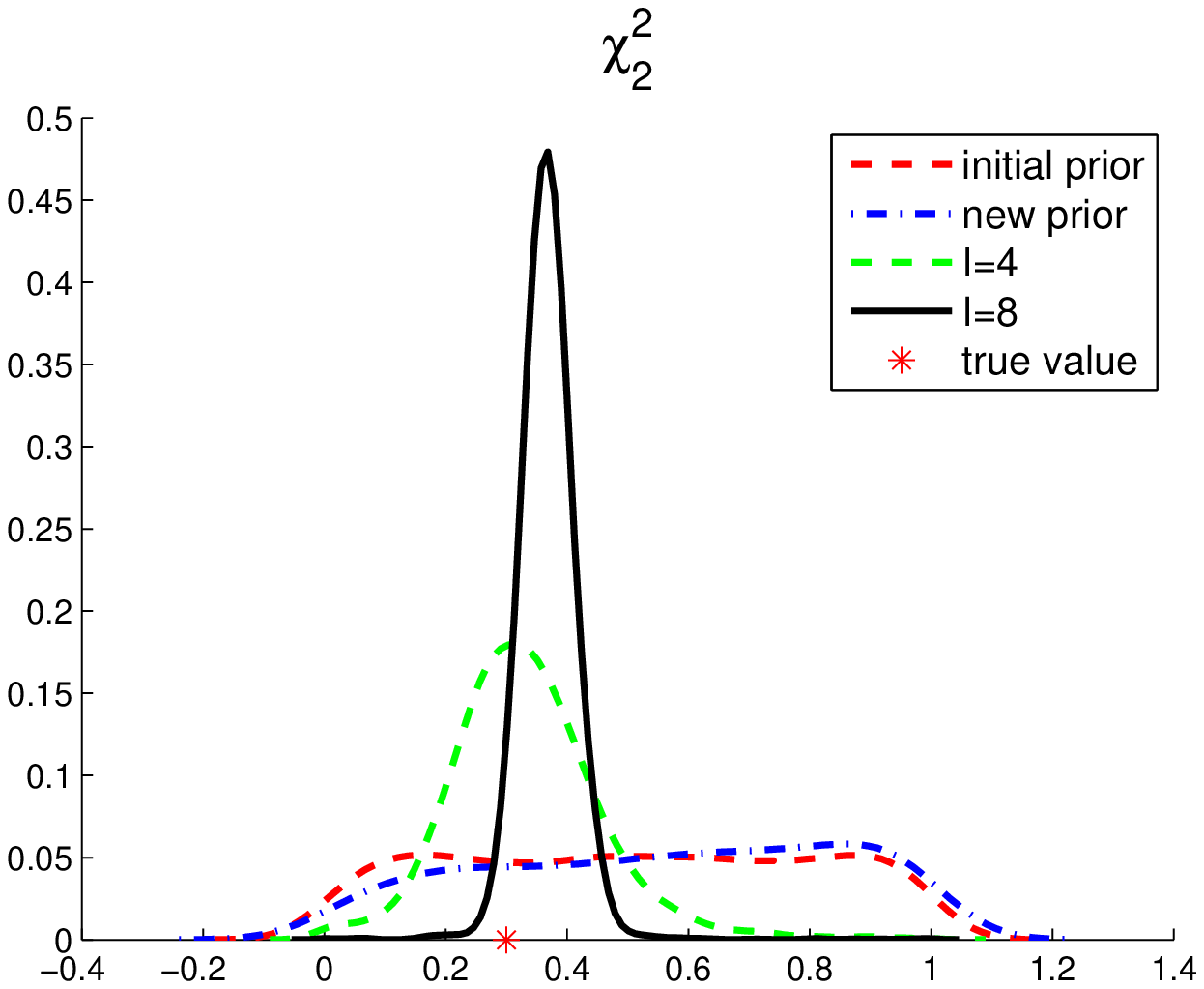}
  \caption{ Marginal posterior density estimation of $\bt$ at different data assimilation steps.}\label{sodensity}
\end{figure}

 To make prediction, we plot the 95\% credible interval and predictive  for the model response with the ensemble members and realizations  by the surrogate model  in Figure \ref{sopredict1} and Figure \ref{sopredict2}. We note the realizations of model response for initial prior is constructed by GMsFE model. The prediction intervals are constructed by taking account of measurement errors. The credible interval and prediction interval, along with the truth, and observation data, are illustrated for $u\big((x, 0); t\big)$ and $u\big((0.9, y); t\big)$ in Figure \ref{sopredict1} and Figure \ref{sopredict2}, respectively. We note that  the credible intervals are tight in the final data assimilation step. Then we can observe from these figures that the uncertainty associated with both the model fit and predictions decreases with respect to assimilation step. We also plot the marginal density for  the initial prior, new prior and posteriors together in  Figure \ref{sodensity}, which clearly shows that the posterior is dynamically  updated. From this figure, we note that the new prior identifies the model belonging to subdiffusion.

\subsection{Estimate  a permeability field in a hierarchical model}\label{coefsc}

\begin{figure}[tbp]
  \centering
  \includegraphics[width=2.5in, height=1.7in]{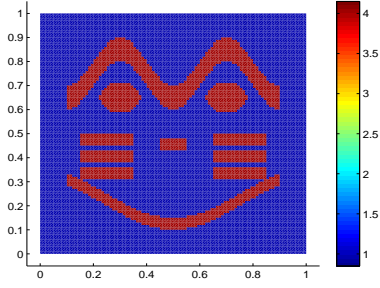}
  \caption{ The spatial distribution of ${\bb E}[\log k (x,\omega)]$.}\label{coh}
 \end{figure}

In this subsection, we consider the fractional superdiffusion equation with homogeneous Dirichlet boundary condition.
The source term is set as $f=20$, and the end time is set as $T=0.11$. In the inversion model, both the permeability filed and observation noise are unknown and we want to estimate them
by the two-stage EnKF.  To parameterize the permeability field,  we set $l_{x}=0.2$, $l_{y}=0.3$, $\sigma^2=1$, and the covariance function is assumed to be
  \[
  C(x,y)=\sigma^2\exp\bigg(-\frac{\|x_1-x_2\|^2}{2l_{x}^2}-\frac{\|y_1-y_2\|^2}{2l_{y}^2}\bigg).
  \]
 We assume that $\log k(x,\omega)$  can be represented by the following truncate the KLE
\[
h(x,\omega):=\log k(x,\omega)={\bb E}[h(x,\omega)]+\sum_{i=1}^N \sqrt\lambda_i\theta_i(\omega)\varphi_i(x), \quad N=20,
\]
where the ${\bb E}[h(x,\omega)]$ is plotted in Figure \ref{coh}.

\begin{figure}[tbp]
  \centering
  \includegraphics[width=2.5in, height=1.7in]{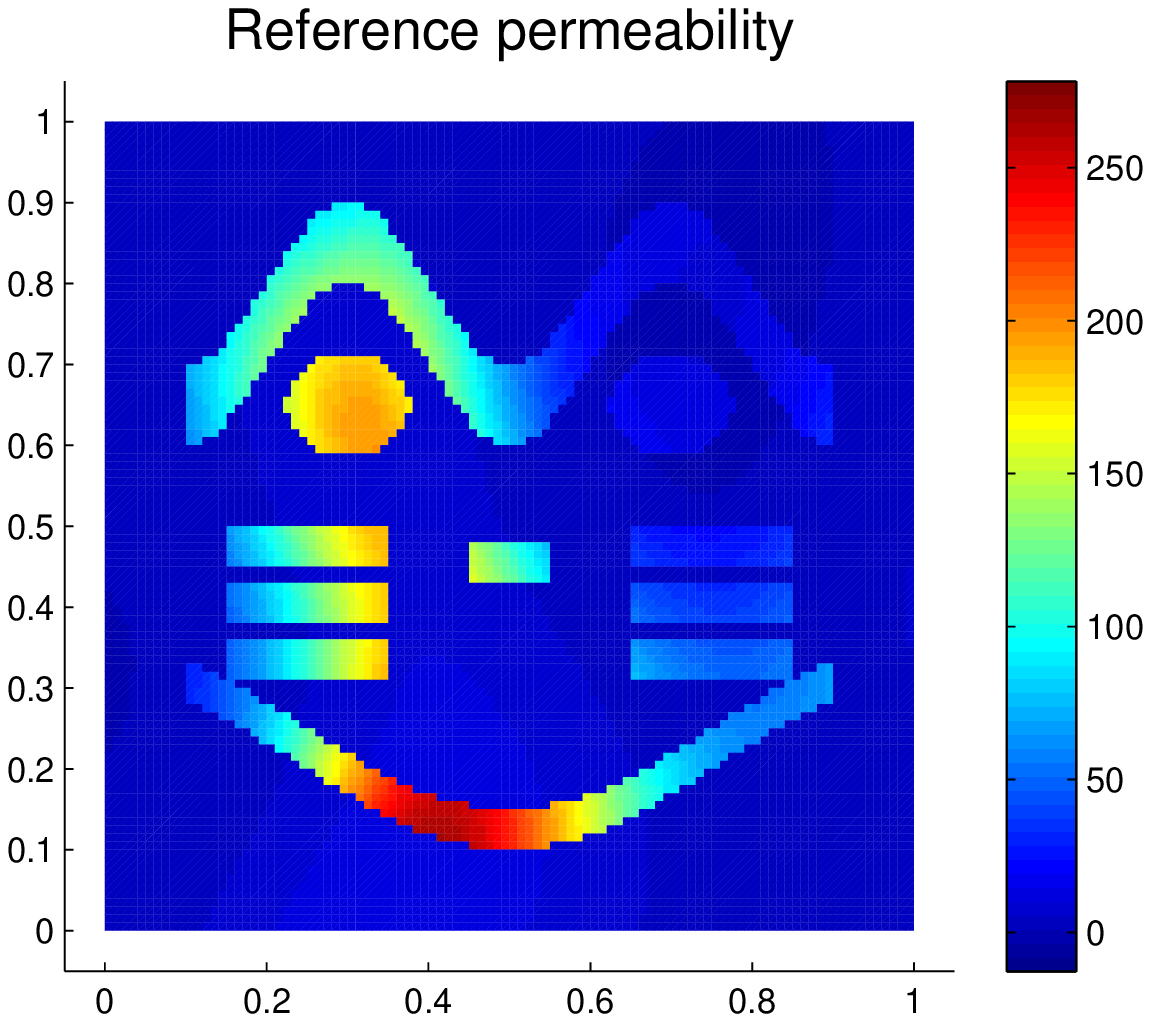}
  \includegraphics[width=2.5in, height=1.7in]{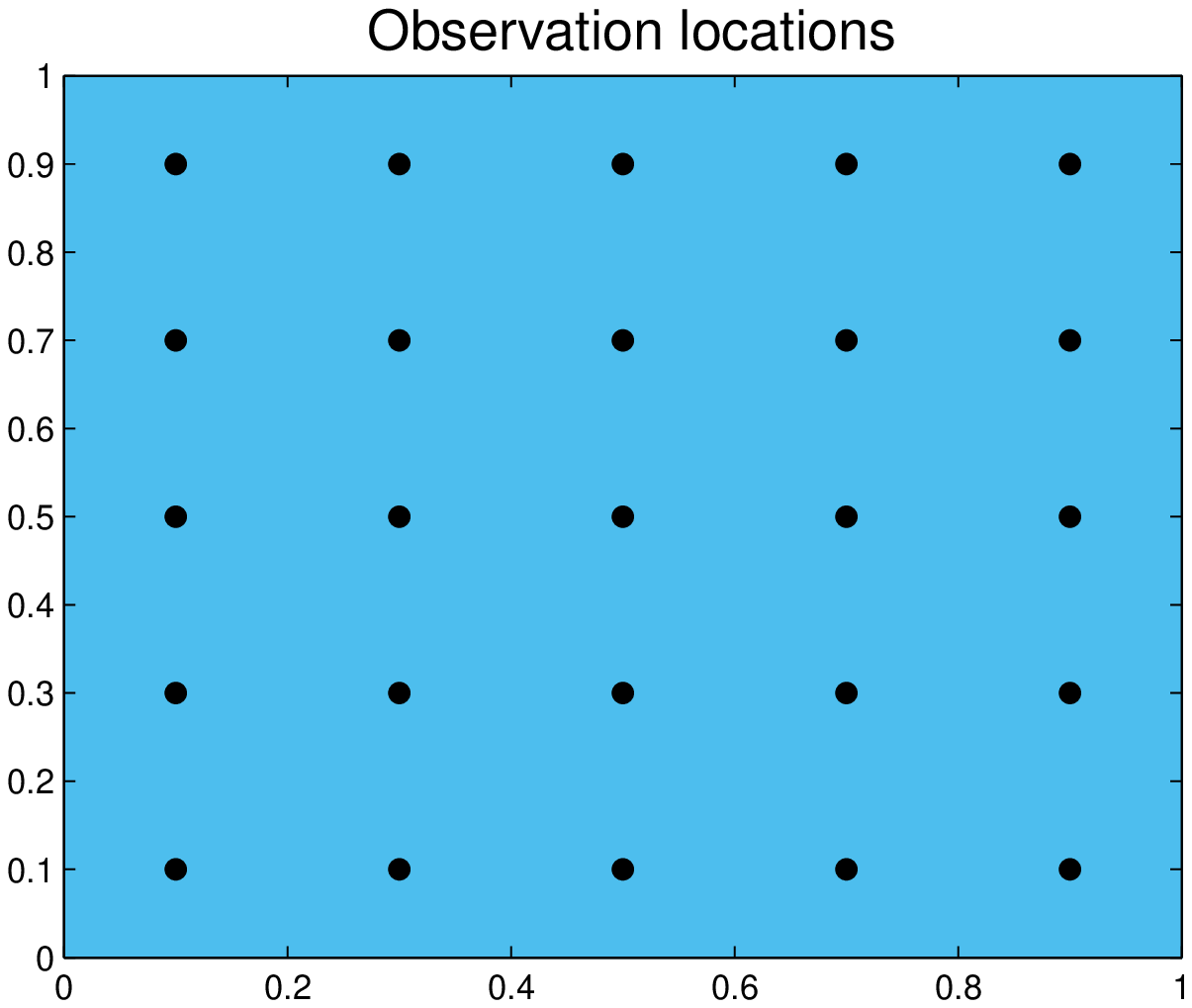}
  \caption{ The left is truth permeability of $k(x)$ and right is observation locations.}\label{coeffhighcon}
 \end{figure}

The ground true parameter $\bt^*$ are randomly drawn from the standard multivariate normal distribution and the truth permeability map is depicted in Figure \ref{coeffhighcon} (left). The total number of  data assimilation steps is set as $I_2=9$. Measurements are taken at time instances $0.012+0.01I:0.002:0.018+0.01I$ in each data assimilation step, where $I\in\{1,2,\cdots,9\}$ and the locations are distributed on the uniform $5\times 5$ grid of the domain $[0.1, 0.9]\times [0.1, 0.9]$ as shown in Figure \ref{coeffhighcon} (right).

The full order forward model is defined  on $100\times100$ uniform fine grid, and GMsFEM is implemented on $5\times 5$ coarse grid. We choose the number of eigenfunctions selected in calculating the snapshot space as $M^i_{\text{snap}}=20$. GMsFEM with $10$ offline local multiscale basis functions (i.e., $M_c=10$) is used to solve optimization problem (\ref{cos}).

\begin{figure}[tbp]
  \centering
  \includegraphics[width=1.5in, height=1.4in]{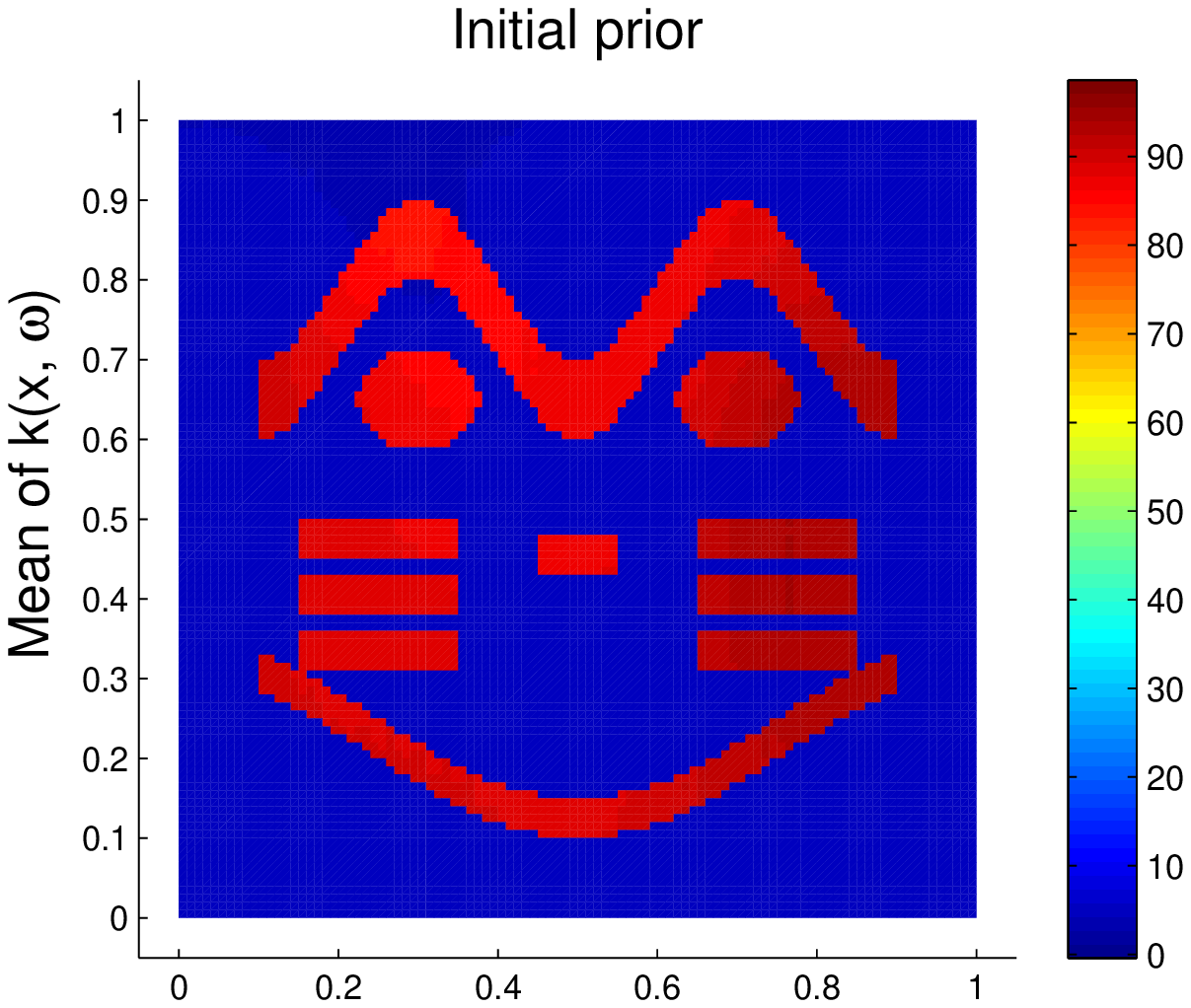}
  \includegraphics[width=1.5in, height=1.4in]{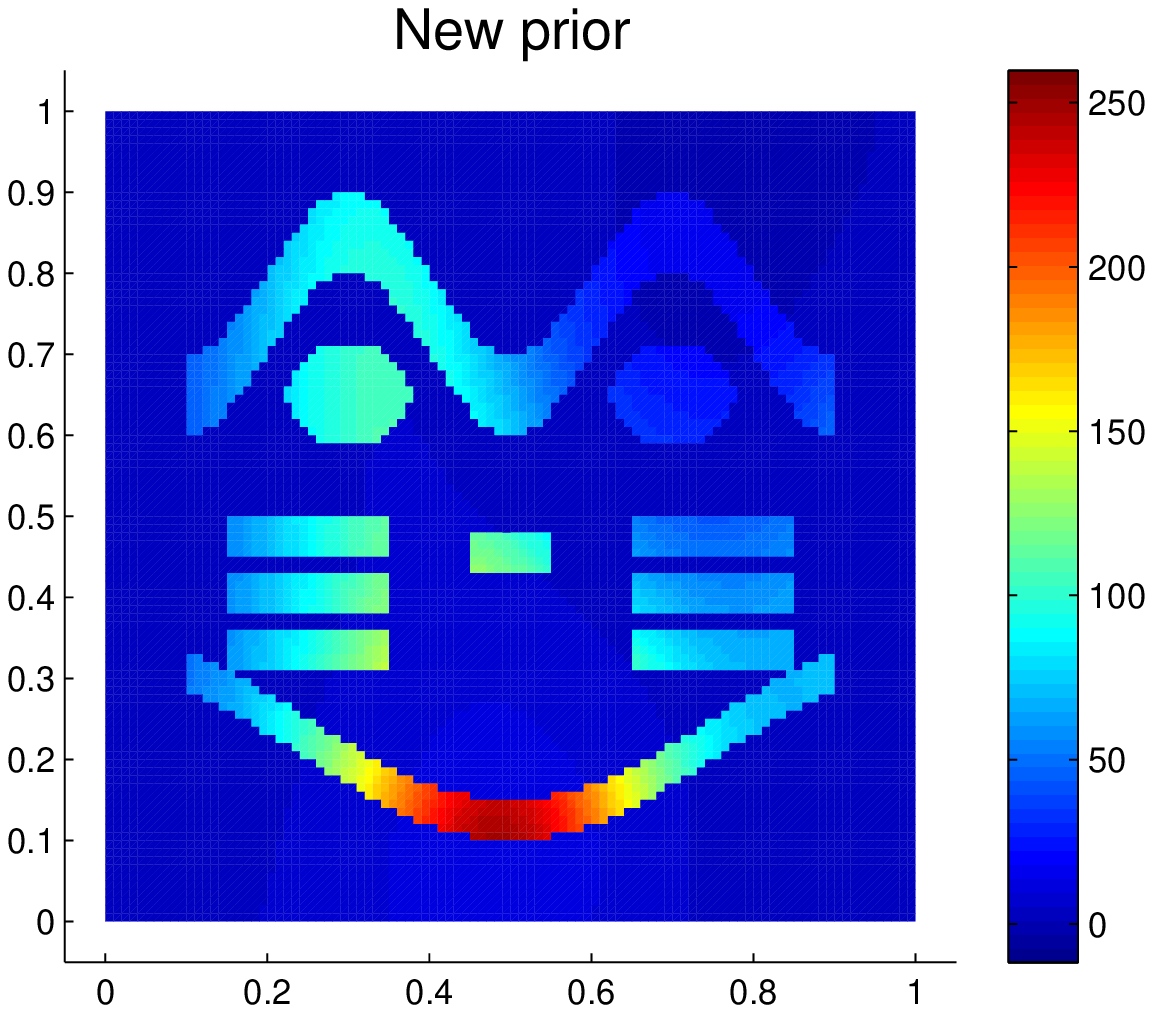}
  \includegraphics[width=1.5in, height=1.4in]{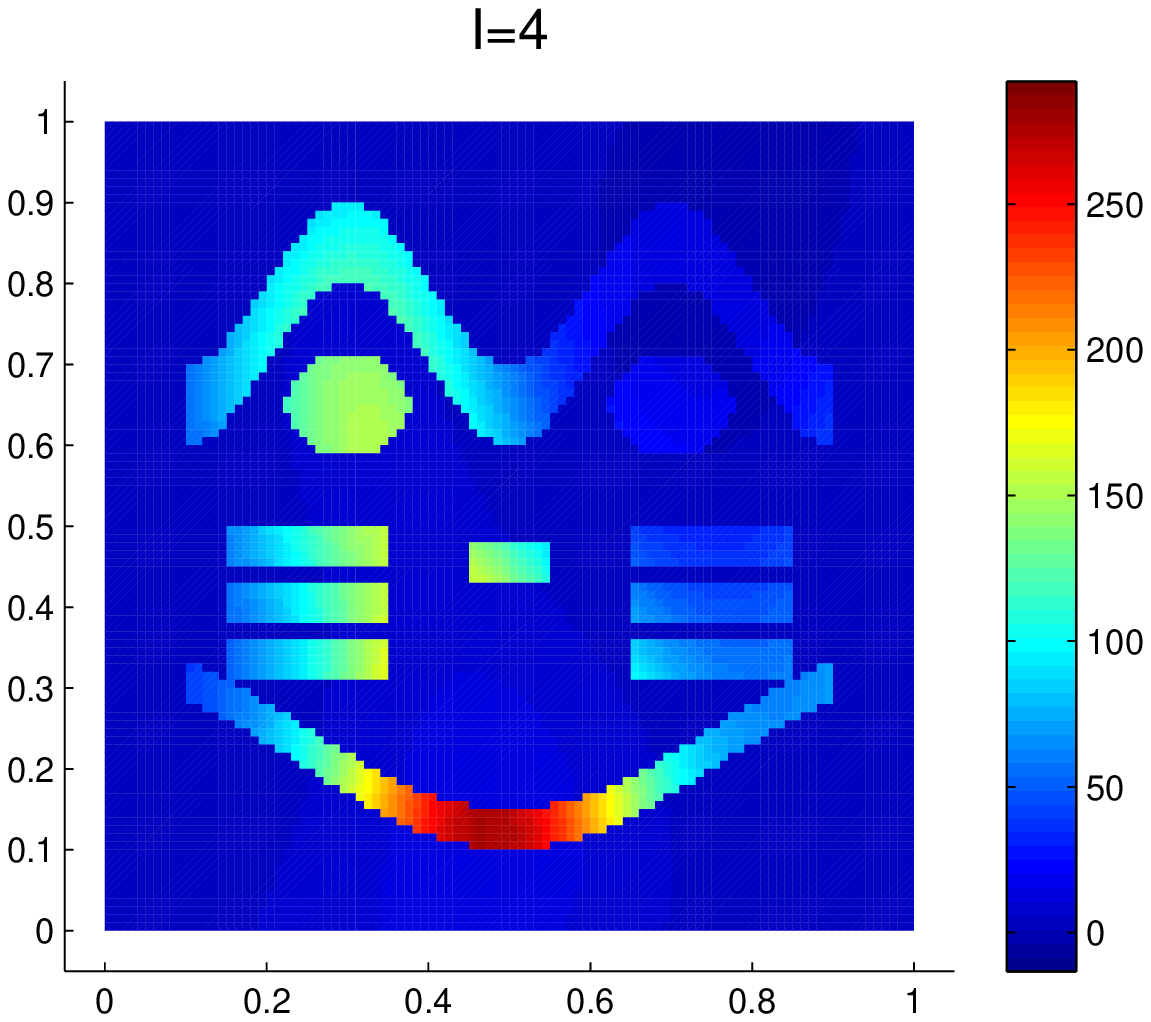}
  \includegraphics[width=1.5in, height=1.4in]{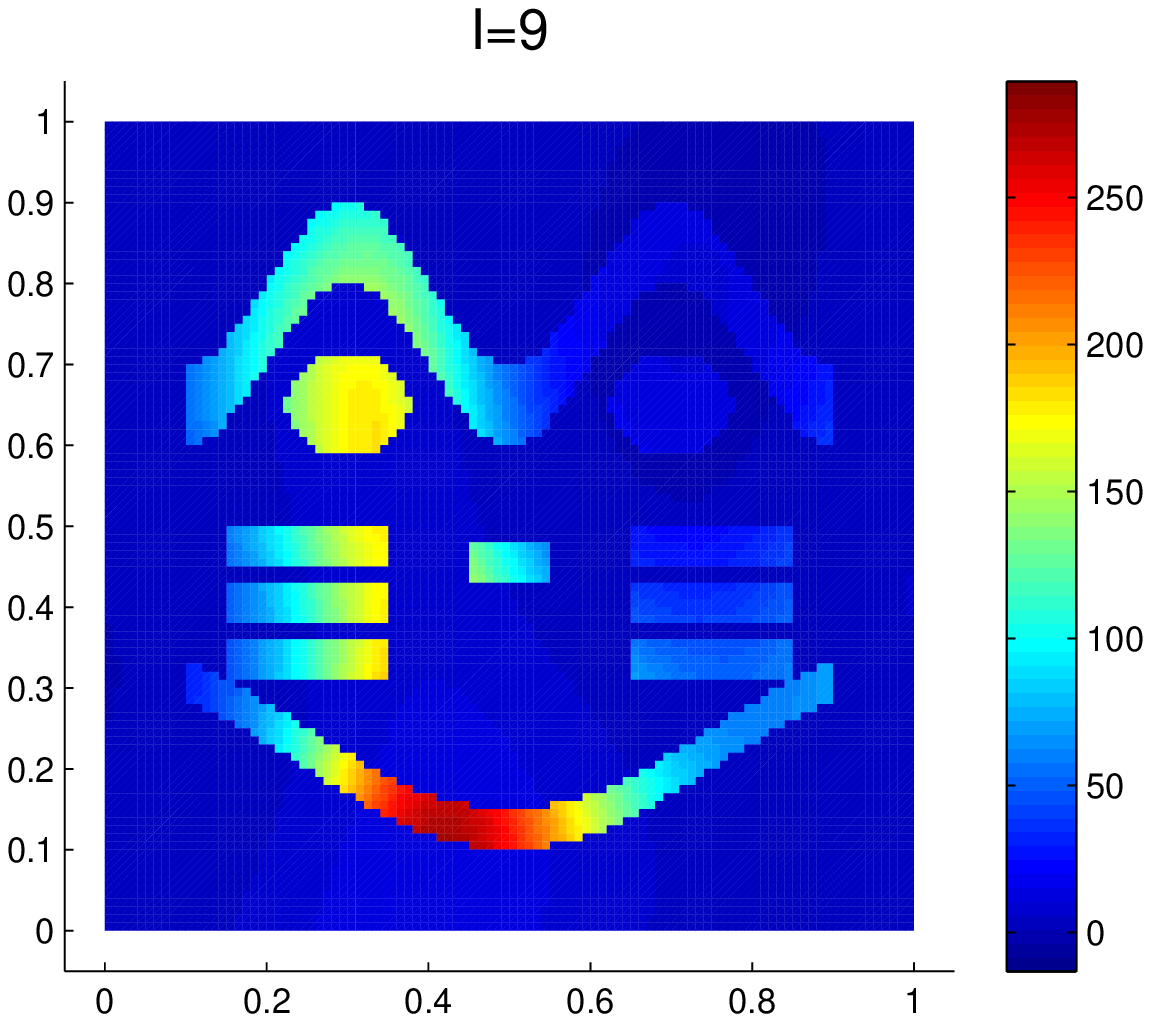}
  \includegraphics[width=1.5in, height=1.4in]{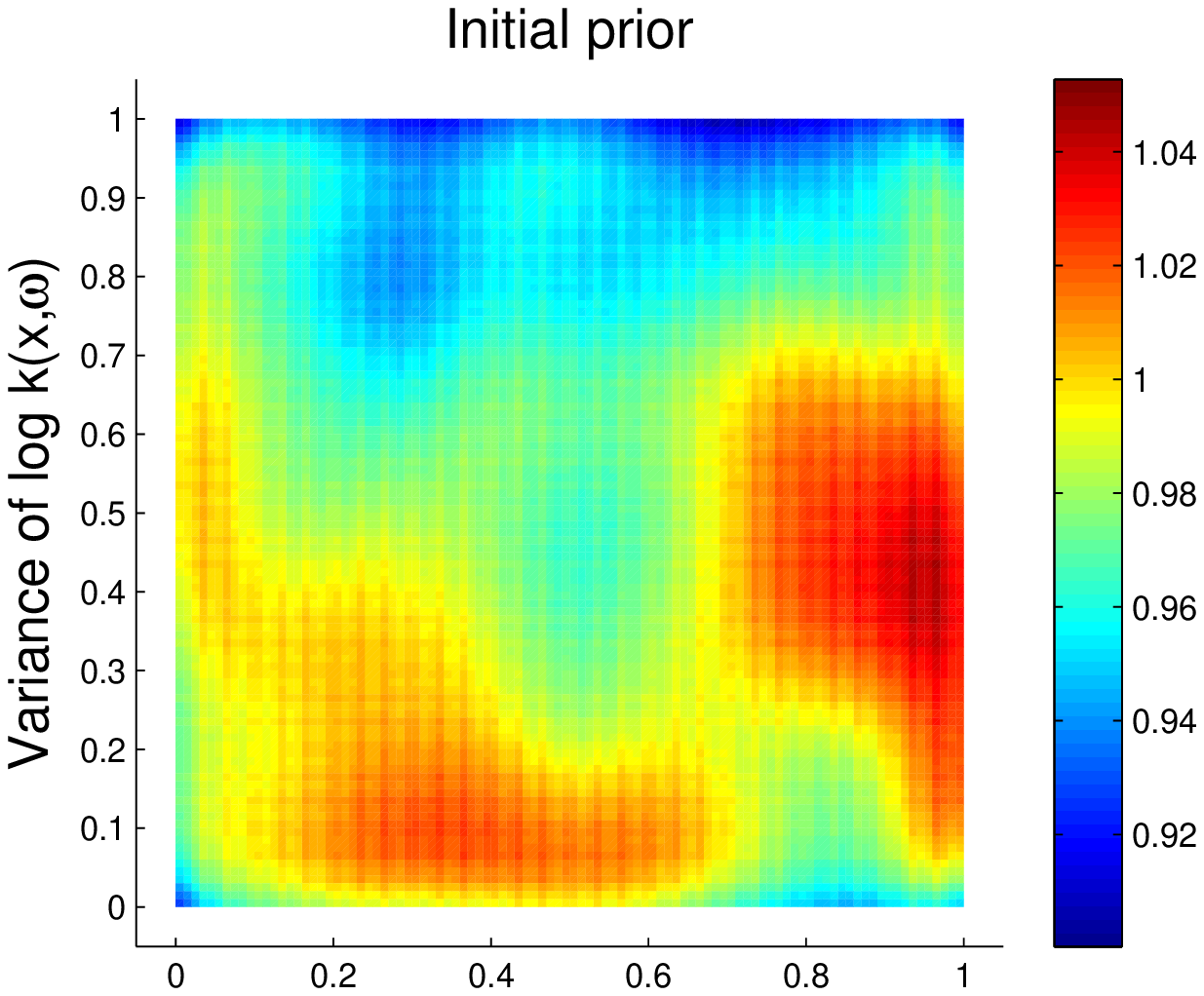}
  \includegraphics[width=1.5in, height=1.4in]{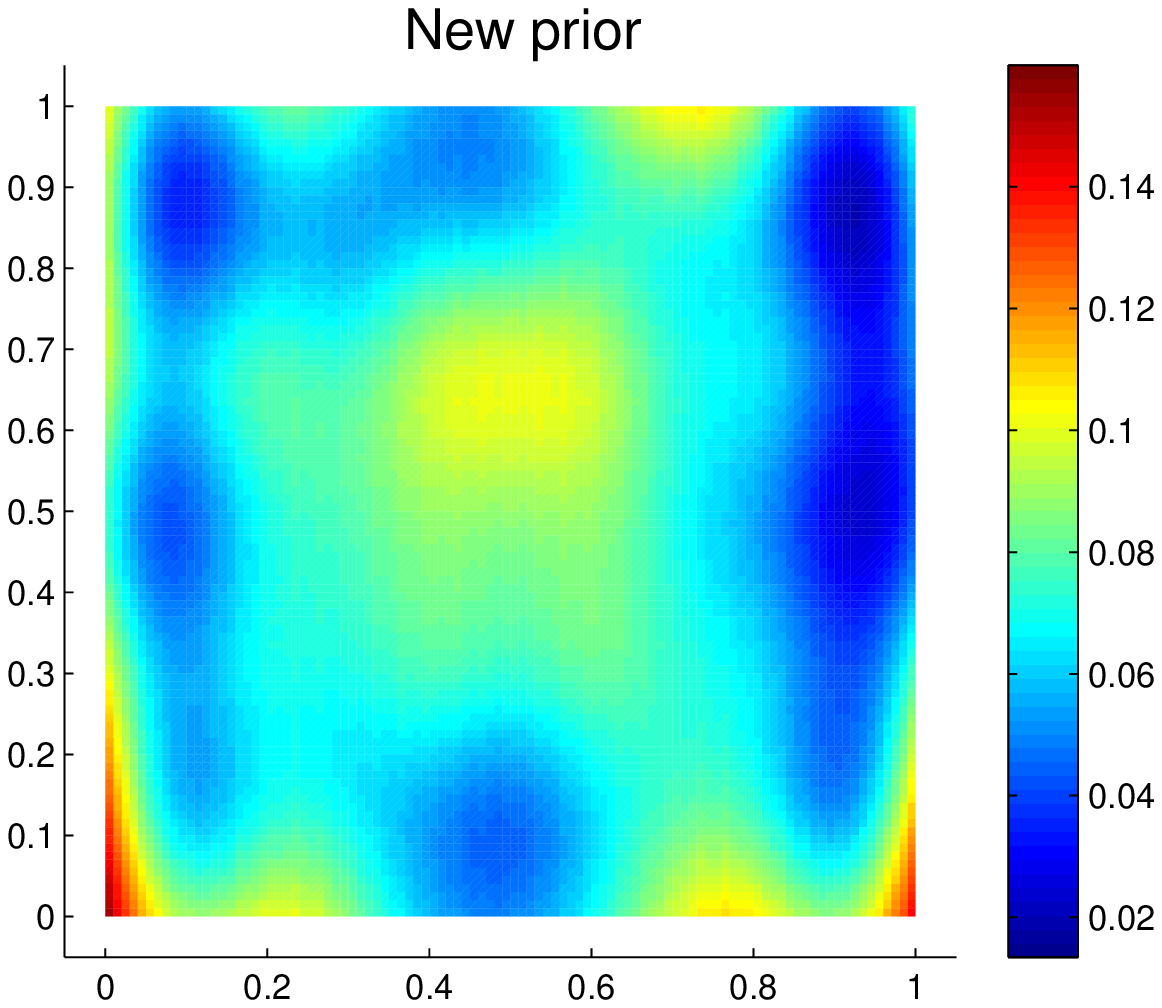}
  \includegraphics[width=1.5in, height=1.4in]{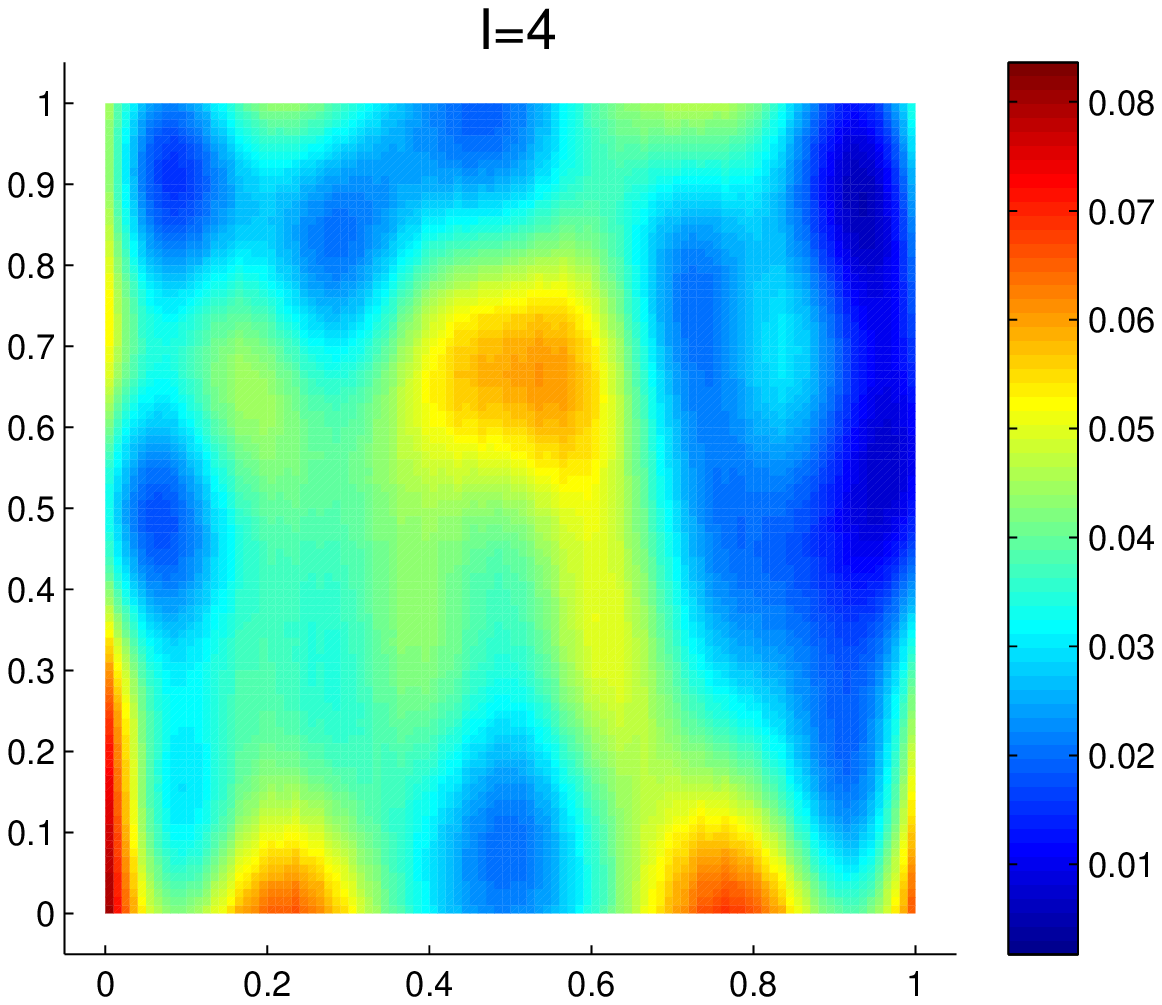}
  \includegraphics[width=1.5in, height=1.4in]{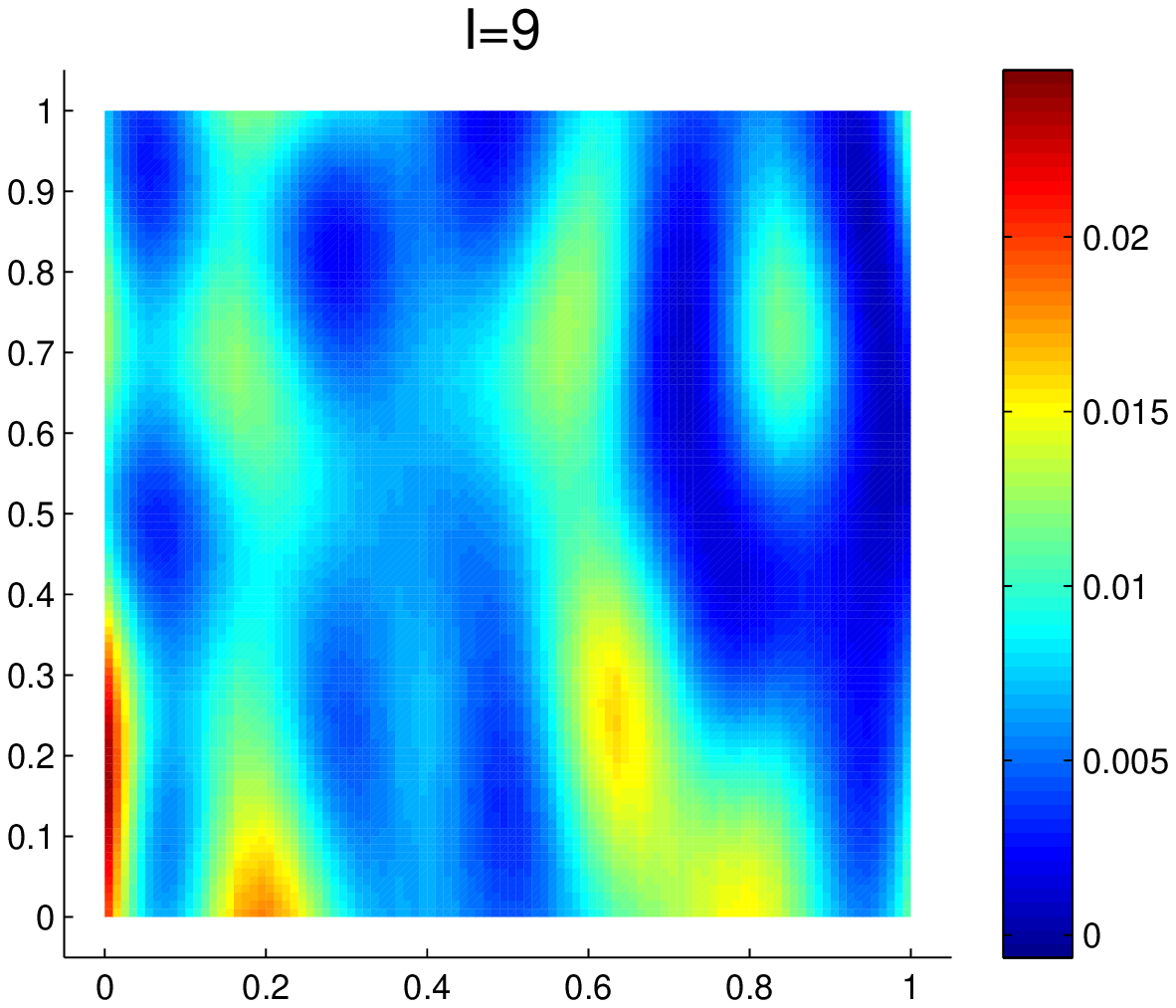}
  \caption{Mean of $k(x,\omega)$ and variance of log $k(x,\omega)$ by two-stage EnKF at different assimilation steps}\label{cme}
 \end{figure}

In the numerical simulation, we use $3\times 10^3$  ensemble members to construct a new prior by standard EnKF in the first stage, where the first two levels of data information is used.
Three online  local multiscale basis functions ($M_i=3$) are selected in constructing the matrix $R$ to construct coarse model to obtain the new prior. Then, we construct the surrogate model after obtaining the new prior ensemble in the second stage. Five online  local multiscale basis functions ( $M_i=5$) are selected  in constructing the multiscale basis  matrix  to construct the gPC surrogate model. When the order of gPC is set as $N_0=3$, we use $1771$ samplers  when computing vectors $\bf b$ and $\bf A$ in Section \ref{SCM}. In this stage, the number of ensemble members is set as $10^4$. By Algorithm \ref{IN}, we need to simultaneously update the hyperparameter $\sigma^2$. We take $n_s=0.05$ in in Algorithm \ref{Seq}.  The  $\sigma^2$ for true measurement noise is set as $10^{-4}$. The  posterior mean and posterior variance via  the data assimilation step are illustrated in Figure \ref{cme}, which shows that the uncertainty mainly lies on high-contrast part and boundary.

 \begin{figure}[tbp]
  \centering
  \includegraphics[width=3.3in, height=2.5in]{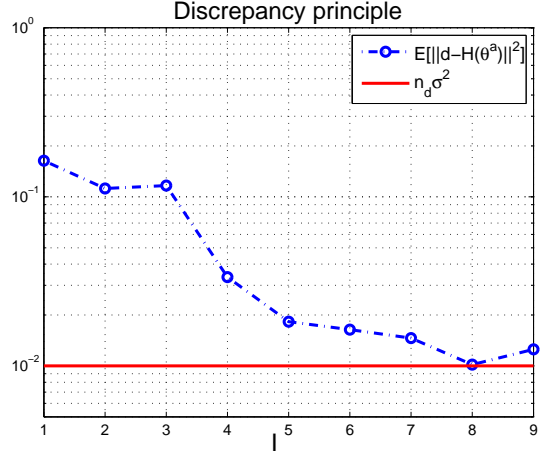}
  \caption{$\bb{E}[\|\bm{d}-\bb{H}(\bt^a)\|^{2}]$ v.s. $n_d\sigma^{2}$ with data assimilation step, in which $I=1$ is the initial prior and $I=2$ is the new prior.}\label{coUnba}
 \end{figure}

\begin{figure}[tbp]
  \centering
  \includegraphics[width=2.5in, height=1.7in]{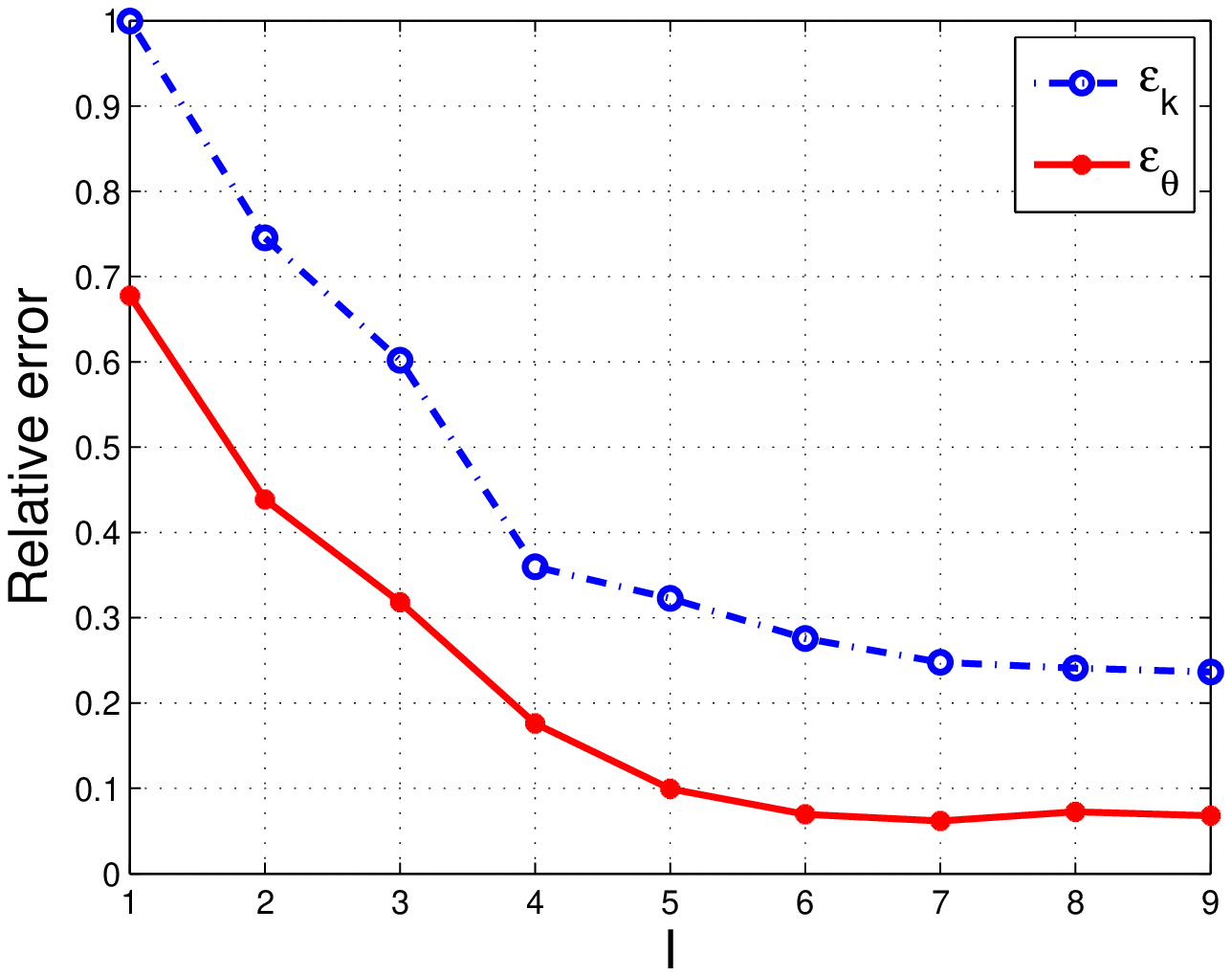}
  \includegraphics[width=2.5in, height=1.7in]{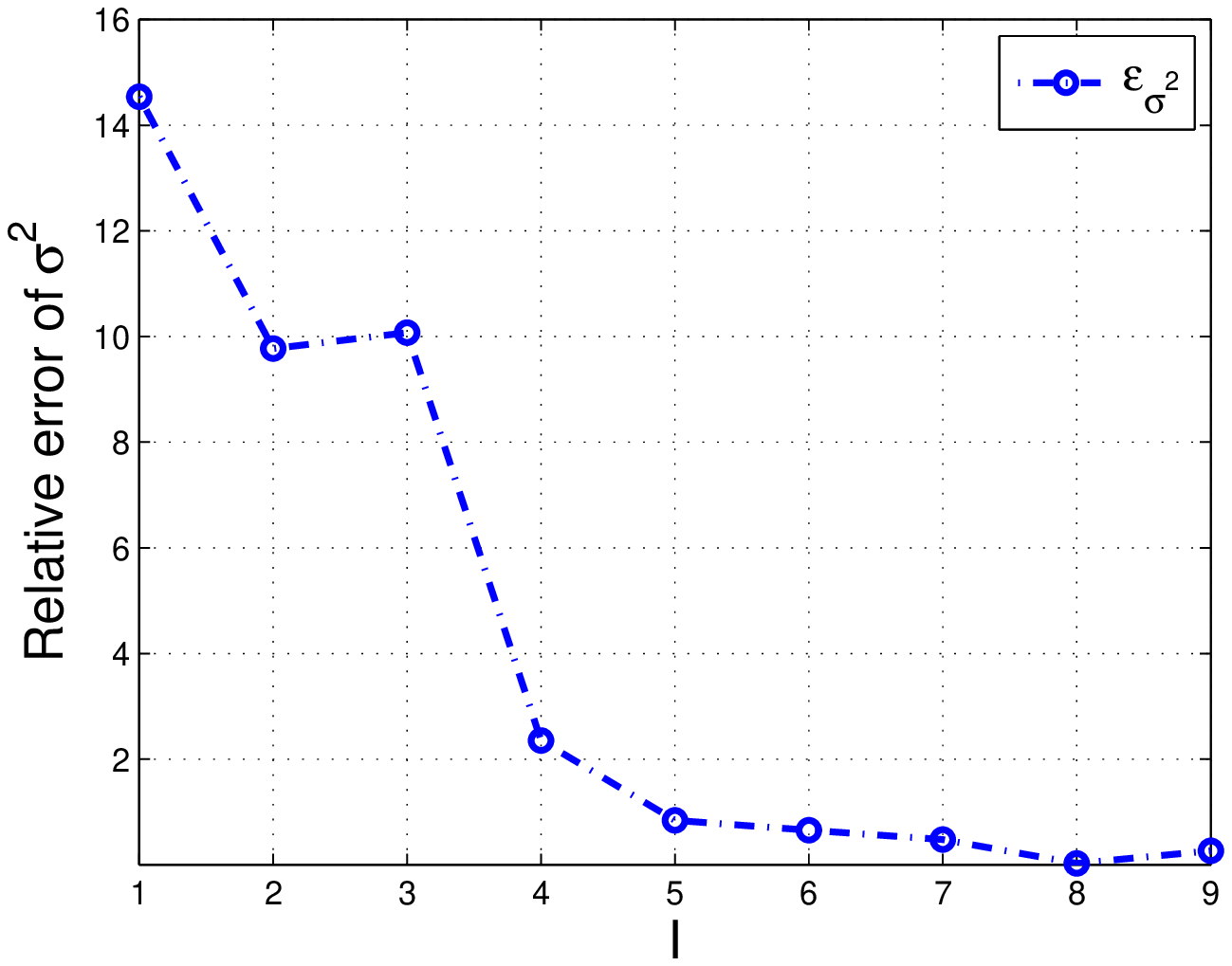}
  \caption{The relative error of permeability field and $\bt$ (left), the  relative error of $\sigma^2$ (right). $I=1$ is the initial prior, and $I=2$ is the new prior.}\label{error}
 \end{figure}

Figure $\ref{coUnba}$ depicts the discrepancy principle when $\sigma^2$ is unknown in the forecast step. As we have mentioned before in Subsection \ref{channel}, the measurement noise is the additional Gaussian type. Then we can see that  the expectation of the discrepancy tends to $n_d\sigma^2$ as assimilation time moves on.
To measure the estimate accuracy, we  define the relative errors  corresponding to the posterior distribution  by
\begin{equation*}
  \quad \varepsilon_k:= \frac{\|k(x,\bar{\bt})-k(x,\bt^*)\|}{\|k(x,\bt^*)\|},
  \quad \varepsilon_{\sigma^2} := \frac{\|\sigma^2-{\sigma_t^2}\|}{\|\sigma_t^2\|},
  \quad \varepsilon_\theta:= \frac{\|\bar{\bt}-{\bt}\|}{\|{\bt}\|},
\end{equation*}
where $\bar{\bt}$ and $\sigma^2$ are the mean of final assimilation step by the two-stage EnKF and $\sigma_t^2$ is the truth observation noise. We just consider the error of the second stage in the two-stage EnKF method. The relative error is plotted in the Figure $\ref{error}$. As expected, the relative error gradually decreases with the assimilation steps.

\begin{figure}[tbp]
  \centering
  \includegraphics[width=1.5in, height=1.4in]{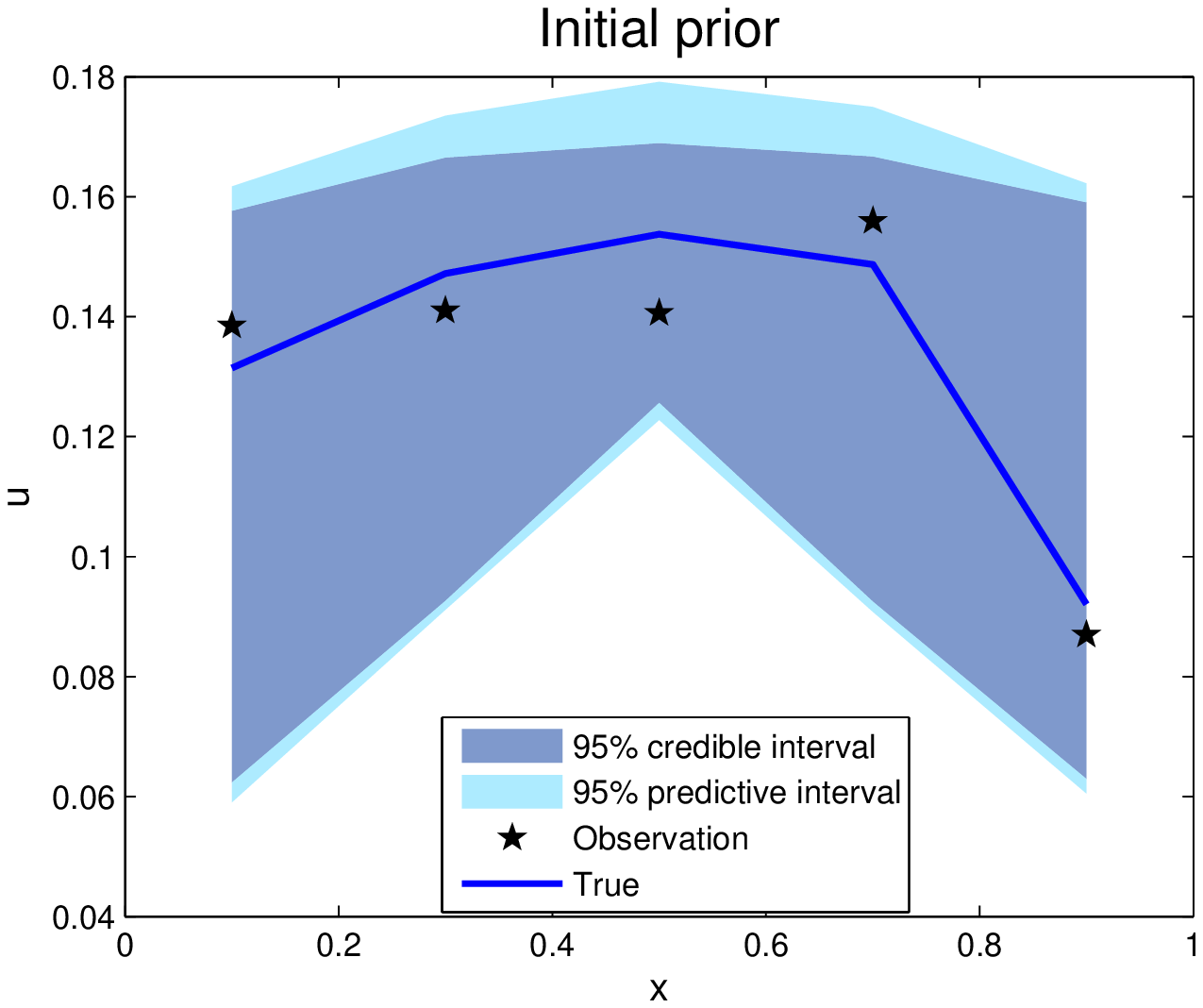}
  \includegraphics[width=1.5in, height=1.4in]{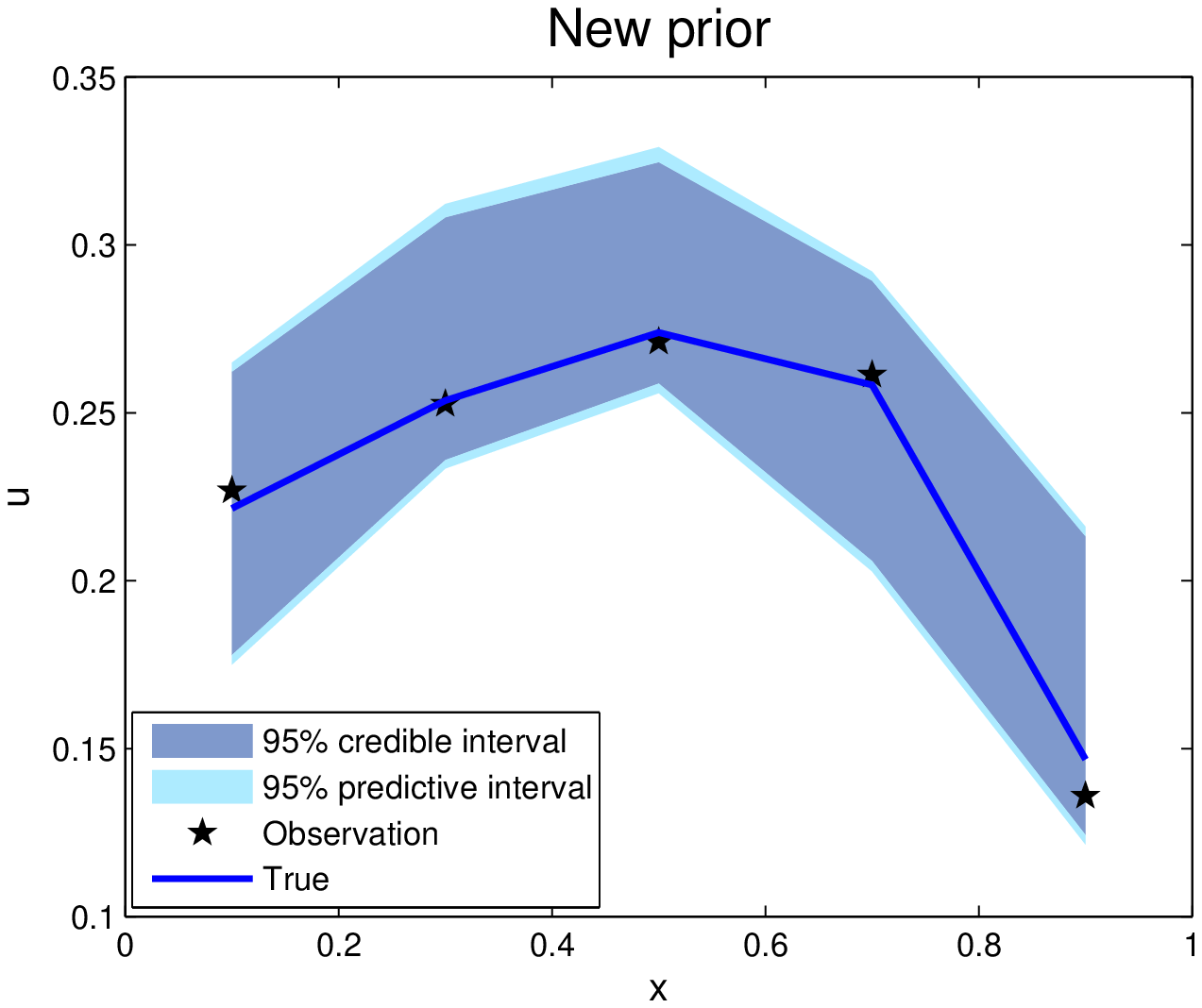}
  \includegraphics[width=1.5in, height=1.4in]{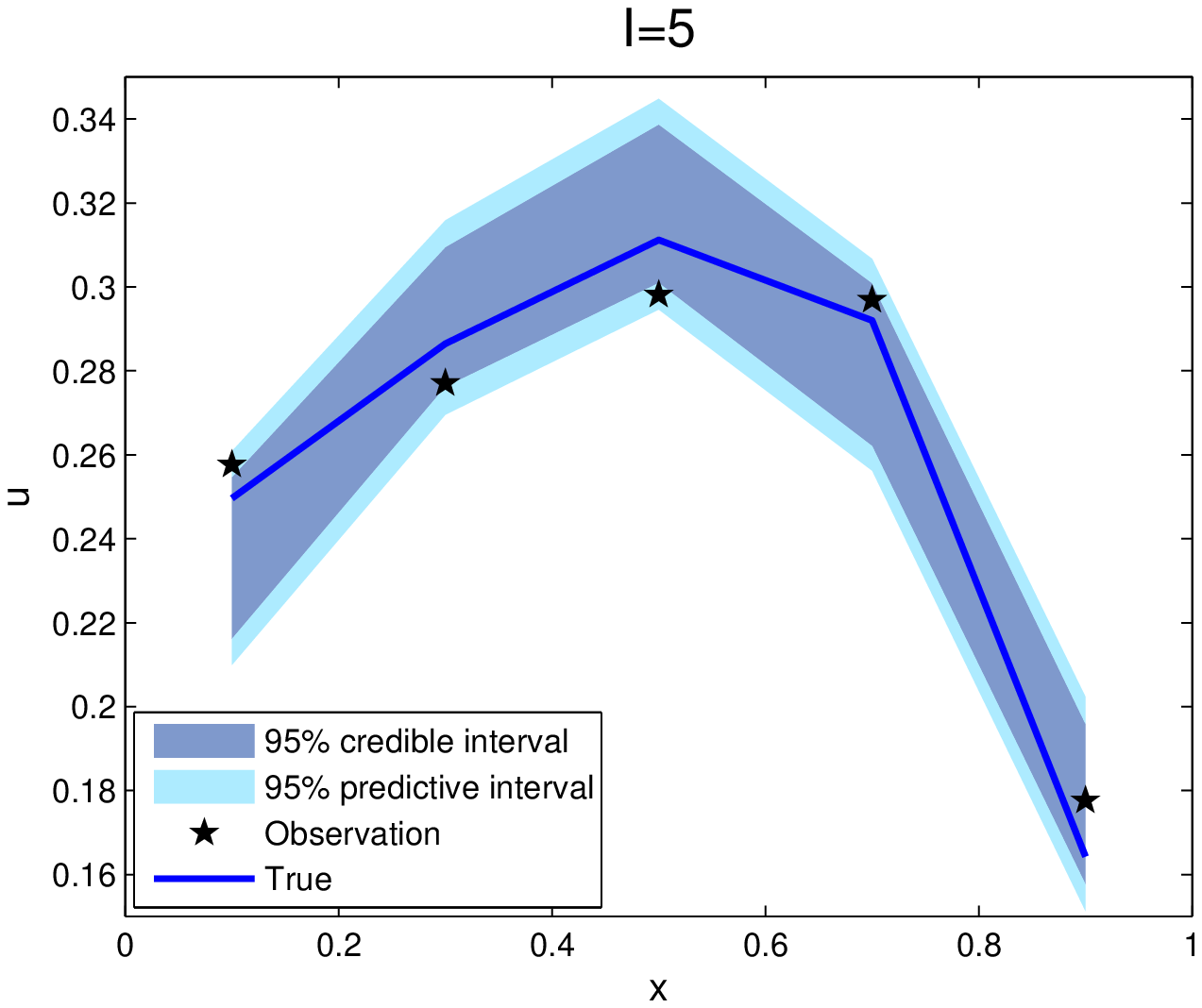}
  \includegraphics[width=1.5in, height=1.4in]{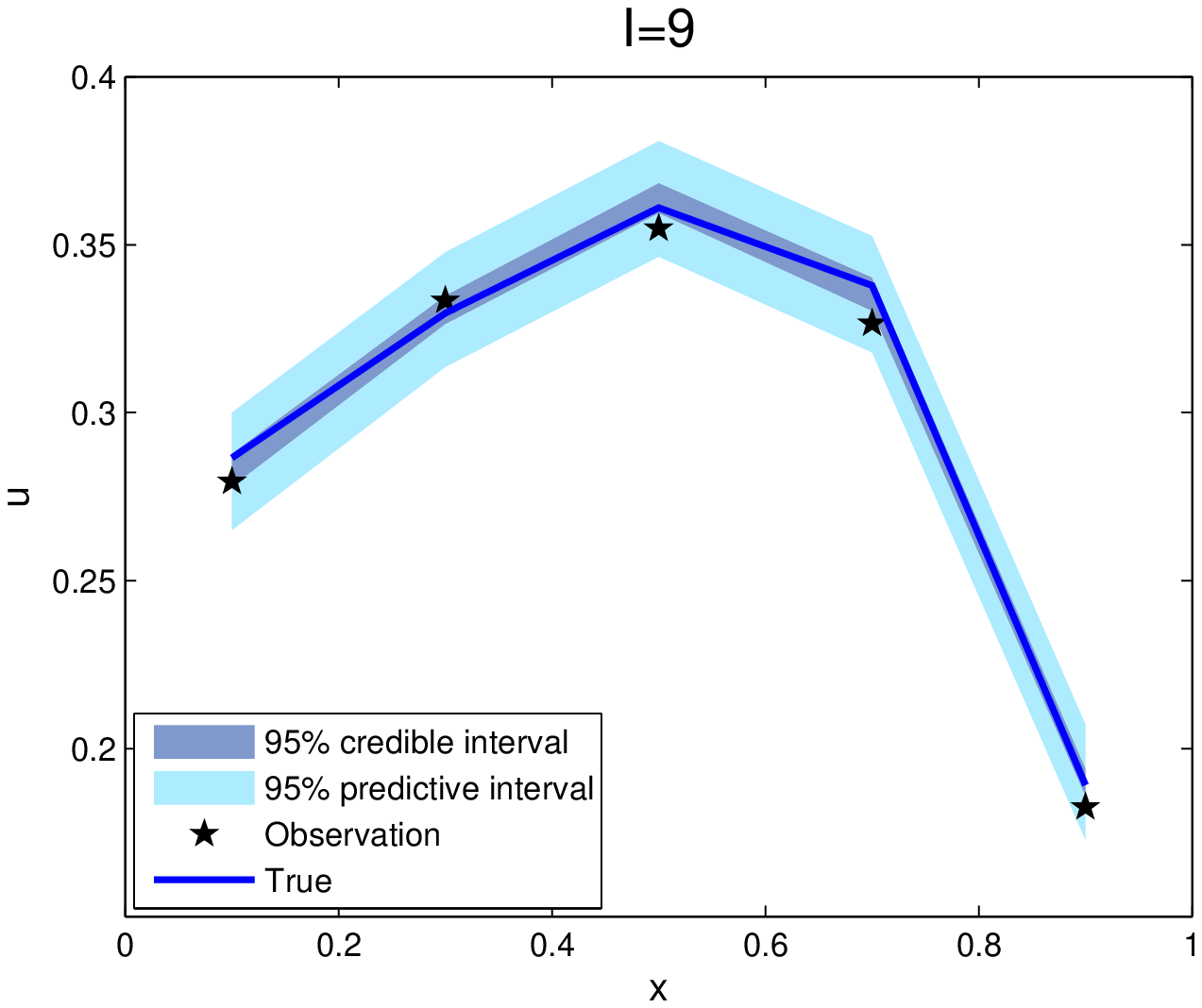}
  \caption{95\% predictive interval, 95\% credible interval, observation and true value of $u((x, 0.5); t)$.}\label{coprx}
\end{figure}

\begin{figure}[tbp]
  \centering
  \includegraphics[width=1.5in, height=1.4in]{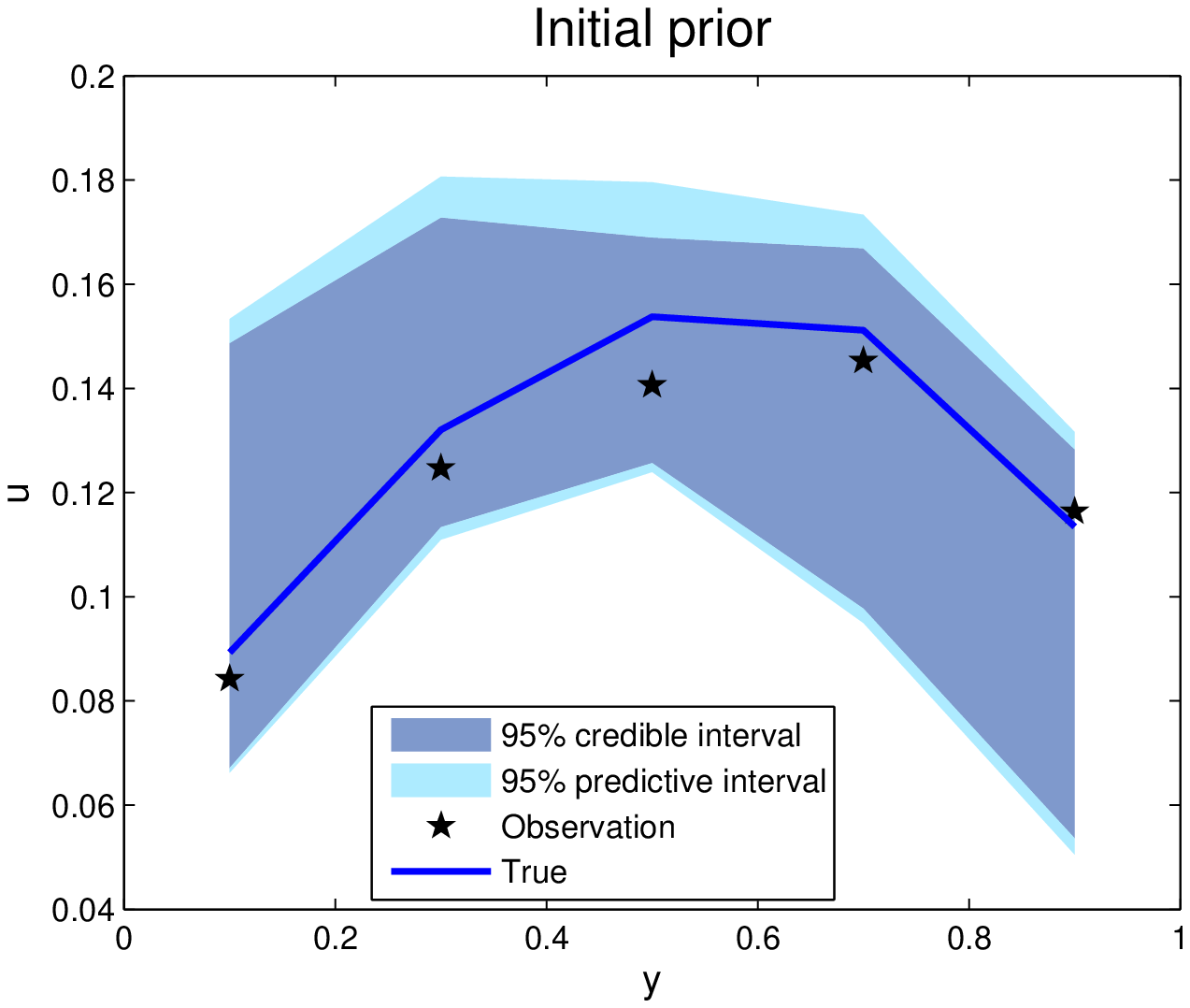}
  \includegraphics[width=1.5in, height=1.4in]{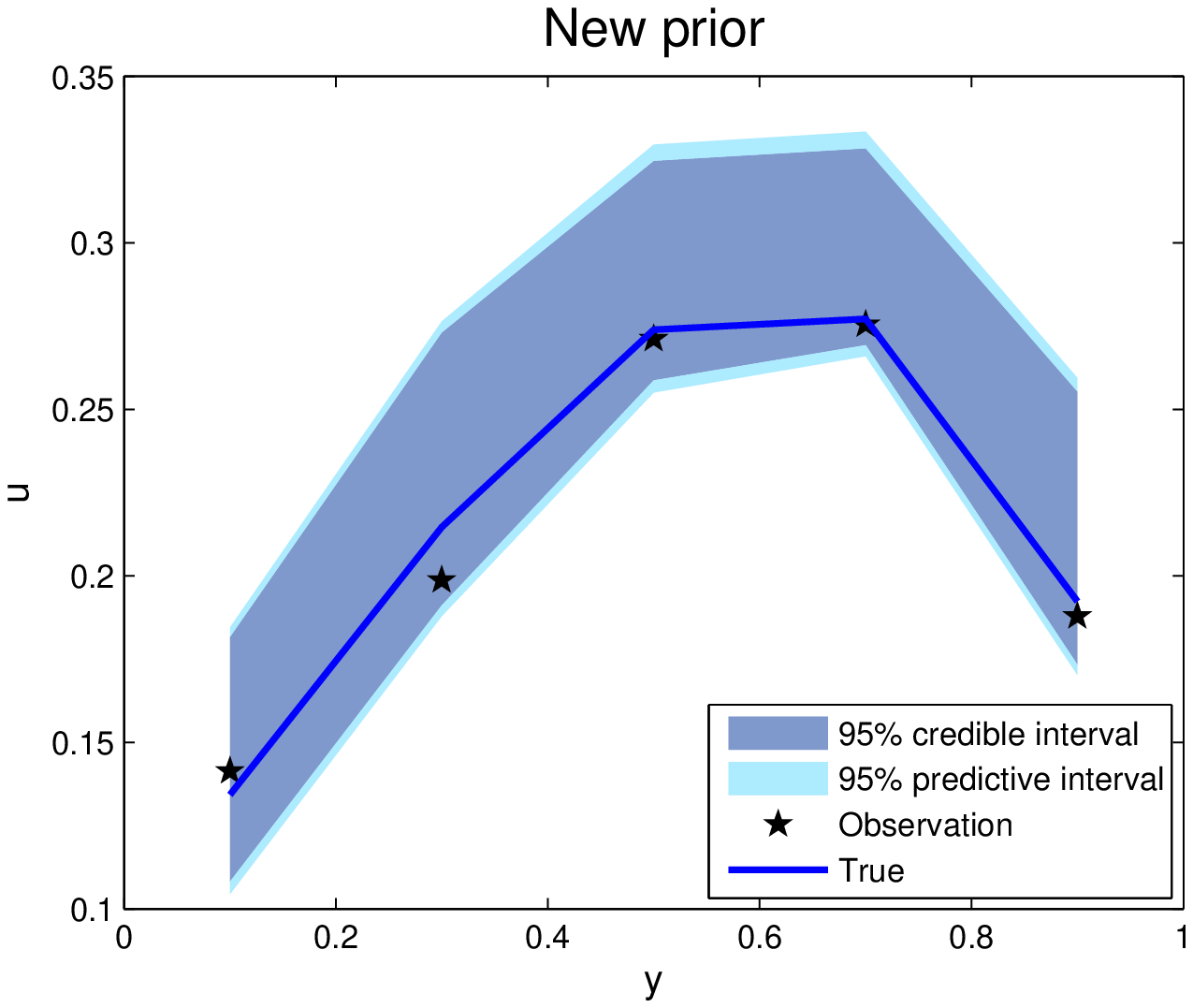}
  \includegraphics[width=1.5in, height=1.4in]{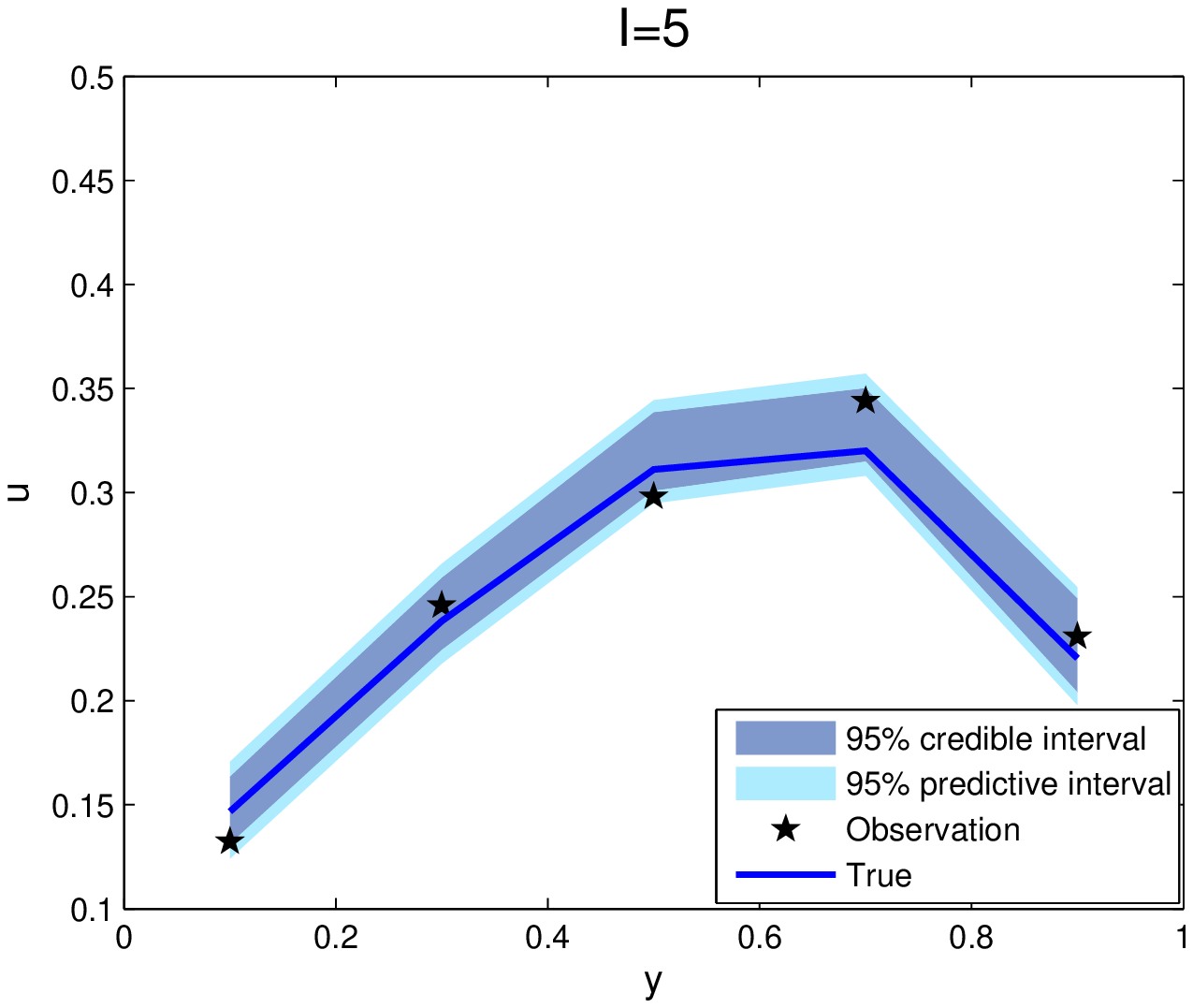}
  \includegraphics[width=1.5in, height=1.4in]{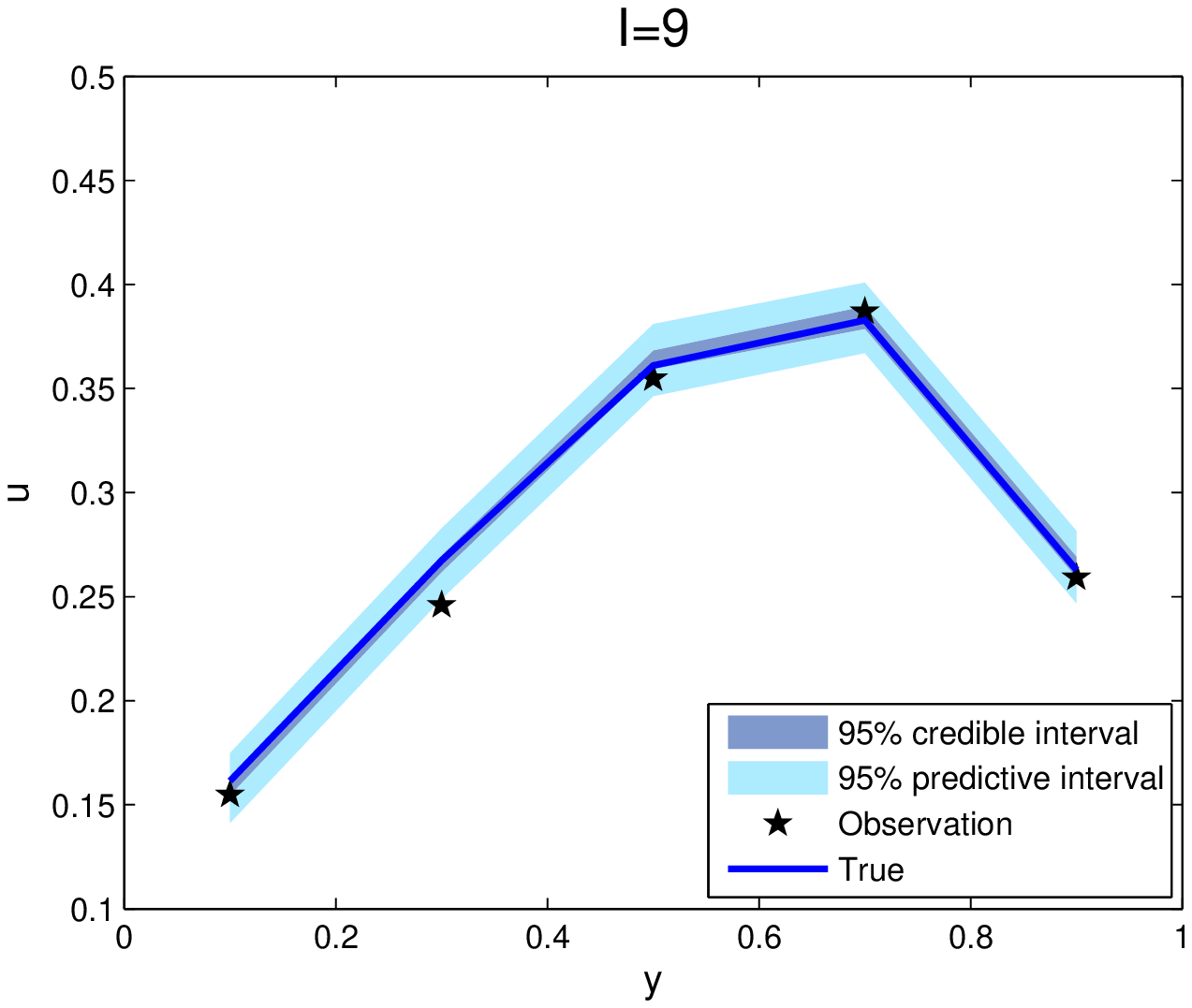}
  \caption{95\% predictive interval, 95\% credible interval, observation and true value of $u((0.5, y); t)$.}\label{copry}
\end{figure}

\begin{figure}[tbp]
  \centering
  \includegraphics[width=2.1in, height=1.7in]{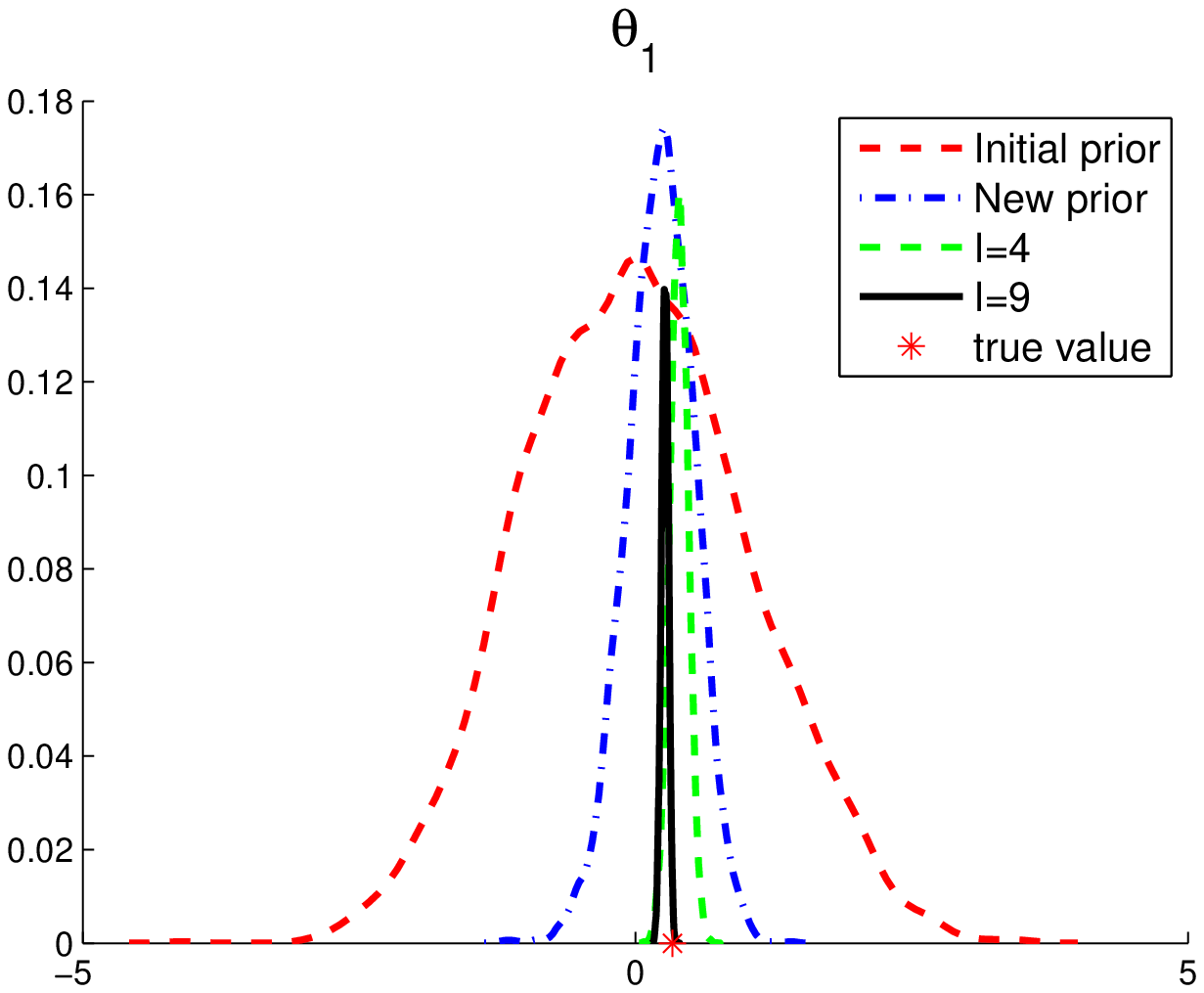}
  \includegraphics[width=2.1in, height=1.7in]{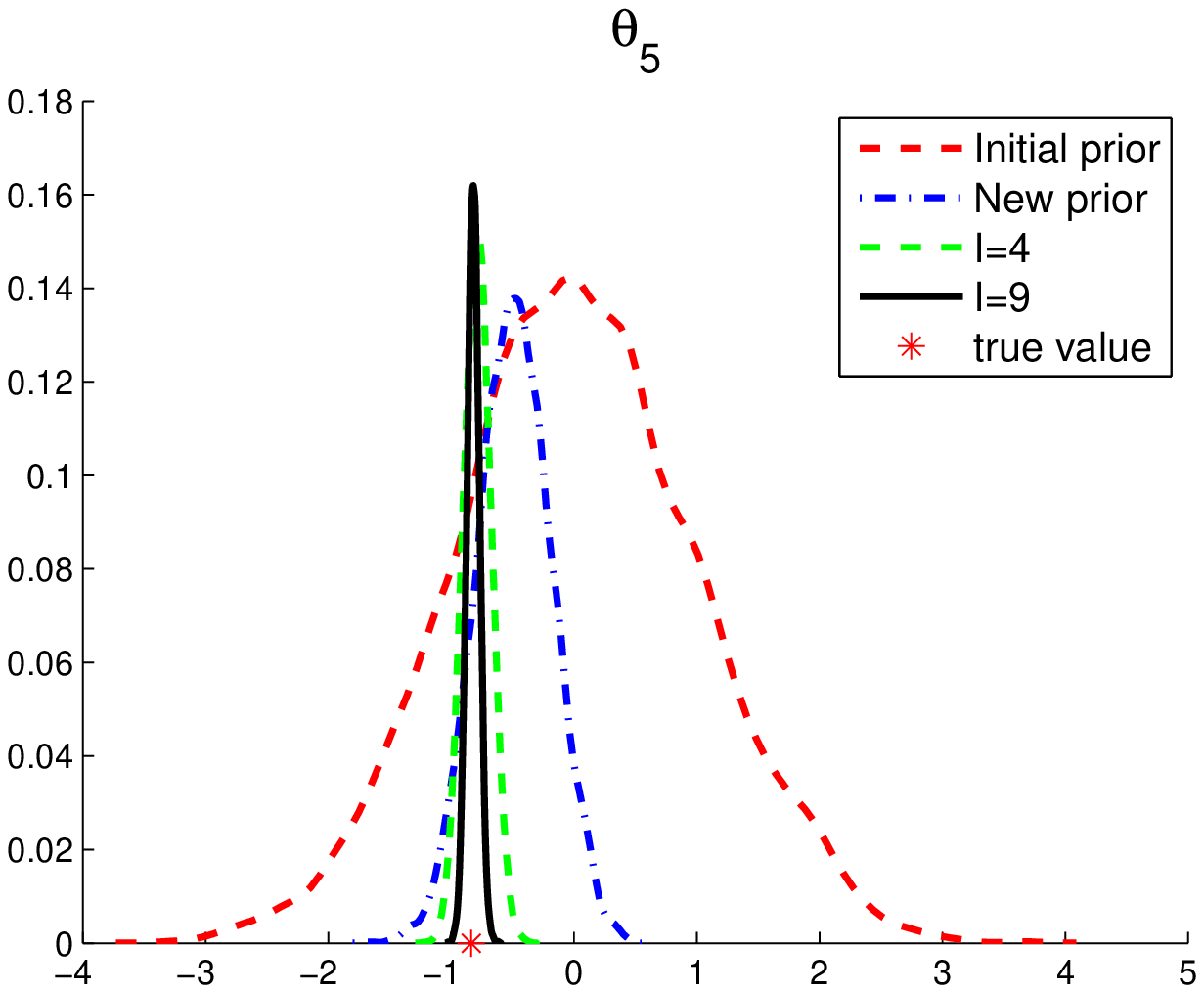}
  \includegraphics[width=2.1in, height=1.7in]{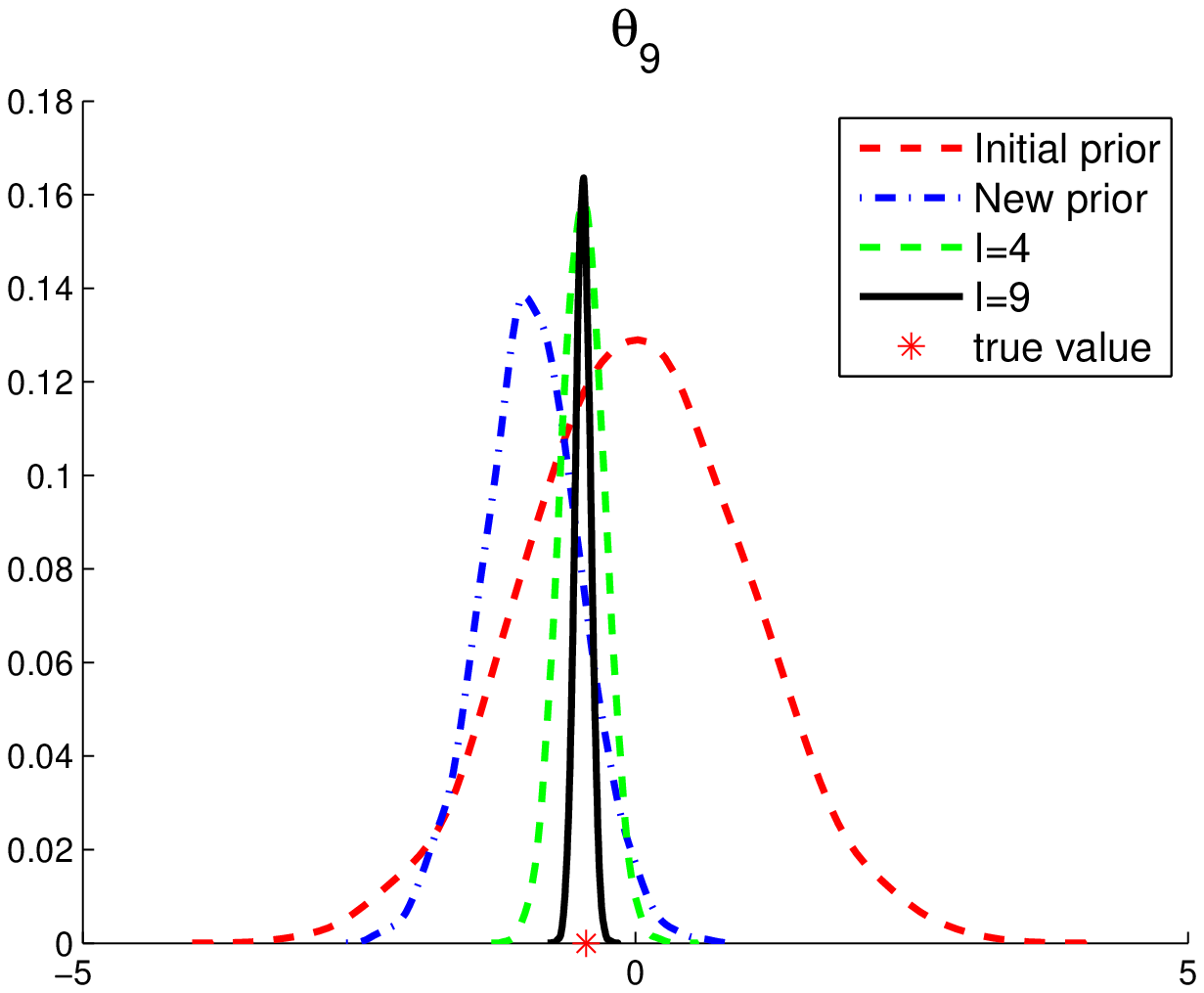}
  \includegraphics[width=2.1in, height=1.7in]{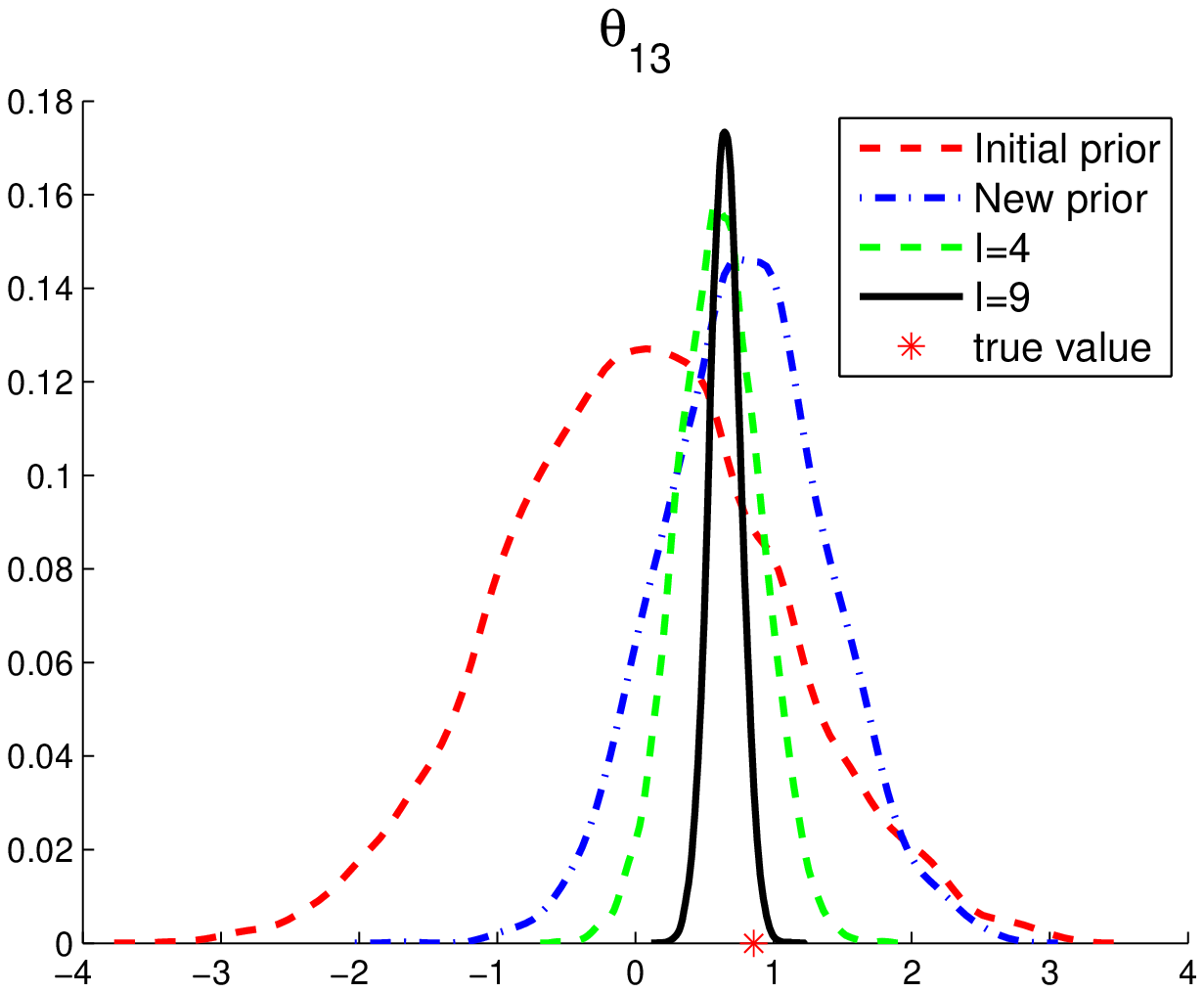}
  \includegraphics[width=2.1in, height=1.7in]{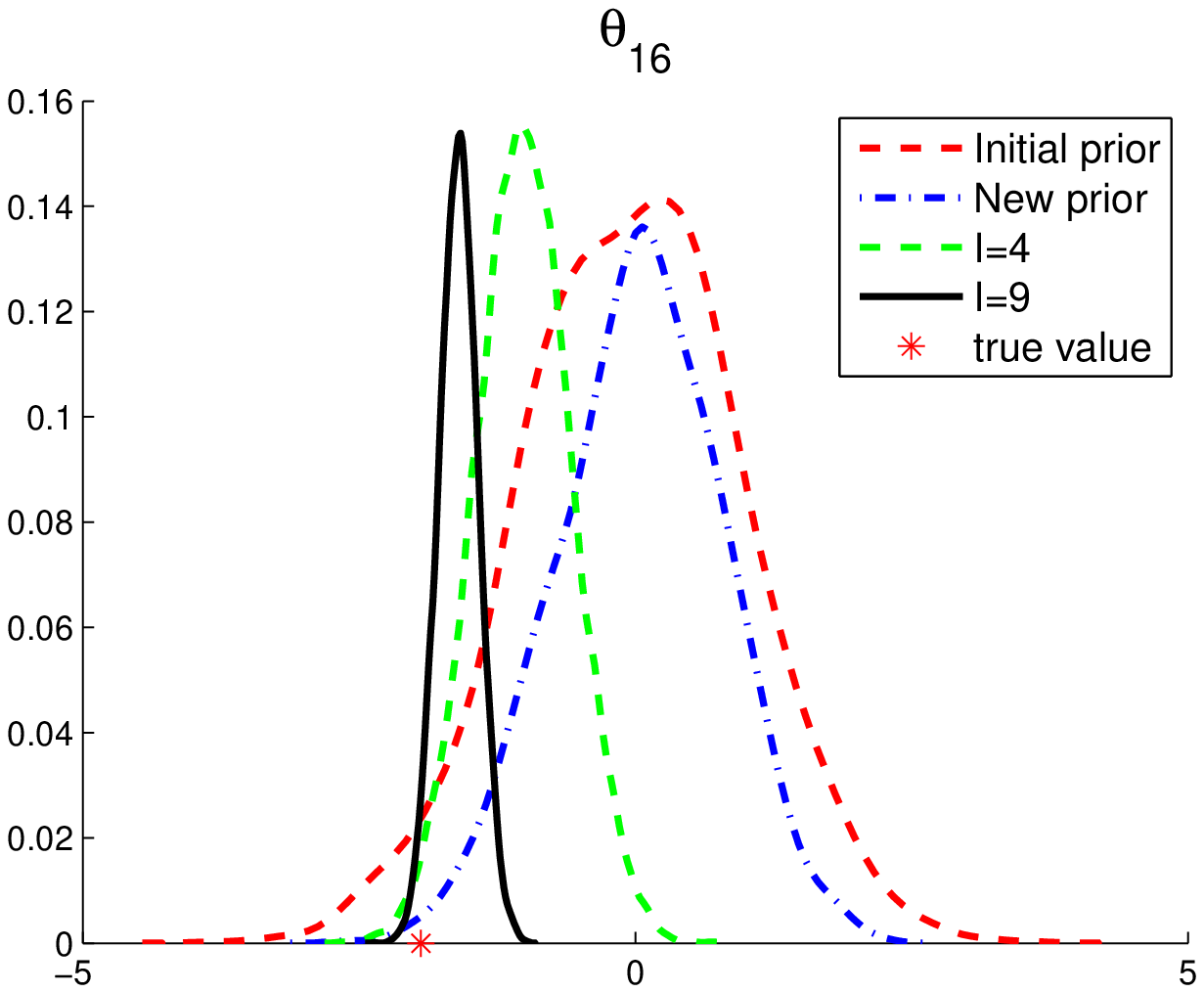}
  \includegraphics[width=2.1in, height=1.7in]{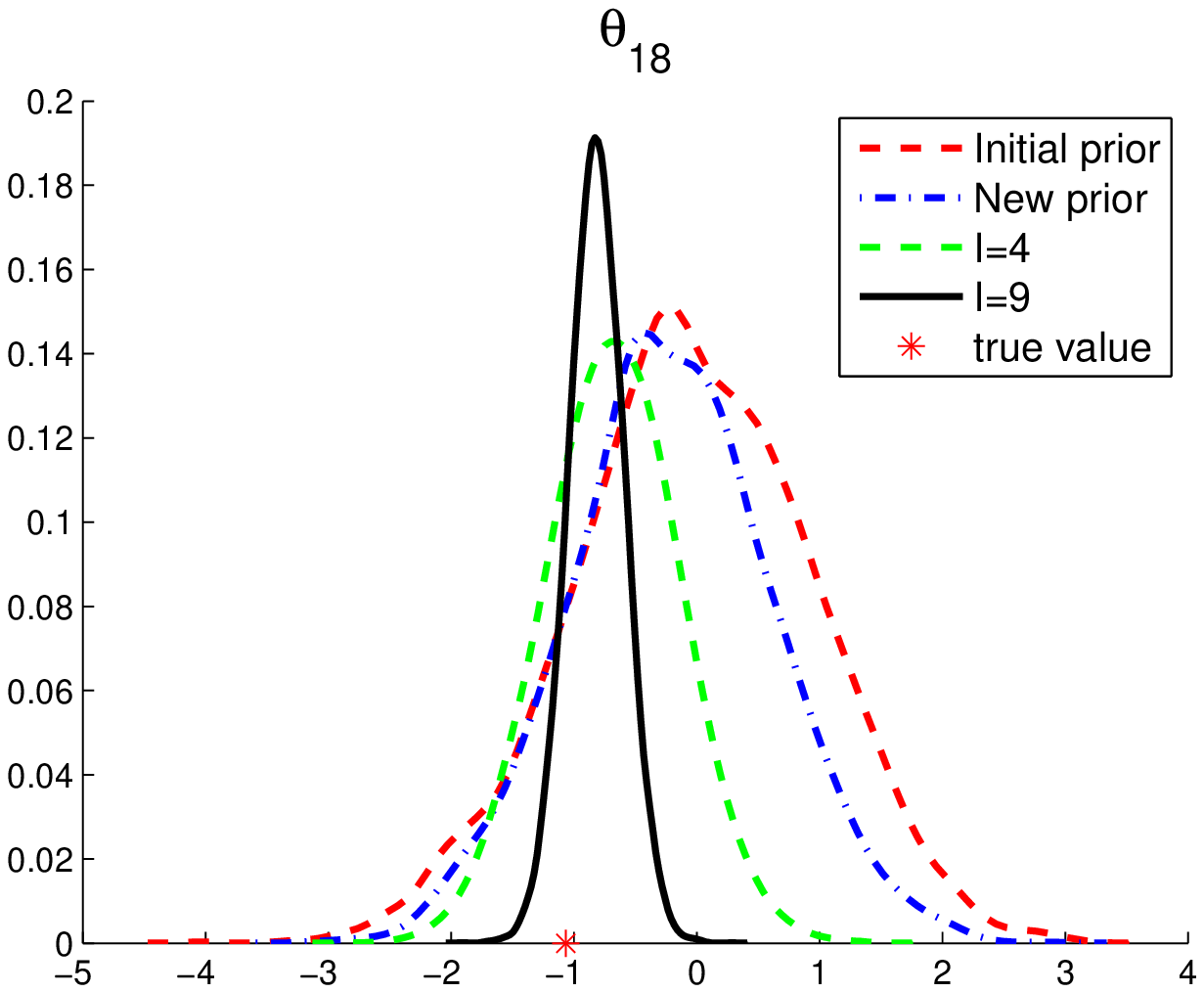}
  \caption{Marginal posterior density estimation of $\theta_1$, $\theta_5$, $\theta_9$, $\theta_{13}$, $\theta_{16}$ and $\theta_{18}$ in  different data assimilation steps.}\label{coeffdensity}
 \end{figure}

The credible interval and prediction interval, along with the true, and observation data, are illustrated for $u\big((x,0.5); t\big)$ and $u\big((0,5, y); t\big)$ in Figure \ref{coprx} and Figure \ref{copry}, respectively.  We find  that  the final credible intervals are tight and that the uncertainty associated with both the model fit and predictions decreases with respect to data assimilation step. The corresponding uncertainty decreases as $x$ and $y$ get closer to boundary, which is due to the deterministic Dirichlet boundary condition.
The marginal densities  for the unknown parameters in different assimilation steps are  plotted in Figure \ref{coeffdensity},  which shows the the support of density distribution becomes narrower as more data information is used in the posterior exploration.


\section{Conclusion}

We presented a two-stage EnKF  using coarse GMsFEM models in the paper. In the first stage, we used a very coarse GMsFEM model  to approximate the forward model, and the corresponding misfit-to-observed  problem is solved by using GMsFEM. We then constructed the new prior based on the ensembles by traditional EnKF, and then  used  GMsFEM and gPC to obtain a compact  representation for the model response based on the new prior. The two-stage EnKF  was employed to explore the surrogate posterior density, which was incorporated by the surrogate likelihood and the updated prior. It showed that the proposed method leads to the approximate posterior with better efficiency and accuracy  than traditional EnKF method.

The deterministic and statistical methods were combined  together to solve inverse problems.  We obtained not only the point estimate and  confidence interval but also the statistical properties of the unknowns.    A new prior was constructed for Bayesian inference  using  coarse GMsFEM models.  The new prior  contains the significant  region or support of the posterior and  is incorporated with the likelihood to be explored. For non-Gaussian models, we presented a two-stage EnKF using normal score transformation.

\begin{ack}
L. Jiang acknowledges the support of Chinese NSF 11471107.
\end{ack}

\end{document}